\documentclass[fleqn,10pt]{wlscirep}
\usepackage[utf8]{inputenc}
\usepackage[T1]{fontenc}
\usepackage{subcaption}

\newcommand{\R}{ \ensuremath{\mathbb{R}}}
\newcommand{\C}{ \ensuremath{\mathbb{C}}}
\newcommand{\Z}{ \ensuremath{\mathbb{Z}}}
\newcommand{\Id}{ \ensuremath{\mathbb{I}}}
\newcommand{\Han}{ \ensuremath{\mathbb{H}}}
\newcommand{\Vand}{ \ensuremath{\mathbb{V}}}
\newcommand{\hb}{ \ensuremath{\mathbf{h}}}
\newcommand{\Aop}{ \ensuremath{\mathbb{A}}}
\newcommand{\rv}{ \ensuremath{\mathbf{v}}}

\newcommand{\rez}{ \ensuremath{\mathbf{r}}}
\newcommand{\dimn}{ \ensuremath{{d}}}
\newcommand{\eb}{\ensuremath{\mathbf{{\overrightarrow{\mathbf{1}}}}}}
\newcommand{\WW}{\mathbf{W}}
\newcommand{\RV}{\mathbf{V}}
\newcommand{\mcV}{\mathcal{V}}

\newcommand{\X}{\mathbf{X}}
\newcommand{\x}{\mathbf{x}}
\newcommand{\bfr}{\mathbf{r}}
\newcommand{\bfu}{\mathbf{u}}

\newcommand{\bfv}{\mathbf{v}}
\newcommand{\Y}{\mathbf{Y}}
\newcommand{\y}{\mathbf{y}}
\newcommand{\z}{\mathbf{z}}
\newcommand{\bfg}{\mathbf{g}}
\newcommand{\bfh}{\mathbf{h}}
\newcommand{\bfe}{\mathbf{e}}
\newcommand{\bff}{\mathbf{f}}
\newcommand{\DD}{\mathbf{D}}

\newcommand{\bfV}{\mathbf{V}}
\newcommand{\bfW}{\mathbf{W}}
\newcommand{\bfS}{\mathbf{S}}
\newcommand{\bfSS}{\mathbf{Q}}

\newcommand{\bfss}{\mathbf{q}}
\newcommand{\snap}{\mathbf{S}}
\newcommand{\DOl}{\mathbf{y}}

\newcommand{\0}{\mathbf{0}}
\newcommand{\CDS}{\mathbf{F}}
\newcommand{\DDS}{\mathbf{T}}
\newcommand{\win}{\mathsf{w}}
\newcommand{\now}{\mathsf{p}}
\newcommand{\Koop}{\mathcal{U}}
\newcommand{\KoopM}{\mathbb{U}}
\newcommand{\Vanderm}{\mathbb{V}}
\newcommand{\Proj}{\mathbb{P}}
\newcommand{\bfalpha}{{\boldsymbol\alpha}}

\newcommand{\bfvarepsilon}{{\boldsymbol\varepsilon}}
\newcommand{\bfdelta}{{\boldsymbol\delta}}
\newcommand{\bfpsi}{{\boldsymbol\psi}}
\newcommand{\bfPsi}{{\boldsymbol\Psi}}

\newcommand{\bfvareps}{{\boldsymbol\varepsilon}}
\newcommand{\bfLambda}{{\boldsymbol\Lambda}}

\newcommand{\bfSigma}{{\boldsymbol\Sigma}}
\newcommand{\flow}{{\boldsymbol\varphi}}
\newcommand{\imunit}{\mathfrak{i}}

\newtheorem{theorem}{Theorem}[section]
\newtheorem{definition}[theorem]{Definition}

\usepackage{algorithm}
\usepackage{algpseudocode}

\usepackage[symbol]{footmisc}

\title{A Koopman Operator-Based Prediction Algorithm and its Application to COVID-19 Pandemic\footnote[2]{Distribution Statement A: Approved for Public Release, Distribution Unlimited.}}

\author[1,4,*]{Igor Mezi\'{c}}
\author[2]{Zlatko Drma\v{c}} 
\author[3]{Nelida \v{C}rnjari\'{c}-\v{Z}ic}
\author[3]{Senka Ma\'{c}e\v{s}i\'{c}}
\author[4]{Maria Fonoberova}
\author[4]{Ryan Mohr}
\author[1,4]{Allan M. Avila}
\author[5]{Iva Manojlovi\'{c}}
\author[4]{Aleksandr Andrej\v{c}uk}
\affil[1]{University of California, Santa Barbara, CA, 93106, USA}
\affil[2]{Faculty of Science,  University of Zagreb, Croatia}
\affil[3]{University of Rijeka, Croatia}
\affil[4]{AIMdyn Inc., Santa Barbara, CA, 93101, USA}
\affil[5]{Department of Applied Mathematics, Faculty of El. Engineering, Univ Zagreb, Croatia}

\affil[*]{mezic@ucsb.edu}

\begin{abstract}
The problem of prediction of behavior of dynamical systems has undergone a paradigm shift in the second half of the 20th century with the discovery of the possibility of chaotic dynamics in simple, physical, dynamical systems for which the laws of evolution do not change in time. The essence of the paradigm is the long term exponential divergence of trajectories. However, that paradigm does not account for another type of unpredictability: the ``Black Swan" event. It also does not account for the fact that short-term prediction is often possible even in systems with exponential divergence. In our framework, the Black Swan type dynamics occurs when an underlying dynamical system suddenly shifts between dynamics of different types. A learning and prediction system should be capable of recognizing the shift in behavior, exemplified by ``confidence loss".  In this paradigm, the predictive power is assessed dynamically and confidence level is used to switch between long term prediction and local-in-time prediction.   Here we explore the problem of prediction in systems that exhibit such behavior. The mathematical underpinnings of our theory and algorithms are based on an operator-theoretic approach in which the dynamics of the system are embedded into an infinite-dimensional space. 
The dynamical switching from global to local prediction algorithm enabled a successful prediction of influenza cases.
We show that the framework correctly identifies the  2009-2010 flu pandemic as a Black Swan event, that prevented machine learning-based algorithms 
from showing subsequent good performance.
The world has recently experienced a Black swan event that lead to the COVID-19 pandemic. We deployed our algorithm to assess its evolution. The results show that, despite being capable of capturing the dynamics of the observed cases of the disease locally, in states and counties, the prediction algorithm is robust to perturbations of the available data, induced for example by delays in reporting or sudden increase in cases due to increase in testing capability. This is achieved in an entirely data-driven fashion, with no underlying mathematical model of the disease.
We  discuss the prediction problem in other complex dynamics datasets, such as signature indices of geomagnetic substorms. In addition, fundamental limits on predictability that our theory implies are discussed. 
\end{abstract}

\begin{document}

\flushbottom
\maketitle
%
%
\thispagestyle{empty}


\section{Introduction}


Ability for  prediction of
events is one of the key differentiators of homo sapiens. The key element of prediction is reliance on collected data over some time interval for estimation of evolution over the next time period. Mathematicians have long worked on formal aspects of prediction theory, and separate streaks such as the Wiener-Kolmogorov, \cite{doob1953stochastic}, Furstenberg \cite{furstenberg1960stationary} and Bayesian prediction \cite{pole2018applied} have emerged.
However, all of these are concerned with prediction of future events based on a, typically long, sequence of prior observations. This is rooted in assumptions on statistical stationarity of the underlying stochastic process. 

The point of view on prediction in this paper is quite different: we view the process over a short (local) time scale and extract its coarse-grained ingredients.  We proceed with prediction of the evo\-lution based on these, learning the process and building a global time-scale on which such prediction is valid. Then, we monitor for~ the change in such coarse-grained ingredients, detect if a substantial change (a \emph{``Black Swan"} event \cite{Taleb-2007}, see \textbf{Supplementary Information} section S1 for the mathematical definition that we use.) is happening, and switch back to local learning and prediction. In this way, we accept the limitations on predictability due to, possibly finite time, nonstationarity, and incorporate them into the prediction strategy. In principle, such strategy is valuable even in the case where over a long time interval the system is indeed a stationary stochastic process with respect to some invariant measure. An example is the Lorenz system, the prototypical system exhibiting chaotic dynamics and the Butterfly Effect \cite{lorenz1963deterministic,luzzatto2005lorenz}, studied in \textbf{Supplementary Information} section S2.1.
For  typical learning algorithms Black Swan events are devastating: the learning algorithm has to be restarted as otherwise it would learn the deviation as normal. Our  method  includes a technique of retouching the Black Swan event data, wherein its replacement with a realization of the normal process evolution obviates the need for a restart.

However, the case of Black Swan unpredictability and Butterfly Effect unpredictability are ontologically different. Namely, the Butterfly Effect is the consequence of dynamics inherent to the system, while the Black Swan arises from an action external to system dynamics.

We approach the prediction problem from the  perspective of Koopman operator theory \cite{koopman1931hamiltonian,mezic:2005,hua2017high,giannakis2020extraction,korda2018linear,khodkar2018data} in its recently developed form that is applicable to nonstationary stochastic processes \cite{vcrnjaric2017koopman,mezic2019spectrum}. The Koopman operator theory is predicated on existence of the composition operator that dynamically evolves all the possible observables on the data, enabling the study of nonlinear dynamics by examining its action on a linear space of observables. The key ingredients of this approach become eigenvalues and eigenfunctions of the Koopman operator and the associated Koopman Mode Decomposition (KMD) of the observable functions, which is then approximated numerically using Dynamic Mode Decomposition (DMD). The numerical approach used in this work relies on lifting the available data to higher dimensional space using Hankel-Takens matrix and on the improved implementation of DMD algorithm for discovering the approximations of the Koopman modes with small residuals. The obtained Koopman mode approximations and the related eigenvalues, called Ritz pairs, are crucial for obtaining satisfactory predictions using KMD. One of the main advantages of the method is that it completely data-driven, i.e., model-free.

%

The retouching method is presented in \S \ref{S=Influenza}.  In \S \ref{S=COVID-19} we present the application of the prediction algorithm to infection cases data of the current COVID-19 pandemic.

In \textbf{Supplementary Information}, we provide mathematical background of all algorithms, with worked examples    that clarify all technical details and offer more in depth discussions. 
 Further,  we test the predictive potential of the proposed method using diverse case studies with the goal to stress the data driven (model free) nature of the method.  For instance, an application of the method to  prediction of physiological process data, that additionally validates our~ methodology, could be of interest in personalized medical treatment context. As further  case studies,  we use the Lorenz system and the \textsf{AL} index of geomagnetic substorms.

\section{Results}

\subsection{Prediction with \textsf{Koopman Mode Decomposition }}\label{SS=KMD-intro}
Our approach starts with the Koopman operator family $\Koop^t$, which acts on observables $f$ by composition 
    $\Koop^t f (\x) = f(\x( t))$. 
$\Koop^t$ is a linear operator that allows studying the nonlinear dynamics by examining its action on a linear space $\mathcal{F}$ of observables. In a data driven setting, usually the states $\z_i\approx \x(t_i)$ of the dynamical system at discrete time moments $t_i$ are known. They are governed by the
discrete dynamical system $\z_{i+1}=\DDS(\z_i)$, for which the Koopman operator reads $\Koop f = f\circ\DDS$. \\
\indent The key of the spectral analysis of the dynamical system is a representation of a vector valued  observable $\bff=(f_1,\ldots,f_{d})^T$ as a linear combination of the eigenfunctions $\bfpsi_j$ of $\Koop$. Under certain assumptions,  each observable $f_i$ can be approximated as
$f_i(\z) \approx \sum_{j=1}^{m}\bfpsi_j(\z) (\mathbf{v}_j)_i$, where $\bfpsi_1, \ldots, \bfpsi_m$ are selected eigenfunctions and
$\mathbf{v}_j =\left((\mathbf{v}_j)_1 \ldots (\mathbf{v}_j)_{d}\right)^T$  are the vectors of the coefficients. 
\noindent Then, we can predict the values of the observable $\bff$ at the \emph{future} states $\DDS(\z)$, $\DDS^2(\z), \ldots$ by numerically evaluating
\begin{equation}\label{eq:predict-intro}
 (\Koop^k\bff)(\z) \stackrel{\mathrm{\tiny def}}{=} \bff(\DDS^k(\z))\approx  \sum_{j=1}^{m}  \lambda_j^k \bfpsi_j(\z) \mathbf{v}_j , \;\;k=1,2,\ldots
\end{equation}
The decomposition of the observables is called the \emph{Koopman Mode Decomposition} (\textsf{KMD}); the scalars $\lambda_j$ are the Koopman eigenvalues, and the $\mathbf{v}_j$'s are the Koopman modes. Their numerical approximations can be computed based on the supplied  data pairs $(\bff(\z_i),\bff(\DDS(\z_i)))$, $i=0,\ldots, M$, using e.g. the \emph{Dynamic Mode Decomposition} (\textsf{DMD}) \cite{schmid:2010}, \cite{DDMD-RRR-2018}, \cite{DMM-2019-DDKSVC-DFT}.
\noindent In this work, the numerical algorithm is fed with the sequence of successive data snapshots, called active window, to compute the approximation of KMD, and then (\ref{eq:predict-intro}) is used  for prediction. If we use wider windows 
and longer forecasting lead time, we speak of \emph{global prediction}; for narrower windows with locally adapted widths and shorter forecasting lead time we have \emph{local prediction (nowcasting)}.

\subsection{Case study: Influenza epidemics.}\label{S=Influenza}
As first example for showing our prediction methodology, we use the set of data associated with influenza epidemics. Clearly, not driven by an underlying deterministic dynamical system, the influenza time series exhibits substantial regularity in that it occurs typically during the winter months, thus enabling coarse-grained prediction of the type ``we will see a very small number of cases of influenza occurring in summer months". However, predicting the number of influenza cases accurately is a notoriously hard problem \cite{lazer2014parable}, exacerbated by the possibility that a vaccine designed in a particular year does not effectively protect against infection. {Moreover, the H1N1 pandemic that occurred in 2009 is an example of a Black Swan event.}

The World Health Organiza\-ti\-on's \textsf{FluNet} is a global web-based tool for influenza virological surveillance. \textsf{FluNet} makes publicly available data on the number of specimens with the detected influenza viruses of type \textsf{A} and type \textsf{B}. The data have been collected from different countries, starting with the year 1997, and are updated weekly by the National Influenza Centers (\textsf{NICs}) of the Global Influenza Surveillance and Response System (\textsf{GISRS}) and other national influenza reference laboratories, collaborating actively with \textsf{GISRS}. We use the weekly reported data for different countries, which consist of  the number of received specimens in the laboratories,  the distribution of the number of specimens with confirmed viruses of type~ \textsf{A}.

The Koopman Mode Decomposition was used in the context of analyzing  the dynamics of the flu epidemic from different  - Google Flu - data in \cite{proctor}. We remark that the authors of that paper have not attempted prediction, and have analyzed only ``stationary" modes - e.g. the yearly cycles, thus making the paper's goals quite different from the  nonstationary prediction pursued here.

\begin{figure}[ht]
\begin{subfigure}{\textwidth}
	\includegraphics[width=\linewidth,height=2.1in]
	{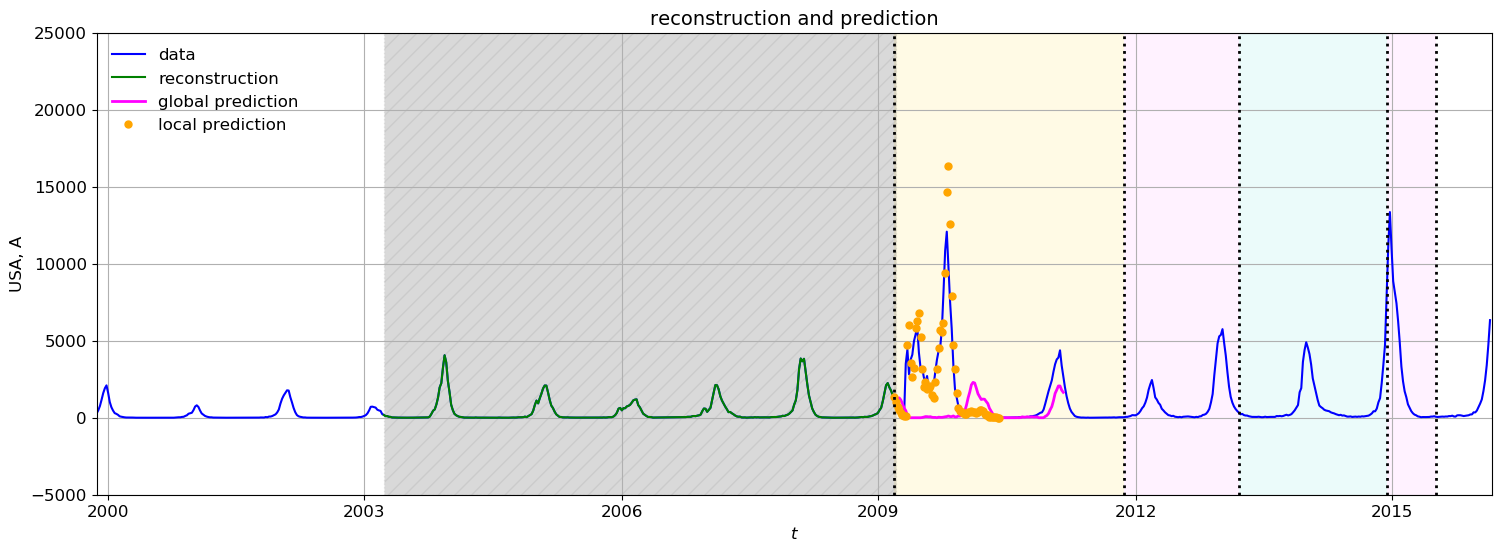}
	\caption{Local prediction}
	\label{fig:influenza-US-global}
	\end{subfigure}
\begin{subfigure}{\textwidth}
	\includegraphics[width=\linewidth, height=2.1in]
	{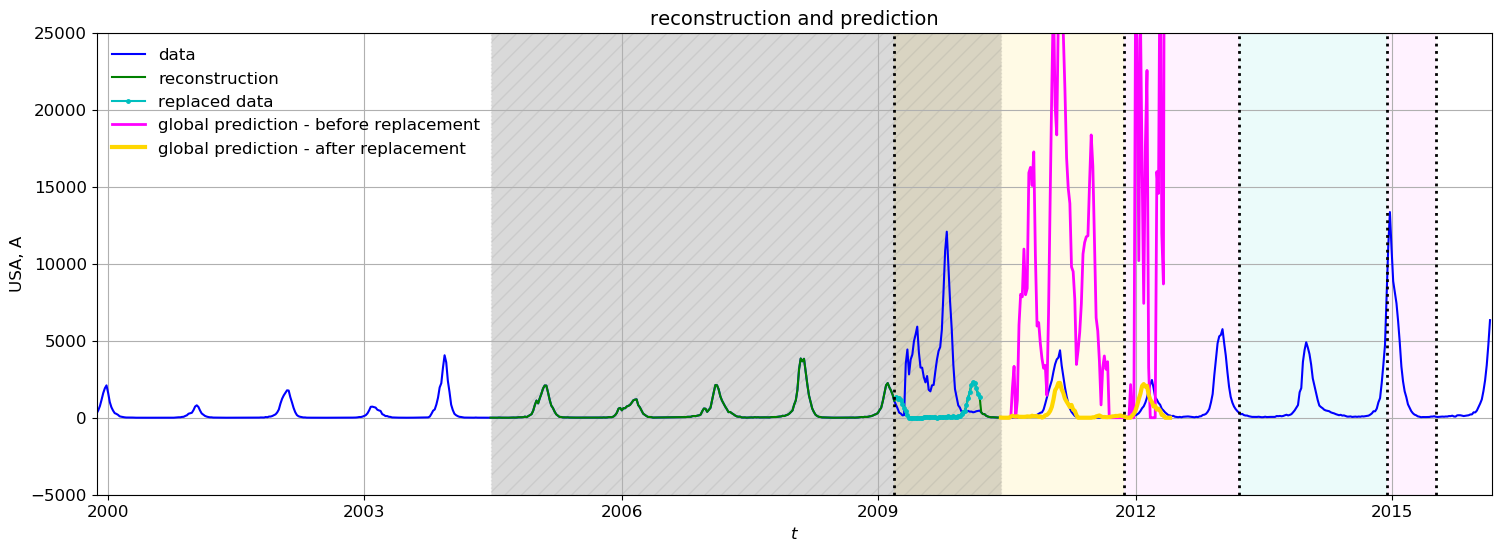}
	\caption{Global prediction with retouching the Black Swan event data}
	\label{fig:influenza-US-replacement}
	\end{subfigure}
    \caption{
    Influenza data (USA). (\ref{fig:influenza-US-global}): The data are collected in the  window April 2003 -- April 2009 (shadowed rectangle) and then the dynamics is predicted for $104$ weeks ahead. 
    The local prediction algorithm recovers the prediction capability by forgetting the old data and using narrower learning windows. {The local prediction algorithm delivers prediction for one week ahead.}
    (\ref{fig:influenza-US-replacement}): The active window (shadowed rectangle) is July 2004 -- July 2010, and  the dynamics is predicted for  $104$ weeks ahead. The global prediction fails due to the Black Swan data in the learning window. ({Some predicted values were even negative; those were replaced with zeros.})  The global prediction algorithm recovers after the retouching the Black Swan event data, which allows for using big learning window. Compare with positions of the corresponding colored rectangles in Figure \ref{fig:influenza4}.  
    }	
\end{figure}

We first compare  the global and the local prediction algorithms. The \textsf{KMD} is computed using active windows of size $\win = 312$, and the $208 \times 104$ Hankel-Takens matrices.
In Fig.\ref{fig:influenza-US-global}, we show the performances of both algorithms, using the learning data from the window April 2003 -- April 2009 (shadowed rectangle). In the global prediction algorithm the dynamics is predicted for $104$ weeks ahead. The first type of failure in the global prediction algorithm and forecasting appears after the Black Swan event occurred in the years 2009 and 2010. This is recognized by the algorithm, so that it adapts by using the smallest learning span and, with this strategy, it allows for reasonably accurate forecasting, at least for shorter lead times.   This data, in addition to those from \textbf{Supplementary Information} section S2.4 show the benefits of monitoring the prediction error and switching to local prediction. The initial Hankel-Takens matrix is $3\times 2$, and the threshold for the local prediction  relative error in \textbf{Supplementary Information} Algorithm S4 is $0.005$.

\subsubsection*{Retouching the Black Swan event data}\label{SS=Retouching} Next, we  introduce an approach that robustifies the global algorithm in the presence of disturbances in the data, including the missing data scenario. We use the data window July 2004--July 2010, which contains a Black Swan event in the period 2009--2010.  As shown in Figure \ref{fig:influenza-US-replacement}, the learned \textsf{KMD} failed to predict the future following the active training window. This is expected because the perturbation caused by the Black Swan event resulted in the computed Ritz pairs that deviated from the precedent ones (from a learning window before disturbance), and, moreover, with most of them having large residuals. This can be seen as a second type of failure in the global prediction. 

The proposed Black Swan event detecting device, built in the prediction algorithm (see \textbf{Supplementary Information} Algorithm S3), checks for this anomalous behaviour of the Ritz values and pinpoints the problematic subinterval.  Then, the algorithm replaces the corresponding supplied data with the values obtained as predictions based on the time interval preceding the Black Swan event.  Figure \ref{fig:influenza-US-replacement} shows that such a retouching of the disturbance allows for a reasonable global prediction.

Note that in  a realistic situation, global predictions of this kind will trigger response from authorities and therefore prevent its own accuracy and induce loss of confidence, whereas local prediction mechanisms need to be deployed again. 

%


\begin{figure} [ht]
	\centering
	\includegraphics[width=0.8\linewidth,height=4in]
	{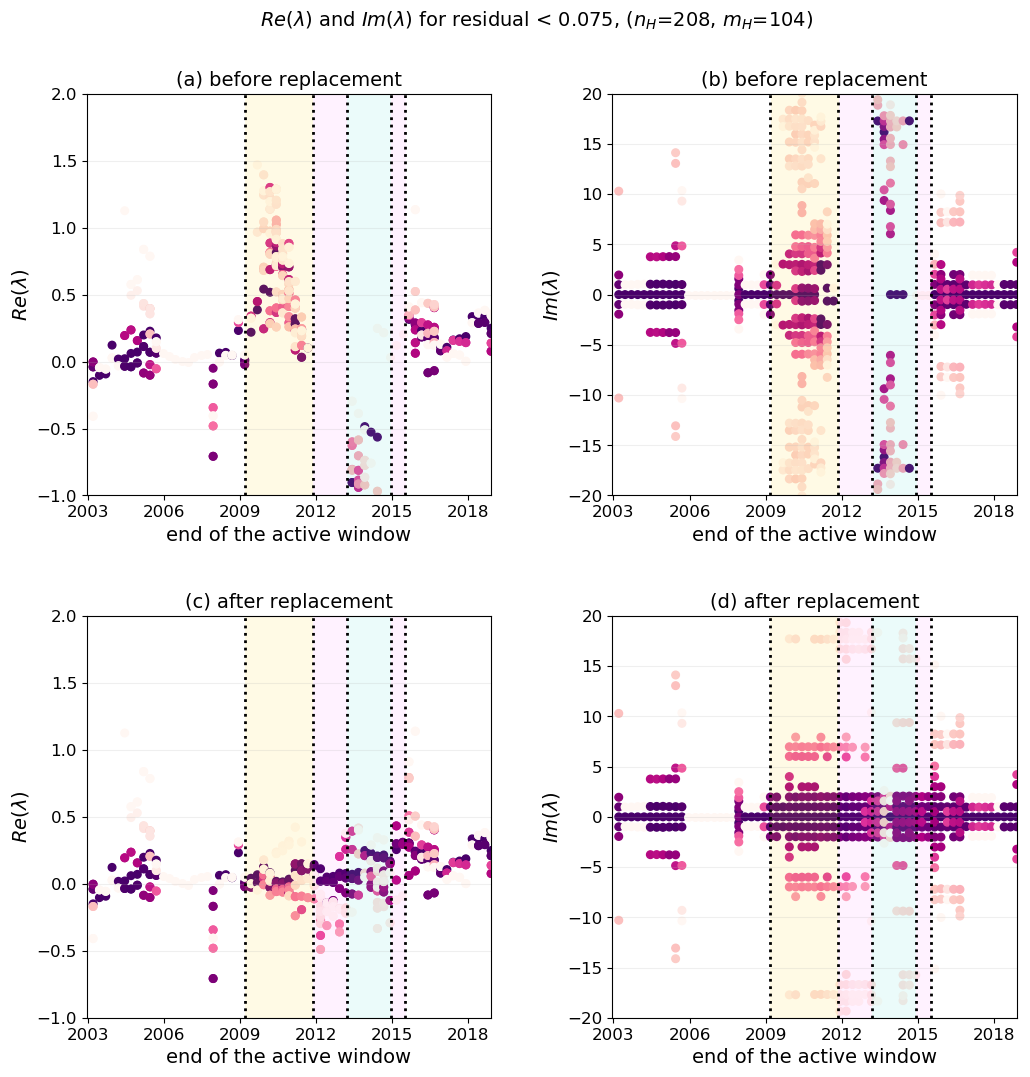}
	\caption{
	The real and imaginary parts of Ritz values with residuals bellow $\eta_r=0.075$ for sliding active windows. 
	The color intensity of eigenvalues indicates the amplitudes of the corresponding modes. Pink rectangles mark ends of training windows with no acceptable Ritz values.
	Note how the unstable eigenvalues ($\Re(\lambda)>0$) impact the  prediction performance, and how the retouching moves them towards neutral/stable -- this is shown in the {yellow rectangle in} panels (a) and (c). Also influenced by the disturbance are the {eigenvalues in the light blue rectangles in panels (a), (b);  retouching moves the real parts of eigenvalues towards neutral/stable and rearranges them in a lattice-like structure \cite{Mezic-2019-Spectrum-JNS}, as shown in panels (c), (d).} Compare with Figure \ref{fig:influenza-US-replacement}.
	}
	\label{fig:influenza4}
\end{figure}

\subsubsection*{Monitoring and restoring the Ritz values}
We now discuss the effect of  the Black Swan event and its retouching to the computed eigenvalues and eigenvectors.
We have observed that, as soon as a disturbance starts entering the training windows, the Ritz values start exhibiting atypical behavior, e.g. moving deeper into the right half plane (i.e. becoming more unstable), and having larger residuals because the training data no longer represent the Krylov sequence of the underlying Koopman operator.\\   
\indent This is illustrated in the panels (a) and (b) in Figure \ref{fig:influenza4}, which show, for the sliding training windows, the real and the imaginary parts of those eigenvalues for which the residuals of the associated eigenvectors are smaller than $\eta_r=0.075$.  Note the absence of such eigenvalues in time intervals that contain the disturbance caused by the Black Swan event. \\
\indent On the other hand, the retouching technique that repairs the distorted training data restores the intrinsic dynamics over the entire training window. The distribution of the relevant eigenvalues becomes more consistent, and the prediction error decreases, see panels (c) and (d) in Figure \ref{fig:influenza4}, and
in \textbf{Supplementary Information} Figure S16.

\subsubsection*{Discussion} Our proposed retouching procedure
relies on detecting anomalous behavior of the Ritz values; a simple strategy of monitoring the spectral radius of active windows (absolutely largest Ritz value extracted from the data in that window) is outlined in \textbf{Supplementary Information}. Note that this can also be used as a \emph{litmus test} for switching to the local  prediction algorithm.  In \textbf{Supplementary Information}, we provide further examples, with the influenza data, that confirm the usefulness of the retouching procedure. In general, this procedure can also be adapted to the situation when the algorithm receives a signal that the incoming data is missing or corrupted.

\subsection{COVID-19 prediction.}\label{S=COVID-19}

The second set of data we consider is that associated with the ongoing COVID-19 pandemic. Because the virus is new, the whole event is, in a sense, a ``Black Swan". 
However, as we show below, the prediction approach advanced here is capable of adjusting quickly to the new incoming, potentially sparse data and is robust to inaccurate reporting of cases.

At the beginning of the spread of COVID-19, we have witnessed at moments rather chaotic situation in gaining the knowledge on the new virus and the disease. The development of COVID-19 diagnostic tests made tracking and modelling feasible, but with many caveats: the data itself is clearly not ideal, as it depends on the reliability of the tests, testing policies in different countries (triage, number of tests, reporting intervals,   
reduced testing during the weekends), contact tracing strategies,  using surveillance technology, credit card usage and phone contacts tracking, the number of asymptomatic transmissions  etc. 
Many different and unpredictable exogenous factors can distort it. So, for instance the authors of \cite{COVID-data-countries}
comment at \url{https://ourworldindata.org/coronavirus-testing} that e.g. "The Netherlands, for instance, makes it clear that not all labs were included in national estimates from the start. As new labs get included, their past cumulative total gets added to the day they begin reporting, creating spikes in the time series." For a prediction algorithm, this creates a Black Swan event that may severely impair prediction skills, see \S \ref{SS=Retouching}. \\
\indent 
This poses challenging problems to the compartmental type models of  (SIR, SEIR) which in order to be useful in practice have to be coupled with data assimilation to keep adjusting the key parameters, see e.g. \cite{Nadler-SIR-model-DA}. Our technique of retouching (\S \ref{SS=Retouching}) can in fact be used to assist data assimilation by detecting Black Swan disturbance and thus to avoid assimilating disturbance as normal.\\
%
\indent In the \textsf{KMD} based framework, the changes in the dynamics are automatically assimilated on-the-fly by recomputing the \textsf{KMD} using new (larger or shifted) data snapshot windows. This is different from the compartmental type models of infectious diseases, most notably in the fact that the procedure presented here does not assume any model and, moreover, that it is entirely oblivious to the nature of the underlying process. 

\subsection*{An example: European countries}
\begin{figure} [ht]
	\includegraphics[width=0.33\linewidth,height=2in]{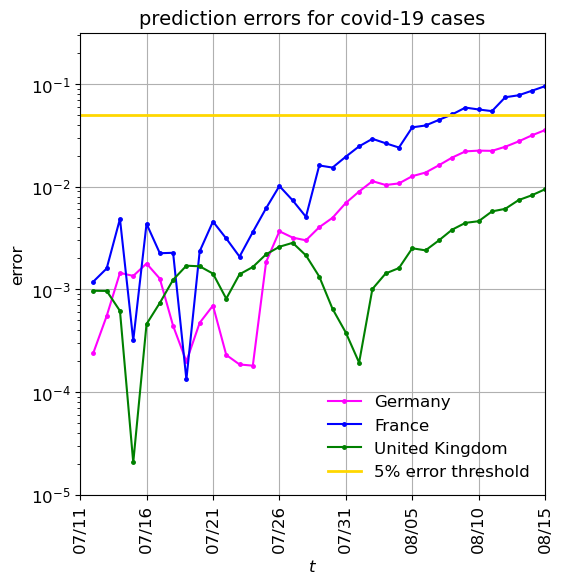}
	\includegraphics[width=0.33\linewidth,height=2in]{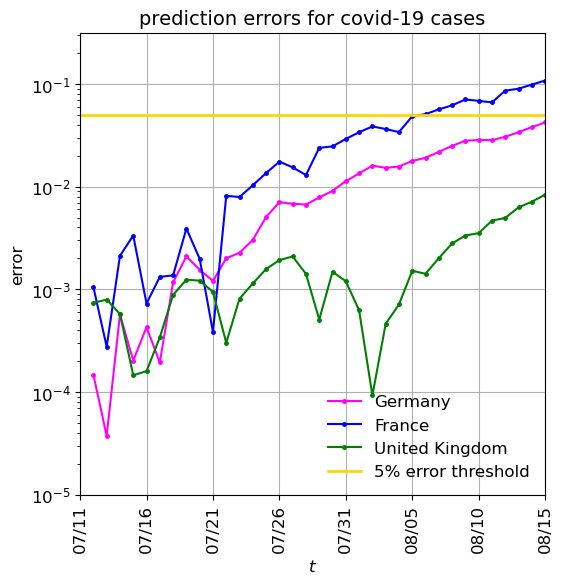}
\includegraphics[width=0.33\linewidth,height=2in]{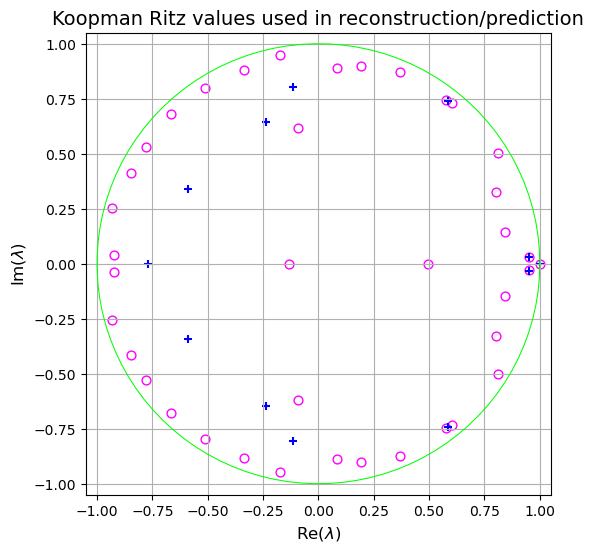}	
	\hfill
	\caption{
Prediction of COVID-19 cases ($35$ days ahead, starting July 11) for Germany, France and United Kingdom. \emph{Left panel}: The Hankel-Takens matrix $\Han$ is $282 \times 172$,  the learning data consists of $\hb_{1:40}$. The \textsf{KMD} uses $39$ modes. \emph{Middle panel}: 
The matrix $\Han$ is $363 \times 145$, the learning data is $\hb_{1:13}$. The \textsf{KMD} uses $12$ modes. \emph{Right panel}: The Koopman-Ritz values corresponding to the first (magenta circles)	and the middle (blue plusses) panel. Note how the three rightmost values nearly match.}
	\label{fig:covid-de-fr-uk-1}
\end{figure}

\begin{figure} [ht]
	\includegraphics[width=0.33\linewidth,height=2in]{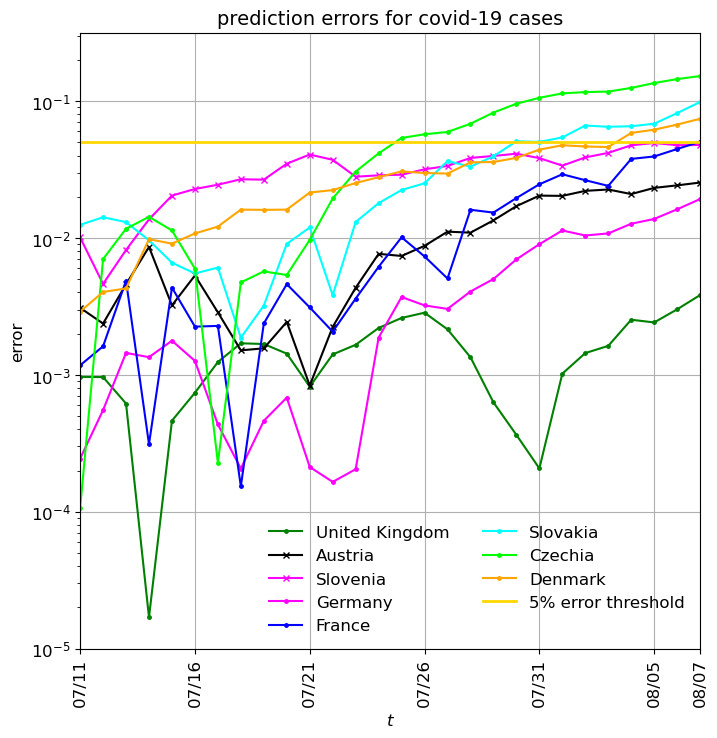}
	\includegraphics[width=0.33\linewidth,height=2in]{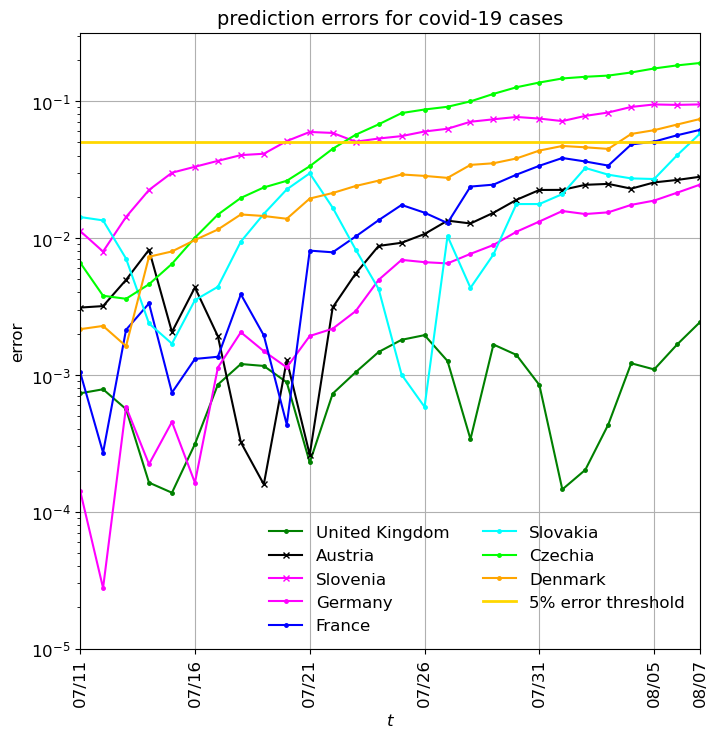}
    \includegraphics[width=0.33\linewidth,height=2in]{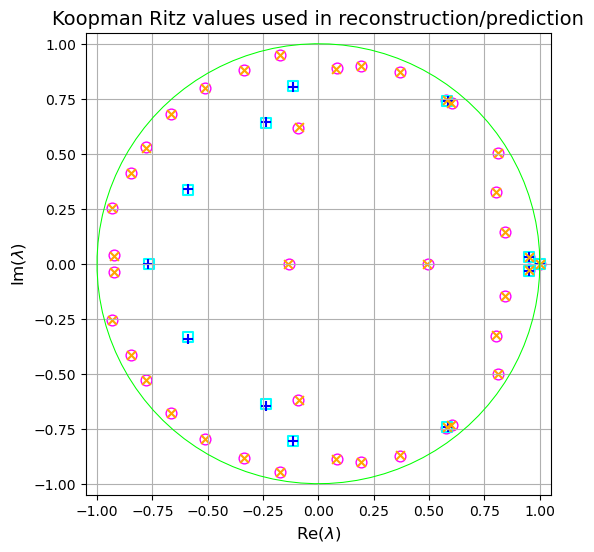}	
	\hfill
	\caption{
Prediction errors and KMD spectrum of COVID-19 cases ($28$ days ahead, starting July 11) for Germany, France, United Kingdom, Denmark, Slovenia, Czechia, Slovakia and Austria. \emph{Left panel}: The Hankel-Takens matrix $\Han$ is $752 \times 172$,  the learning data consists of $\hb_{1:40}$. The \textsf{KMD} uses $39$ modes. \emph{Middle panel}: 
The matrix $\Han$ is $968 \times 145$, the learning data is $\hb_{1:13}$. The \textsf{KMD} uses $12$ modes. \emph{Right panel}: The Koopman-Ritz values corresponding to the first two computations in Figure \ref{fig:covid-de-fr-uk-1} (magenta circles and blue pluses, respectively) and the the first two panels in this Figure (orange x-es and cyan squares, respectively). Note how the corresponding Koopman-Ritz values nearly match for all cases considered.}
	\label{fig:covid-de-fr-uk-five-more-1}
\end{figure}

\begin{figure}[ht]
	\centering
\includegraphics[width=0.32\linewidth]{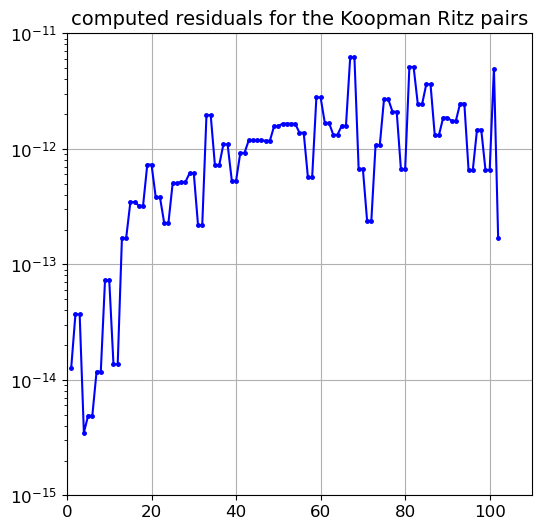}	
\includegraphics[width=0.32\linewidth]{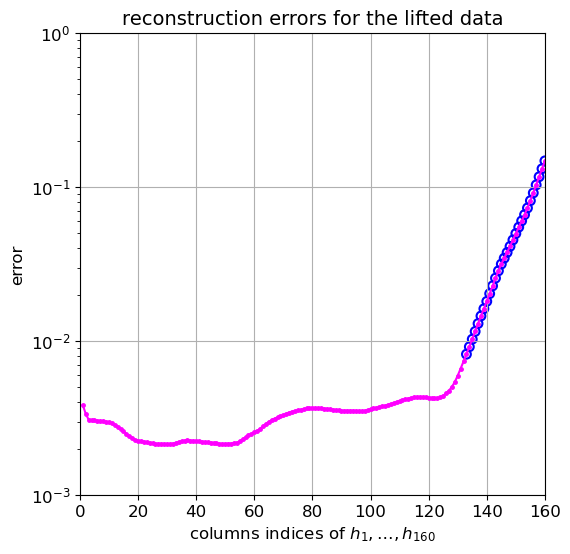}
\includegraphics[width=0.32\linewidth]{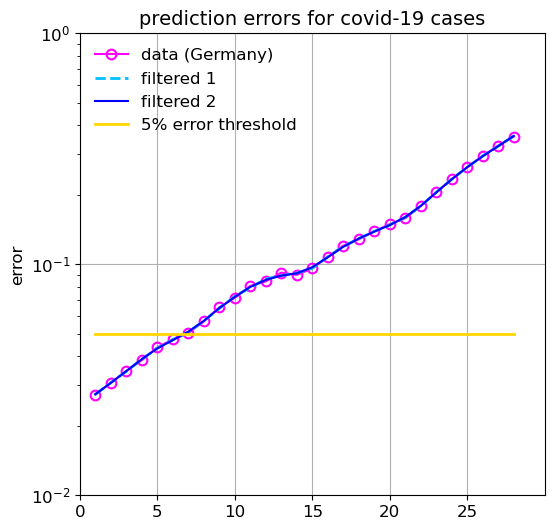}	
\caption{Prediction experiment with data from Germany. \emph{Left panel}: the computed residuals for the computed $102$ Koopman Ritz pairs (extracted from a subspace spanned by $132$ snapshots $\hb_{1:132}$). Note that all residuals are small. The corresponding Ritz values are shown in the first panel in Figure \ref {fig:main:DE_132_28_W}. \emph{Middle panel}: \textsf{KMD} reconstruction error for $\hb_{1:132}$ and the error in the predicted values $\hb_{133:160}$ (encircled with $\textcolor{blue}{\circ}$). The reconstruction is based on the coefficients $(\alpha_j)_{j=1}^r=\mathrm{arg\min}_{\alpha_j}\sum_{k} \| \hb_k - \sum_{j=1}^{r} \lambda_j^{k}\alpha_j \rv_j\|_2^2$.
\emph{Right panel}: Prediction errors for the period October 11 -- November 7. 
}
\label{fig:main:DE_132_28_UW}
\end{figure}

\begin{figure}[ht]
	\centering
\includegraphics[width=0.32\linewidth]{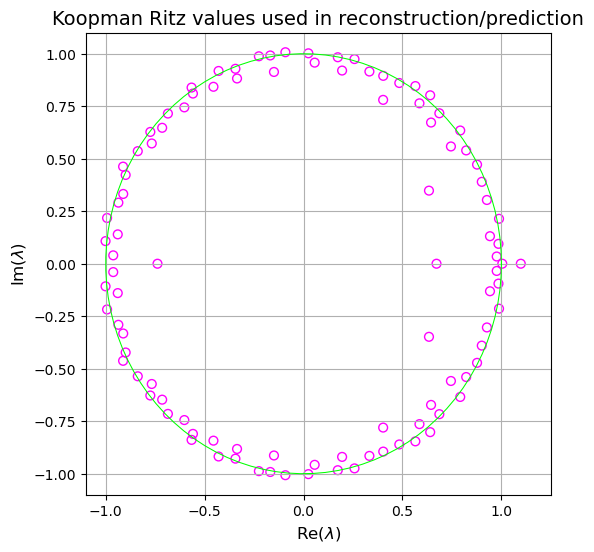}	
\includegraphics[width=0.32\linewidth]{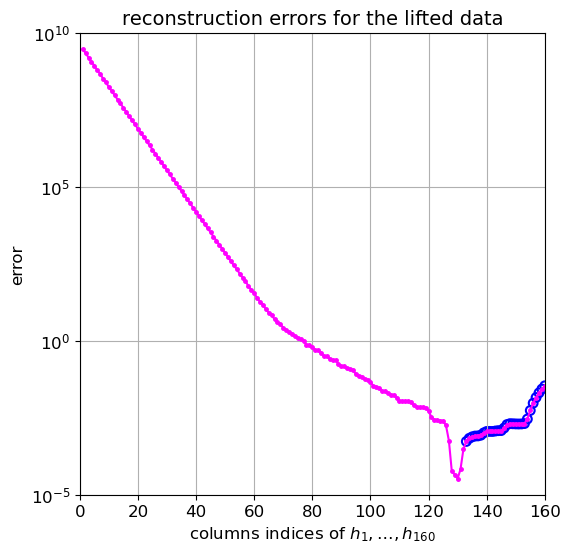}
\includegraphics[width=0.32\linewidth]{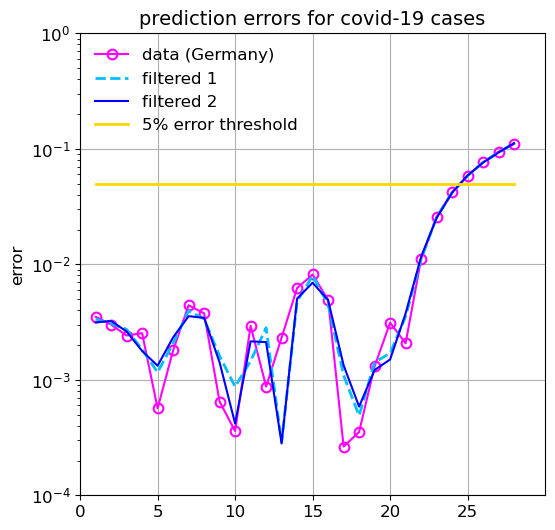}	
\caption{Prediction experiment with \textbf{DS3} with data from Germany. \emph{Left panel}: the computed $102$ Koopman Ritz values (extracted from a subspace spanned by $132$ snapshots $\hb_{1:132}$). The corresponding residuals are shown in the first panel in Figure \ref{fig:main:DE_132_28_UW}. \emph{Middle panel}: \textsf{KMD} reconstruction error for $\hb_{1:132}$ and the error in the predicted values $\hb_{133:160}$ (encircled with $\textcolor{blue}{\circ}$). The reconstruction is based on the coefficients $(\alpha_j)_{j=1}^r=\mathrm{arg\min}_{\alpha_j}\sum_{k} w_k^2 \| \hb_k - \sum_{j=1}^{r} \lambda_j^{k}\alpha_j \rv_j\|_2^2$.
\emph{Right panel}: Prediction errors for the period October 11 -- November 7. Compare with  the third graph in Figure \ref{fig:main:DE_132_28_UW}.
}
\label{fig:main:DE_132_28_W}
\end{figure}

As a first numerical example, we use the reported cumulative daily cases in European countries. In \textbf{Supplementary Information} section S1.5, we use this data for a detailed worked example that shows all technical details of the method. 
This is a good test case for the method -- using the data from different countries in the same vector observable poses an additional difficulty for a data driven revealing of the dynamics, because the countries independently and in an uncoordinated manner impose different restrictions, thus changing the dynamics on local levels.
For instance, at the time of writing these lines, a new and seemingly more infectious strain of the virus circulating in some parts of London and in south of England prompted the UK government to impose full lockdown measures in some parts of the United Kingdom. Many European countries reacted sharply and immediately suspended the air traffic with the UK.\\
\indent In the first numerical experiment, we use two datasets from the time period February 29 to November 19. and consider separately two sets of countries: Germany, France and the UK in the first, and Germany, France, UK, Denmark, Slovenia, Czechia, Slovakia and Austria in the second. The results for a particular prediction interval are given in Figure \ref{fig:covid-de-fr-uk-1} and Figure \ref{fig:covid-de-fr-uk-five-more-1}. For more examples and discussion how the prediction accuracy depends on the Government Response Stringency Index (GRSI \cite{GRSI, GRSI-2}) see \textbf{Supplementary Information} section S1.5.
 In the above examples, the number of the computed modes was equal to the dimension of the subspace of spanned by the training snapshots, so that the \textsf{KMD} of the snapshots themselves was accurate up to the errors of the finite precision arithmetic. In general, that will not be the case, and the computed modes will span only a portion the training subspace, meaning that the \textsf{KMD} of the snapshots might have larger representation error. (Here we refer the reader to \textbf{Supplementary Information} section S1.3, where all technical details are given.) This fact has a negative impact to the extrapolation forward in time and the problem can be mitigated by giving more importance to reconstruction of more recent weights. This is illustrated in Figures \ref{fig:main:DE_132_28_UW} and \ref{fig:main:DE_132_28_W}, where the observables are the raw data (reported cases) for Germany, extended by a two additional sequence of filtered (smoothened) values.

The figures illustrate an important point in prediction methodology, that we emphasized in the introduction: a longer dataset and a better data reconstruction ability (i.e. interpolation) does not necessarily lead to better prediction. Namely, weighting more recent data more heavily produces better prediction results. This was already observed in \cite{avila2020data} for the case of traffic dynamics, and the method we present here can be used to optimize the prediction ability.
\subsubsection*{An example: USA and worldwide data}
\begin{figure}[ht]
\begin{subfigure}{\textwidth}
	\centering
	\includegraphics[width=0.98\textwidth,height=2.4 in]{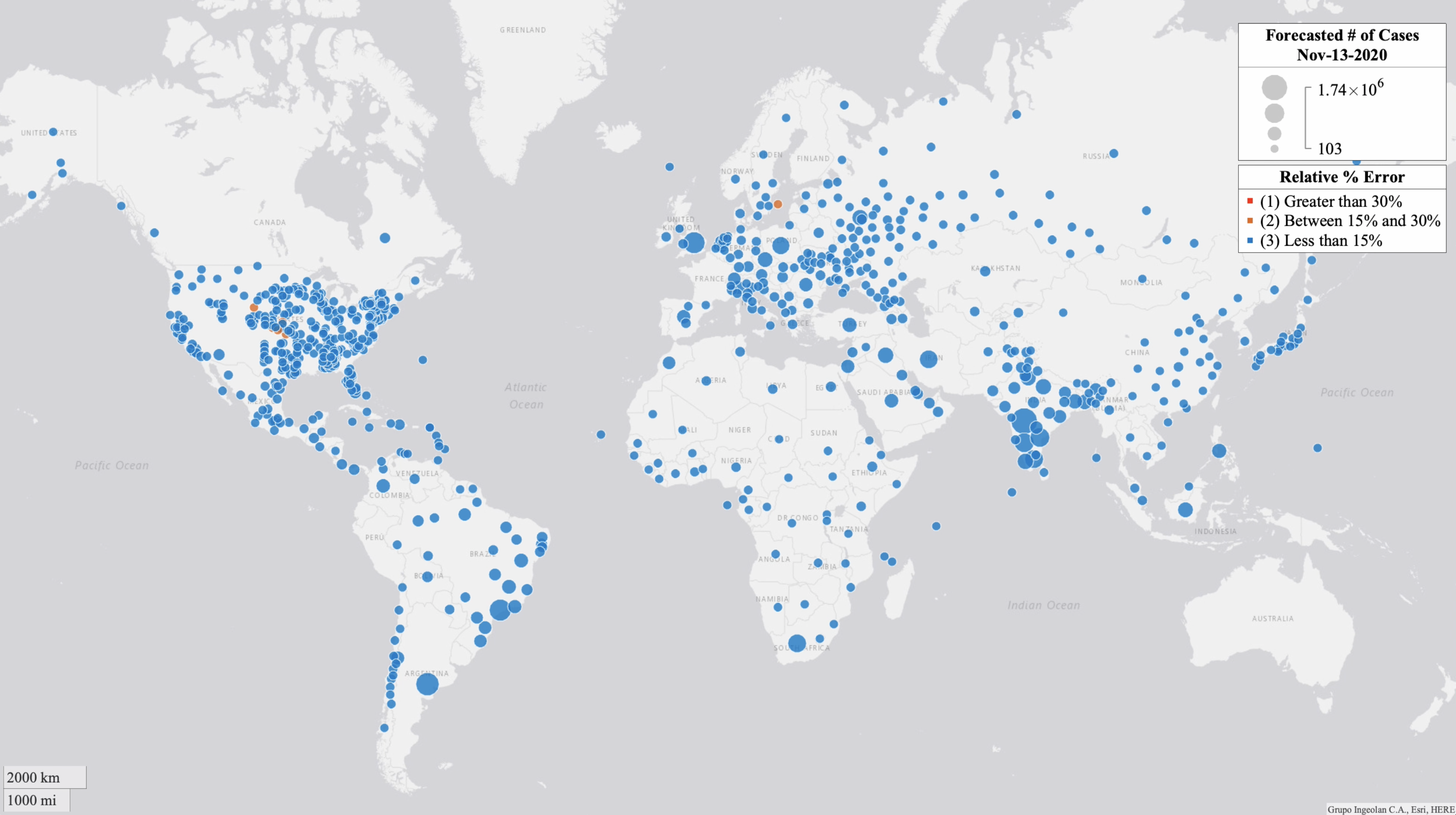}	
	\caption{Worldwide predicted cases and prediction error for COVID-19 pandemic on November 13, 2020.}
	\label{worlda}
\end{subfigure}
\begin{subfigure}{.49\textwidth}
\includegraphics[width=\linewidth,height=1.6in]{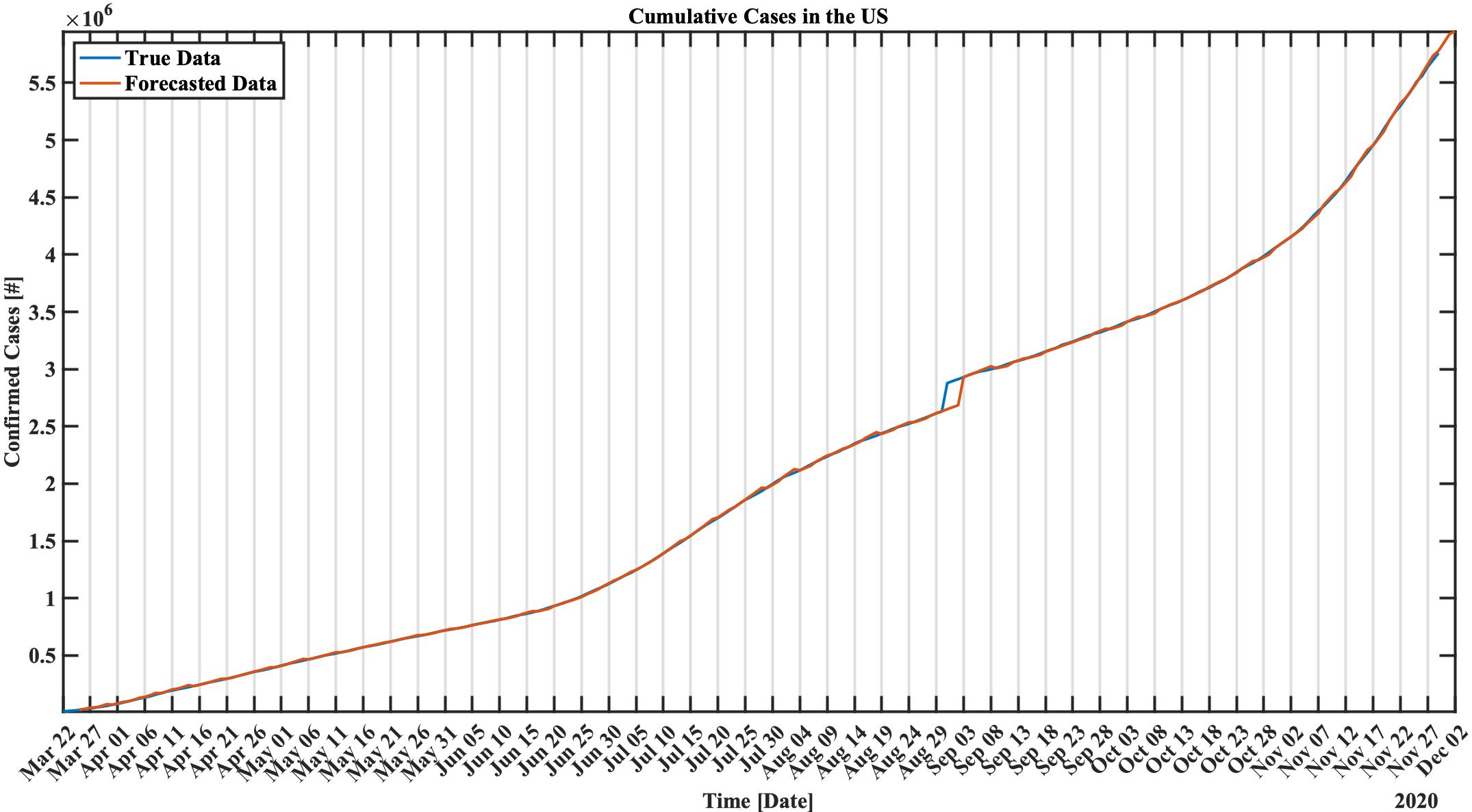}
\caption{Data and prediction for US number of COVID-19 cases.}
\label{us-pred}
\end{subfigure}
\begin{subfigure}{.49\textwidth}
\includegraphics[width=\linewidth, height=1.6in]{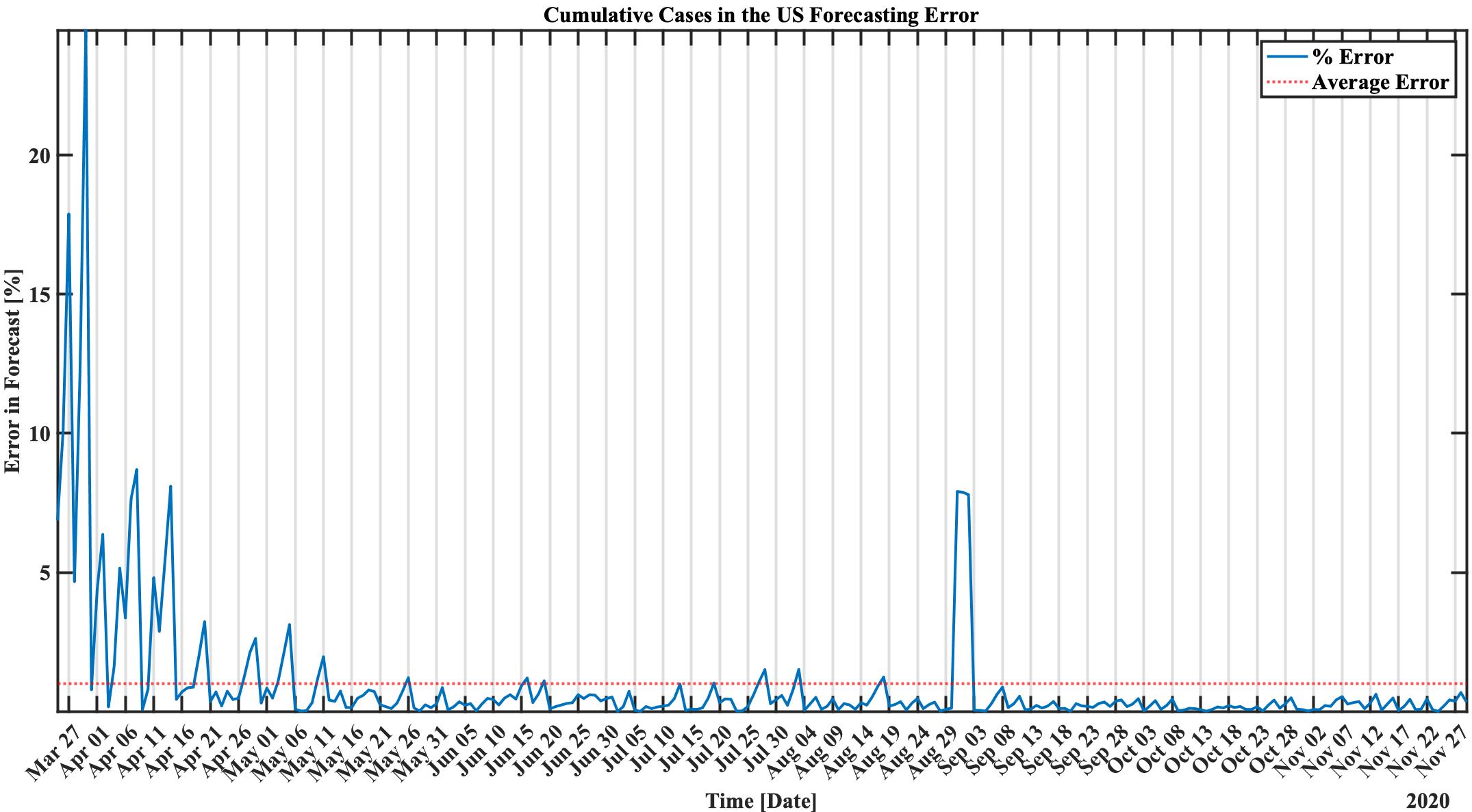}
\caption{Prediction error for US number of COVID-19 cases.}
\label{us-err}
\end{subfigure}
\caption{Prediction of confirmed COVID-19 cases utilizing the publicly available COVID-19 data repository provided by Johns Hopkins. The true data ranges between March 22nd, 2020 and November 29th, 2020. We utilize the last three days of data to forecast the following three days of data.
(\ref{worlda}) Predicted conditions and prediction error worldwide on {November 13}. The  widths of the bubbles  represent the number of cases in a region; only regions with more that 100 cases are used and the bubbles are colored according to their relative percent error.
(\ref{us-pred}) Comparison of true and forecast data for cumulative confirmed cases in the US for April to December 2020. The cumulative forecasts shown here were obtained by summing the forecasts of the individual locations, indicating that the region specific forecasts were sufficiently accurate for tracking the cumulative dynamics of the virus in the US.  
(\ref{us-err}) Percent error for the forecasts of the cumulative confirmed cases in the US. On average the percent error is less than 5 percent and although spikes occur, which could be due to changes in testing availability, the algorithm adjusts and the error stabilizes within a short amount of time. Furthermore, Johns Hopkins provided data for around 1787 locations around the United States and we produced forecasts for each of those locations. 
}
\label{world}
\end{figure}

We have deployed the algorithm to assess the global and United States evolution of the COVID-19 pandemic. The evolution of the virus is rapid, and "Black Swans" in the sense of new cases in regions not previously affected appear with high frequency. Despite that, the Koopman Mode Decomposition based algorithm performed well. \\
\indent
In Figure \ref{worlda} we show the worldwide forecast number  of confirmed cases produced by the algorithm for November 13th, 2020. The forecasts were generated by utilizing the previous three days of data to forecast the next three days of data for regions with higher than 100 cases reported. The bubbles in figure \ref{worlda} are color coded according to their relative percent error. As can be observed, a majority of the forecasts fell below below 15\% error. The highest relative error for November 13th, 2020 was 19.8\% which resulted from an absolute error of 196 cases. The mean relative percent error, produced by averaging across all locations, is 1.8\% with a standard deviation of 3.36\% for November 13th, 2020. Overall, the number of confirmed cases are predicted accurately and since the forecasts were available between one to three days ahead of time, local authorities could very well utilize our forecasts to focus testing and prevention measures in hot-spot areas that will experience the highest growth.

A video demonstrating the worldwide forecasts for March 25, 2020 - November 29, 2020 is provided in the \textbf{Supplementary Information} online (Figure \ref{worlda} is a snapshot from that video). Lastly, it is well known that the ability to test people for the virus increased throughout the development of the pandemic and thus resulted in changes in the dynamics of reported cases. Although it is impossible for a data-driven algorithm to account for changes due to external factors, such as increased testing capabilities, it is important that the algorithm be able to adjust and relearn the new dynamics. For this reason, we encourage the reader to reference the video and note that although periods of inaccuracy due to black swan events occur, the algorithm is always able to stabilize and recover. In contrast, since this is at times a rapidly (exponentially) growing set of data,  methods like naive persistence forecast do poorly. \\
\indent {In Figures \ref{us-pred}, \ref{us-err} we show the performance of the prediction for the cumulative data for the US in March-April 2020. It is of interest to note that the global curve is obtained as a sum of local predictions shown in figure \ref{worlda}, rather than as a separate algorithm on the global data. Again, the performance of the algorithm on this nonstationary data is good.}
%
%

\section{Discussion}


In this work, we have presented a new paradigm for prediction in which the central tenet is understanding of the confidence with which the algorithm is capable of predicting the future realizations of a non-stationary stochastic process. 
{Our methodology is based on Koopman operator theory \cite{mezic:2005}. Operator-theoretic methods have been used for detection of change in complex dynamics in the past, based on both Koopman \cite{mezicandbanaszuk:2000,MezicandBanaszuk:2004} and Perron-Frobenius operators \cite{prinz2011markov}. Other methods include variational finite element techniques combined with information theoretic measure (Akaike's information criterion) and maximum entropy principle \cite{metzner2012}.} \\
\indent Our approach to the  problem of prediction of nonstationary processes has several key  ingredients. First, the Koopman operator on the space of the observables is used as a global linearization tool, whose eigenfunctions provide a coordinate system suitable for representation of the observables. Second, in a numerical computation, we lift the available snapshots to a higher dimensional Hankel-Takens structure, which in particular in the case of abundance of data, allows for better numerical (finite dimensional) Rayleigh-Ritz approximation of eigenvalues and eigenvectors of the associated  Koopman operator, as well as the \textsf{KMD}. Third, using our recent implementation of the \textsf{DMD}, we  select the Koopman modes that have smallest residuals, and thus highest confidence, which is the key for the prediction capabilities of the \textsf{KMD}. In the absence of enough modes with reasonably small residuals, i.e. low confidence, we switch to local prediction, with narrower learning windows and shorter lead time. By monitoring the prediction error, the algorithm may return back to global prediction. 

Our methodology is entirely consistent with the typical training/test dataset validation techniques in machine learning. Namely, the globally learned model on the training data is applied to test data for the next time interval. The novelty in our approach is that we constantly check for how well the learned model generalizes, and if it does not generalize well, we restart the learning. One can say that we implemented a feedback loop, within which the machine learning algorithm's generalizability from training to test dataset is constantly checked, and the system adapts to new conditions. Evidence for effectiveness of this procedure  is presented for the COVID-19 prediction example, where we  show how the generalization error  diminishes over time.

\section{Methods}\label{S=Tools-of-trade} 


Our starting assumption is that observed data is generated by a dynamical process realized on some underlying state space.
This is a broad enough assumption to cover data generated by both deterministic and stochastic dynamical systems \cite{korda2018linear}.
 The (internal) state is often inaccessible; instead, an observable (output) is given as a function $f(\x(t))$ of the state vector $\x(t)$.
 \subsection{The Koopman operator and the \textsf{KMD}}\label{SS=KMD-methods}
The Koopman operator family $\Koop^t$, acts on observables $f$ by composition
    $\Koop^t f (\x) = f(\x( t))$. 
It is a global linearization tool: $\Koop^t$ is a linear operator that allows studying the nonlinear dynamics by examining its action on a linear space $\mathcal{F}$ of observables. In data analysis, for the discrete time steps $t_i$, the discrete sequence $\z_i\approx \x(t_i)$, generated as  numerical software output, is then a discrete dynamical system $\z_{i+1}=\DDS(\z_i)$, for which the Koopman operator reads $\Koop f = f\circ\DDS$. \\
\indent The key of the spectral analysis of the dynamical system is a representation of a vector valued  observable $\bff=(f_1,\ldots,f_{d})^T$ as a linear combination of the eigenfunctions $\bfpsi_j$ of $\Koop$. In a subspace spanned by eigenfunctions each observable $f_i$ can be written  as
$f_i(\z) \approx \sum_{j=1}^{\infty}\bfpsi_j(\z) (\mathbf{v}_j)_i$ and thus ({see e.g. \cite{mezic:2005,mezic_annual_reviews}})
\begin{equation}\label{eq:KMD-intro}
\!\bff(\z) = \!\left(\begin{smallmatrix} f_1(\z) \cr \vdots\cr f_{d}(\z)\end{smallmatrix}\right) \approx  \sum_{j=1}^{\infty}  \bfpsi_j(\z) \mathbf{v}_j,\; \mbox{where}\;\; 
\mathbf{v}_j =\! \left(\begin{smallmatrix} (\mathbf{v}_j)_1 \cr \vdots \cr (\mathbf{v}_j)_{d}\end{smallmatrix}\right)\!\!,   
\end{equation}
\noindent then, since $\Koop\bfpsi_j=\lambda_j\bfpsi_j$, we can envisage the values of the observable $\bff$ at the \emph{future} states $\DDS(\z)$, $\DDS^2(\z), \ldots$ by 
\begin{equation}
 (\Koop^k\bff)(\z) \stackrel{\mathrm{\tiny def}}{=} \bff(\DDS^k(\z))\approx  \sum_{j=1}^{\infty}  \lambda_j^k \bfpsi_j(\z) \mathbf{v}_j , \;\;k=1,2,\ldots
\end{equation}
The numerical approximation of KMD can be computed using for example DMD algorithms. Different versions  of the algorithm used in this work are described in details in \textbf{Supporting Information-Methods}.

\subsection{Finite dimensional compression and Rayleigh-Ritz extraction} For practical computation, $\Koop$ is restricted to a finite dimensional space $\mathcal{F}_{\mathcal{D}}$ spanned by the dictionary of suitably chosen functions $\mathcal{D}=\{ f_1, \ldots, f_d\}$, and we use a matrix representation $\KoopM$ of the compression 
$\bfPsi_{\mathcal{F}_{\mathcal{D}}}\Koop_{|\mathcal{F}_{\mathcal{D}}} : \mathcal{F}_{\mathcal{D}}\rightarrow \mathcal{F}_{\mathcal{D}}$, where $\bfPsi_{\mathcal{F}_{\mathcal{D}}}$ is a $L^2$ projection e.g. with respect to the empirical measure defined as the sum of the Dirac measures concentrated at the $\z_i$'s. Since $\KoopM$ is the adjoint of the \textsf{DMD} matrix $\Aop$ associated with the snapshots $\z_i$, the approximate (numerical) Koopman  modes and the eigenvalues are the Ritz pairs (Ritz eigenvalues and eigenvectors) of $\Aop$, computed using the Rayleigh-Ritz method. The residuals of the Ritz pairs can be computed and used to check the accuracy \cite{DDMD-RRR-2018}. See \textbf{Supporting Information-Methods}.

\subsection{The Hankel-\textsf{DMD}\label{SS=Hankel-intro} (\textsf{H-DMD})}
The data snapshots (numerical values of the observables) can be rearranged in a Hankel-Takens matrix structure: for a subsequence (\emph{window}) of $\win$ successive snapshots $\bff_b, \bff_{b+1}, \ldots, \bff_{\win-1}$,  split $\win = m_H + n_H$ and then define new snapshots as the columns $\hb_i$ of the $n_H\times m_H$ Hankel-Takens matrix (see \cite{Tu2014jcd, arbabiandmezic:2017,mezicnum}, and \textbf{Supporting Information})
\begin{displaymath}
\Han = \!\!	
	\begin{pmatrix}
		\bff_{b} & \bff_{b+1} & \cdots & 
		\bff_{b+m_H} \\
		\bff_{b+1} & \bff_{b+2} & \cdots & 
		\bff_{b+m_H+1} \\
		\vdots & \vdots & \ddots &  
		\vdots \\ 
		\bff_{b+n_H-1} & \bff_{b+n_H} & \cdots & 
		\bff_{b+n_H+m_H-1} \\
	\end{pmatrix}
	 \!= \!\begin{pmatrix} \hb_1 & \ldots & \hb_{m_H+1} \end{pmatrix} \! .
\end{displaymath}
Then, for this data we compute  the \textsf{KMD} and use (\ref{eq:predict-intro}) for prediction. 
Predictions of the observables $\bff_i$ are then extracted from the predicted values of the observables $\bfh_i$.

The introduction of Hankel-Takens matrix alleviates issues that arise from  using a basis on a potentially high dimensional space: namely, taking products of basis elements on $1$-dimensional subspaces - for example Fourier basis on an interval in $\mathbb{R}$. Such constructions lead to an exponential growth in the number of basis elements, and the so-called curse of dimensionality. The Hankel-Takens matrix is based on the dynamical evolution of a one or more
observables - functions on state space - that span a Krylov subspace. The idea is that one might start even with a single observable, and due to its evolution span an invariant subspace of the Koopman operator (note the  connection of such methods with the Takens embedding theorem ideas \cite{MezicandBanaszuk:2004,arbabiandmezic:2017,mezicnum}). Since the number of basis elements is in this case equal to the number of dynamical evolution steps, in any dimension, Krylov subspace-based methods do not suffer from the curse of dimensionality.



\section*{Acknowledgements}


This work was partially supported under DARPA contract HR001116C0116, DARPA contract HR00111890033, NIH/NIAAA grant R01AA023667, and DARPA SBIR Contract No. W31P4Q- 21-C-0007. Any opinions, findings and conclusions or recommendations expressed in this material are those of the authors and do not necessarily reflect the views of the DARPA SBIR Program Office. The support of scientific research of the University of Rijeka, project No. uniri-prirod-18-118-1257, and the  Croatian Science Foundation through grant IP-2019-04-6268.

Distribution Statement A: Approved for Public Release, Distribution Unlimited.

\section*{Author contributions statement}

I.M. conceptualized the prediction algorithm, 
analyzed data, and wrote parts of the paper. Z.D. worked on the numerical algorithms and writing of some parts of the paper. N.C.Z. and S.M. designed parts of the prediction algorithm,
participated in developing methodology and in the results, analysis,
contributed to the preparation of the paper. M.F. participated in methodology development, data preparation and analysis, contributed to the preparation of the paper. R.M. helped develop the algorithm and prepared parts of the paper. A..M.A. helped write a part of the manuscript and produced the COVID forecasting results. I.V. and A.A. were responsible for numerical experiments.


\section*{Data Availability Statement}

The raw COVID-19 data is made publicly available by the Center for Systems Science and Engineering (CSSE) at Johns Hopkins University at \url{https://github.com/CSSEGISandData/COVID-19} .
The raw Influenza data is made publicly available by the World Health Organization at \url{https://www.who.int/influenza/gisrs_laboratory/flunet/en/}.
The raw geomagnetic storm data is made publicly available by the National Aeronautics and Space Administration at \url{https://omniweb.gsfc.nasa.gov/form/omni_min.html}.

\section*{Competing interests}

The authors declare no competing interests.

\clearpage

\renewcommand{\theequation}{S\arabic{equation}}
\renewcommand{\thesection}{S\arabic{section}}
\renewcommand{\thetable}{S\arabic{table}}   
\renewcommand{\thefigure}{S\arabic{figure}}
\renewcommand{\thealgorithm}{S\arabic{algorithm}}


\begin{center}
{\Large Supplementary Information for}
\\ \bigskip

{\Large \bf A Koopman Operator-Based Prediction Algorithm and its Application to COVID-19 Pandemic}
\end{center}

\begin{center}
{\large Igor Mezi\'{c}, Zlatko Drma\v{c}, Nelida \v{C}rnjari\'{c}-\v{Z}ic, Senka Ma\'{c}e\v{s}i\'{c}, Maria Fonoberova, Ryan Mohr, \\ Allan M. Avila, Iva Manojlovi\'{c}, Aleksandr Andrej\v{c}uk }
\end{center}

\setcounter{section}{0}
\section{Methods}\label{S=Supplement}
In this section we provide technical details of the numerical spectral analysis of dynamical systems, which is at the core of the prediction algorithms presented in this work. The tools of trade are the Koopman operator, the \textsf{KMD} and the \textsf{DMD} decompositions.
We first present a compact tutorial on the numerical aspects of the Koopman modal analysis of nonlinear dynamical systems, and, then, we provide theoretical underpinnings for the methods presented in the paper. 

{ We do not go into the details of convergence of the Koopman operator approximations utilized in the paper, but we mention the associated work. For example,  for on-attractor evolution, the properties of the Generalized Laplace Analysis (GLA) method acting on $L^2$ functions were studied in \cite{MezicandBanaszuk:2004,levnajic2010ergodic}. The off-attractor case was pursued in \cite{mohrandmezic:2014} in Hardy-type spaces. This study was continued in \cite{Mezic:2019} to construct dynamics-adapted Hilbert spaces. A study of convergence of DMD-type approximations utilized in this paper is provided in \cite{Korda_andMezic:2018}.  There are two types of results presented in these papers 1) Convergence of the spectral objects over an infinite time interval and 2) rate of convergence to spectral objects. Since we are pursuing a finite-time analysis, the results of type 2) are more relevant. The gist of these results is that the convergence rate is $1/n$ for regular dynamics (limit cycles, limit tori), where $n$ is the number of snapshots,  and $1/\sqrt(n)$ for irregular (chaotic) dynamics. This can be improved to even exponential convergence under some conditions \cite{das2018super}. While the proofs indicated here are for GLA, and depend on convergence time averages over trajectories, they can be extended for DMD methods since in general they can be related to time averages over trajectories \cite{arbabiandmezic:2017}.

Since the prediction methods in the paper are tightly connected to detection of ``Black Swan" events, and these are often defined in imprecise terms, we provide a mathematical definition that we utilize in this work for the orientation of the reader:
\begin{definition}
 Let $f:M\rightarrow \mathbb{C}$, $f \in H$ and $K:H\rightarrow H$ an operator from the Hilbert space $H$ to itself. Consider the dynamics given by
$
f'=Kf.
$
Let $\mu$ be an ergodic invariant measure for $K$, i.e. 
\begin{equation}
  \int_M fd\mu=\lim_ {n\rightarrow \infty}\frac{1}{n}\sum_{j=0}^{n-1}K^jf(x),  
\end{equation}
for almost all $x\in M$. Assume the support $S$ of $\mu$ is such that $S \neq M$. A Black Swan event for an observable $f$ is $g \notin f(S)$.
The magnitude of the Black Swan event for observation $g$ is $d(g,f(S))$ where $d$ is a metric, and $f(S)$ is the range of $f$ on $S$.
\end{definition}
The set $M$ does not necessarily need to be the state space. In the examples shown in the paper, the observables are spectral (e.g. the spectral radius), $M=\mathbb{C}$ and the operator $K$ is the operator acting on the spectral objects induced by the Koopman evolution.

For a non-degenerate stochastic process with an ergodic measure $\mu$ it can be unlikely that the support of  $\mu$ is different from $M$. In that case, the "BlackSwannes" of a value can be defined as $d(g,\nu)$, where $d$ is e.g. the Wasserstein distance of a delta distribution at $f$ and $\nu(E)=\mu(f^{-1}(E))$ is the pushforward measure under $f$ \cite{mezic2008uncertainty}. Note that this particular $d$ can be used in the deterministic case described previously.
}

\subsection{Organization of the section} This section  is organized as follows. First, in \S \ref{SS=DS+Koopman} we review the continuous and the discrete autonomous dynamical systems, and the definition of the Koopman operator $\Koop$ on the space of observables. 
The review material is based on the papers on the theory and applications of the Koopman operator \cite{BudisicMezic_ApplKoop}, \cite{Susuki-Mezic-Raak-Hikihara-2016}, \cite{williams2015-EDMD}, and our recent work \cite{DDMD-RRR-2018}, \cite{DMM-2019-DDKSVC-DFT}, \cite{LS-Vandermonde-Khatri-Rao-2018}.
The spatio-temporal  representation of the evolution of the observables using the eigenvalues and eigenvectors of the Koopman operator (Koopman mode decomposition, \textsf{KMD}) is discussed in \S \ref{SS=Spectral-decomposition-observables}.  The next key ingredient, numerical computation of approximate eigenvalues and eigenvectors, is reviewed in \S \ref{S=Numerical-Computation}. In particular, the details of the matrix representation of a compression of $\Koop$ to the subspace of observables are worked out in \S \ref{SS=Compr-dual}; numerical realization of the Koopman mode decomposition is presented in detail in \S \ref{SSS=Eigfun+Eigmode} and \S \ref{SSS=Kryl-Compr}.
In \S \ref{SS=Schmid-DMD} we provide the key elements of the Schmid's \textsf{DMD} algorithm, and in \S \ref{S=DDMD-RRR-2018} its recent enhancement that allows selection of Ritz pairs that can be used for a \textsf{KMD} suitable for prediction. 
In \S \ref{SSS=Space-time-represent-LS} we review the least squares methods for spatio-temporal representation of the snapshots using the selected Koopman modes.
And finally, after having prepared all necessary ingredients, in \S \ref{S=Global_Prediction} we present the global prediction algorithm. The setup of the prediction framework is given in  \S \ref{SS=Prediction-task} which introduces basic notation and \S \ref{SS=Lift+Hankel}, where we lift the data snapshots in a Hankel structure that will be used in the algorithms. In \S \ref{SS=Predict-Idea-Limits}, we discuss the limitations of the numerical realization of the \textsf{KMD} based prediction scheme introduced in  \S \ref{SS=Spectral-decomposition-observables}.  A worked example that illustrates technical details is provided in \S \ref{SS=worked-example}, where we apply the prediction scheme to the spread of the coronavirus disease in European countries.
In \S \ref{SS-BS-detect}, we discuss the problem of Black Swan events
\cite{Taleb-2007} in the training data, and we propose 
a novel technique, based on the method reviewed in \S \ref{S=DDMD-RRR-2018},  for detecting and retouching Black Swan type disturbances of the data. This makes the global prediction scheme resilient to sudden and unpredictable disturbances that have become part of the learning data window. In \S \ref{SS=Local-prediction}, we present a more flexible local prediction scheme that dynamically resizes the learning data windows and the forecasting lead time in an event of sudden changes and large prediction errors.
 \subsection{Setting the scene: Koopman operator}\label{SS=DS+Koopman}
 Consider an autonomous system of differential equations
 \begin{equation}\label{eq:CDS}
 \dot \x(t)=\CDS(\x(t)) \equiv \left(\begin{smallmatrix} \CDS_1(\x(t))\cr \vdots\cr \CDS_N(\x(t))\end{smallmatrix}\right), 
 \end{equation}
 with state space $\mathcal{X}$ and vector-valued nonlinear function $\CDS$. Here $\mathcal{X}$ is a compact smooth $N$-dimensional manifold, endowed with a Borel sigma algebra $\mathcal{B}$, and for simplicity identified with a subset of $\mathbb{R}^N$, with $\CDS : \mathcal{X}\longrightarrow \mathbb{R}^N$.
 The associated flow map $\flow^t : \mathcal{X}\longrightarrow \mathcal{X}$ advances an initial state $\x(t_0)$ forward in time by a time unit $t$, 
 \begin{equation}\label{eq:flow}
 \x(t_0+t)=\flow^t(\x(t_0))=\x(t_0) + \int_{t_0}^{t_0+t}\CDS(\x(\tau))d\tau.
 \end{equation}
 Note that $\flow^{t+s}=\flow^t\circ\flow^s$, where $\circ$ denotes the composition of mappings.
The (internal) state is often inaccessible; instead an observable (output) is given as a function $f : \mathcal{X}\longrightarrow \C$ of the state, where  the class (function space) $\mathcal{F}\ni f$ of observables is appropriately chosen and endowed with a Banach or Hilbert space structure. For more detailed introduction we refer to \cite{book-Comp-op-Fun-spaces}, in particular Chapters II and V, and \cite{book-Op-Th-Erg-Th}. For instance, we can take $\mathcal{F} = L^p(\mathcal{X},\mu)$, $1\leq p\leq \infty$, with an appropriate measure $\mu$ and e.g. for $p=2$ with the corresponding Hilbert space structure.

The Koopman operator semigroup $(\Koop_{\flow^t})_{t\geq 0}$ is defined by
\begin{equation}\label{eq:Koop-C}
    \Koop_{\flow^t} f = f\circ \flow^t,\;\; f\in\mathcal{F}.
\end{equation}
Here we assume that $\flow^t$ preserves sets of measure zero (if $\mu(A)=0$, then $\mu((\flow^t)^{-1})(A)=0$)  and that $\Koop_{\flow^t}$ is defined on the equivalency classes (modulo $\mu$).
It can be considered as a linearization tool for (\ref{eq:CDS}): $\Koop_{\flow^t}$ is a linear operator that allows studying (\ref{eq:CDS}) by examining its action on the infinitely dimensional space $\mathcal{F}$ of observables. If $\flow^t$ is measure-preserving ( $(\forall A\in\mathcal{B})\;\;(\mu((\flow^t)^{-1}(A))=\mu(A)$) then $\Koop_{\flow^t}$ is an isometry. For an introduction to the theory of the Koopman operator on the Banach lattice $L^p$ see 
\cite[Chapter 7]{book-Op-Th-Erg-Th}.

An analogous approach is applicable to a discrete dynamical system 
\begin{equation}\label{eq:DDS}
\z_{i+1}=\DDS(\z_i), 
\end{equation}
where $\DDS : \mathcal{X} \longrightarrow \mathcal{X}$ is a measurable nonlinear map on a state space $\mathcal{X}$ and $i\in\Z$. 
 The Koopman operator $\Koop\equiv \Koop_{\DDS}$ for the discrete system is defined analogously by 
\begin{equation}\label{eq:Koop-D}
    \Koop f = f\circ \DDS,\;\; f\in\mathcal{F}.
\end{equation}
Discrete dynamical systems naturally describe evolution of discrete events, e.g. stock market data, reported cases of influenza illnesses,  or lynx population in Europe, but they are also at the core of numerical analysis of the continuous systems. More precisely, if we run a numerical simulation of (\ref{eq:CDS}) in a time interval $[t_0,t_*]$, the numerical solution is obtained on a discrete equidistant grid  with fixed time lag $\Delta t$: 
\begin{equation}\label{eq:grid-t}
	t_0, \; t_1 = t_0+\Delta t,\; \ldots ,\; t_{i-1}=t_{i-2}+\Delta t,\;  t_{i}=t_{i-1}+\Delta t,\; \ldots 
\end{equation}
In this case, a black-box software toolbox acts as a discrete dynamical system $\z_i=\DDS(\z_{i-1})$ that produces the discrete sequence of $\z_i \approx \x(t_i)$; this is sampling with noise. For $t_i = t_0+i\Delta t$ we have (using (\ref{eq:flow}), (\ref{eq:Koop-C}) and the group property)
\begin{equation}\label{eq:F-Koop-i-C}
f(\x(t_0+i\Delta t)) = (f\circ\flow^{i\Delta t})(\x(t_0)) = (\Koop_{\flow^{i\Delta t}}f)(\x(t_0)) = (\Koop_{\flow^{\Delta t}}^i f)(\x(t_0)),\;\; \Koop_{\flow^{\Delta t}}^i = \Koop_{\flow^{\Delta t}} \circ \ldots \circ \Koop_{\flow^{\Delta t}}.
\end{equation}
On the other hand, using (\ref{eq:Koop-D}), 
\begin{equation}\label{eq:F-Koop-i-D}
f(\z_i) = f(\DDS(\z_{i-1})) = \ldots = f(\DDS^i(\z_{0}))
= (\Koop^i f)(\z_0),
\end{equation}
where $\DDS^{2}=\DDS\circ\DDS$, $\DDS^i = \DDS\circ\DDS^{i-1}$.
Hence, in a software simulation of (\ref{eq:CDS}) with the initial condition $\z_0 = \x(t_0)$,  we have an approximation
\begin{equation}\label{eq:U-approx-Ut}
(\Koop^i f)(\z_0) \approx (\Koop_{\flow^{\Delta t}}^i f)(\z_0), \;\;f\in\mathcal{F},\; \z_0\in\mathcal{X},\;\;i=0, 1, 2, \ldots 
\end{equation}
with the fidelity that depends on the numerical scheme deployed in the software, and which is studied in the shadowing theory,  see e.g. \cite{Pilyugin2011}, \cite{book-Pilyugin-shadowing}.

This can be obviously extended to vector valued observables: for 
$\bfg=(g_1,\ldots, g_d) : \mathcal{X}\longrightarrow \C^d$ we define
\begin{equation}\label{eq:vector-observable}
\Koop_d \bfg = \begin{pmatrix} g_1\circ \DDS \cr \vdots \cr g_d\circ\DDS\end{pmatrix}
= \begin{pmatrix} \Koop g_1 \cr \vdots \cr \Koop g_d\end{pmatrix}.
\end{equation}
The observables can be both purely physical quantities (e.g. temperature, pressure, energy) and mathematical constructs using suitable classes of functions (e.g. (multivariate) Hermite polynomials, radial basis functions). 
In particular, if we set $d=N$, $g_i(\z)=e_i^T \z$, where $\z\in\C^N$, $e_i=(\bfdelta_{ji})_{j=1}^N$, $i=1,\ldots, N$, then $\bfg(\z)=\z$ is called full state observable and $(\Koop_d \bfg)(\z_i)=\z_{i+1}$.


\subsubsection{Spectral decomposition and representation of observables}\label{SS=Spectral-decomposition-observables}
Spectral decomposition of $\Koop$ is the pillar of both theoretical and practical analysis of dynamical systems in the framework of the Koopman linearizations (\ref{eq:Koop-C}), (\ref{eq:Koop-D}). We consider (\ref{eq:Koop-D}) with $\mathcal{F}=L^2(\mathcal{X},\mu)$, where $\mathcal{X}$ is compact and big enough to contain the states.

An eigenpair $(\lambda_j,\psi_j)$ of the eigenvalue $\lambda_j\in\C$ and nonzero function $\psi_j\in\mathcal{F}$ (eigenvector, eigenfunction) satisfies
$$
\Koop \psi_j = \lambda_j \psi_j .
$$
The key of the spectral analysis of the dynamical system is a representation of an observable as a linear combination of the eigenfunctions of $\Koop$. Since the Koopman operator can have continuous spectrum \cite{mezic:2005,mezic2019spectrum}, this might not be always possible. However, an observable that belongs to a subspace spanned by eigenfunctions can be written  as 
\begin{equation}\label{eq:observable-spectral-decomp}
\hb(\z) = \begin{pmatrix} h_1(\z) \cr \vdots\cr h_{\ell}(\z)\end{pmatrix} \approx  \sum_{j=1}^{\infty}  \bfpsi_j(\z) \mathbf{v}_j,\;\; \mbox{where}\;\; 
h_i(\z) \approx \sum_{j=1}^{\infty}\bfpsi_j(\z) (\mathbf{v}_j)_i, \;\;
\mathbf{v}_j = \begin{pmatrix} (\mathbf{v}_j)_1 \cr \vdots \cr (\mathbf{v}_j)_{\ell}\end{pmatrix}. 
\end{equation}
Then we can envisage the values of the observable $\hb$ at the \emph{future} states $\DDS(\z)$, $\DDS^2(\z), \ldots$ by 
\begin{equation}
\hb(\DDS(\z)) = \begin{pmatrix} h_1(\DDS(\z)) \cr \vdots\cr h_{\ell}(\DDS(\z))\end{pmatrix} = \begin{pmatrix} (\Koop h_1)(\z) \cr \vdots\cr (\Koop h_{\ell})(\z)\end{pmatrix} = (\Koop_d\hb)(\z) \approx  \sum_{j=1}^{\infty}  \lambda_j \bfpsi_j(\z) \mathbf{v}_j , \ldots \;\; \hb(\DDS^k(\z)) \approx  \sum_{j=1}^{\infty}  \lambda_j^k \bfpsi_j(\z) \mathbf{v}_j , \ldots
\end{equation}
For more details see e.g.  \cite{Mezic-2019-Spectrum-JNS}, \cite{GIANNAKIS2019338}, \cite{appr-Koop-spectra-MPT}.
If the dynamics is evolving on an attractor, then all eigenvalues are on the unit circle; however we are also interested in an off-attractor analysis.
The mapping $\flow^t$ is thus not  assumed measure preserving, and thus $\Koop$ is not necessarily unitary. 
Detailed analysis of the spectrum in this general case and the function spaces associated with it can be found in \cite{Mezic-2019-Spectrum-JNS}.

\subsection{Numerical computation}\label{S=Numerical-Computation}
For a practical application of the Koopman operator, we need numerical methods to compute its approximate eigenvalues and eigenvectors, and the modes $\mathbf{v}_j$ in  (\ref{eq:observable-spectral-decomp}). We consider the discrete case (\ref{eq:DDS}), (\ref{eq:Koop-D}), (\ref{eq:vector-observable}); for numerical computations with a continuous system we invoke the discretization (\ref{eq:F-Koop-i-C}), (\ref{eq:F-Koop-i-D}), (\ref{eq:U-approx-Ut}).
For more details on the Extended DMD and implementation using the kernel trick see \cite{williams2015-EDMD}.
\subsubsection{The data}\label{SSS-Data}
In a typical data driven setting, we will have a sequence of snapshots, where we use the notion of snapshot as a numerical value of a scalar or vector valued observable at a specific instance in time. We do not assume explicit knowledge of the mappings $\CDS$ (\ref{eq:CDS}) or $\DDS$ (\ref{eq:DDS}). For example, the snapshots may be obtained from high speed camera recording of a combustion process in a turbine \cite{Combustion-Inst-Flame-Images-2016}, or. e.g. as the wind tunnel measurements. A less expensive and less restrictive is a numerical simulation of (\ref{eq:CDS}) represented by (\ref{eq:F-Koop-i-C}), (\ref{eq:F-Koop-i-D}), (\ref{eq:U-approx-Ut}), where we can feed an initial $\z_0$ to a software tool (representing $\DDS$, or its linearization through a numerical scheme encoded in the software toolbox)  to obtain the sequence 
\begin{equation}\label{eq:seq-Uif}
\bff(\z_0)=(\Koop_d^0 \bff)(\z_0),\;\bff(\z_1) = (\Koop_d \bff)(\z_0),\;  \bff(\z_2) = (\Koop_d^2 \bff)(\z_0),\ldots ,  \bff(\z_{M+1}) = (\Koop_d^{M+1} \bff)(\z_0) ,
\end{equation}
where $\bff=(f_1,\ldots, f_d)^T$ is a  vector valued ($d>1$) observable with the action of $\Koop_d$ defined by (\ref{eq:vector-observable}). The time resolution $\Delta t$ can  be set to obtain the desirable numerical accuracy. This can then be repeated for many initial values $\z_0$; if the new initial value is defined as $\widetilde{\z}_0 = p(\z_0)$, then the new simulation data can be  incorporated by adding the new observables $f_i\circ p$ as new components of $\bff$.   In a CFD application, $\bff$ may be the full state observable, and the entries in the state vectors $\z_i$ are then e.g. the values of the pressure and of the components of the velocity at a discrete spatial grid in the physical domain. 

Hence, independent of the underlying process, the numerical data values are of the form of a matrix $\snap$  with columns, respectively, $\bff(\z_0)$, $\bff(\z_{k+1})=(\Koop_d \bff)(\z_{k})$, $k=0,\ldots, M$:
\begin{equation}\label{eq:data-matrix}
\snap = \left(\begin{smallmatrix}
\bff(\z_0) & \bff(\z_1) & \bff(\z_2) & \ldots & \bff(\z_M) & \bff(\z_{M+1})
\end{smallmatrix}\right) = 
\left(\begin{smallmatrix}
f_1(\z_0) & f_1(\z_1) & f_1(\z_2) & \ldots & f_1(\z_M) & f_1(\z_{M+1}) \cr
f_2(\z_0) & f_2(\z_1) & f_2(\z_2) & \ldots & f_2(\z_M) & f_2(\z_{M+1}) \cr
\vdots & \vdots & \vdots & \ldots & \vdots & \vdots \cr
f_d(\z_0) & f_d(\z_1) & f_d(\z_2) & \ldots & f_d(\z_M) & f_d(\z_{M+1}) \cr
\end{smallmatrix}\right) \in \C^{d\times (M+2)},\;\;\z_{k+1}=\DDS(\z_k),\;\; k=0,\ldots, M.
\end{equation}
Although the snapshots are generated by the nonlinear system
 (\ref{eq:DDS}), (\ref{eq:F-Koop-i-D}), the recursion (Krylov sequence) (\ref{eq:seq-Uif}), driven by the linear operator $\Koop_d$ and numerically evaluated along a trajectory initialized at $\z_0$, motivates to seek out a linear operator (matrix) $\Aop\in\C^{d\times d}$ whose action on the available snapshots is given by 
\begin{equation}\label{eq:Af=Uf}
\Aop \bff(\z_k) = (\Koop_d \bff)(\z_k) = \left(\begin{smallmatrix}
(\Koop f_1)(\z_k) \cr \vdots \cr (\Koop f_d)(\z_k)
\end{smallmatrix}\right) = 
\bff(\DDS(\z_{k})),\;\;k=0,\ldots, M.
\end{equation}
Thus, if we set $\X = \snap(1:d,1:M+1)$, $\Y=\snap(1:d,2:M+2)$,  then such an  $\Aop$ would satisfy $\Y = \Aop \X$, and this could be extended linearly to the span of the columns of 
$\X$ by $\Aop (\X v)=\Y v$, $v\in\C^{M+1}$. The action of $\Aop$ outside the range of $\X$ is not specified by the available data.

In general, $\X$ and $\Y$ are not necessarily extracted from a single trajectory (\ref{eq:seq-Uif}, \ref{eq:data-matrix}). The data may consist of several short bursts with different initial conditions, arranged as a sequence of column vector pairs of snapshots $(\x_k,\y_k)$, where $\x_k=\bff(\z_k)$, $\y_k=\bff(\DDS(\z_k))$ column-wise so that a $k$th column in $\Y$ corresponds to the value of the observable in the $k$th column of $\X$ through the action of $\Koop_d$, as in (\ref{eq:Af=Uf}); see \cite{Tu-On-DMD-TheoryAppl-2014}.
Depending on the parameters $d$ and $M$, the matrices $\X$, $\Y$ can be square, tall (more rows than columns), or wide (more columns than rows).
Then, analogously to (\ref{eq:Af=Uf}), we can search for a linear transformation $\Aop$ such that $\Y=\Aop\X$. Such an $\Aop$ may not exist.

However, we can always define a particular matrix $\Aop$ which minimizes $\| \Y - \Aop\X\|_F$. Clearly, if $\X^T$ has a nontrivial null-space, $\Aop$ is not unique; we can choose $B$ so that $B\X=\0$ and thus $(\Aop + B)\X=\Aop\X$. One can impose an additional condition of minimality of $\|\Aop\|_F$, which yields the well known solution $\Aop = \Y \X^\dagger$, expressed using the Moore-Penrose pseudoinverse $\X^\dagger$ of $\X$. This additional constraint, although useful to enforce uniqueness and boundness, has  (to the best of our knowledge) no other useful interpretation in this framework.
If we are interested in approximating some eigenvalues and eigenvectors, then with the restricted information on $\Aop$ (only its action on the range of $\X$ is meaningfully defined), we will use the Rayleigh quotient $\X^\dagger \Aop \X = \X^\dagger(\Aop + B)\X=\X^\dagger \Y$, so this non-uniqueness of $\Aop$ is immaterial. If $\X$ is of full row rank, then the optimal $\Aop$ is unique. If $\X$ is of full column rank, then $\Aop=\Y\X^\dagger$ satisfies $\Y=\Aop\X$ exactly. Throughout this paper we assume that $\X$ is of full (row or column)
rank. However, even when full rank, $\X$ may be severely ill-conditioned so that its numerical rank is lower, which poses notrivial numerical challenges; these issues are addressed in our recent work  \cite{SysId-KIG-2021}.

%
\subsubsection{A dual representation and the compression of \texorpdfstring{$\Koop$}{U}}\label{SS=Compr-dual}
The sequence of vector valued observables (\ref{eq:data-matrix}) naturally describes the (discrete) time dynamics, with the column index representing the timestamp, and row indices representing the scalar observables (e.g. the pressure and the components of the velocity at particular spacial coordinates) that build the vector observable. For instance, in a \textsf{CFD} application, the multi-indexed (2D, 3D) spatial positions are mapped (vectorized) into a column vector in the usual way, and $\snap$ is generated by the software, column by column.

In a dual interpretation (reading) of the data, we can think of each row in (\ref{eq:data-matrix}) as a set of values of an observable, sampled over the spatial domain. In other words, we transpose  the matrix $\snap$ in (\ref{eq:data-matrix}), partition $\widehat{\snap}=\snap^T$ as
\begin{eqnarray}
    \widehat{\snap}(1:M+1,1:d) &=& \left(\begin{smallmatrix}
f_1(\z_0) & f_2(\z_0) & f_3(\z_0) & \ldots & f_d(\z_0) \cr
f_1(\z_1) & f_2(\z_1) & f_3(\z_1) & \ldots & f_d(\z_1) \cr
\vdots & \vdots & \vdots & \ldots & \vdots \cr
f_1(\z_M) & f_2(\z_M) & f_3(\z_M) & \ldots & f_d(\z_M) \cr
\end{smallmatrix}\right)= \X^T, \label{eq:XT}\\
\widehat{\snap}(2:M+2,1:d) &=& \left(\begin{smallmatrix}
f_1(\DDS(\z_0)) & f_2(\DDS(\z_0)) & f_3(\DDS(\z_0)) & \ldots & f_d(\DDS(\z_0)) \cr
f_1(\DDS(\z_1)) & f_2(\DDS(\z_1)) & f_3(\DDS(\z_1)) & \ldots & f_d(\DDS(\z_1)) \cr
\vdots & \vdots & \vdots & \ldots & \vdots \cr
f_1(\DDS(\z_{M})) & f_2(\DDS(\z_{M})) & f_3(\DDS(\z_{M})) & \ldots & f_d(\DDS(\z_{M})) \cr
\end{smallmatrix}\right)  = \Y^T ,\label{eq:YT}
\end{eqnarray}
and consider the action of $\Koop$ on the space $\mathcal{F}_{\mathcal{D}}$ spanned by the dictionary of scalar functions $\mathcal{D}=\{ f_1, \ldots, f_d\}$. 
That is, we seek a matrix representation $\KoopM$ of the compression 
$\bfPsi_{\mathcal{F}_{\mathcal{D}}}\Koop_{|\mathcal{F}_{\mathcal{D}}} : \mathcal{F}_{\mathcal{D}}\longrightarrow \mathcal{F}_{\mathcal{D}}$, where $\bfPsi_{\mathcal{F}_{\mathcal{D}}}$ is a suitable projection with the range ${\mathcal{F}_{\mathcal{D}}}$.
This is the standard construction: we need a representation of $\Koop f_i$ of the form 
\begin{equation}\label{eq:Ufi(x)}
(\Koop f_i)(z)=f_i(\DDS(z))=\sum_{j=1}^d \bfu_{ji} f_j(z) + \rho_i(z),\;\;i=1,\ldots, d, \;\;z\in\mathcal{X}.
\end{equation}
With the data at hand, the projection is feasible only in the discrete (algebraic) sense: we can define the matrix $\KoopM=(\bfu_{ji})\in\C^{d\times d}$ column-wise by minimizing the residual $\rho_i(z)$ in (\ref{eq:Ufi(x)}) over the states $z=\z_k$, using the values
\begin{equation}\label{eq:Ufi(zj)}
(\Koop f_i)(\z_k) = f_i(\DDS (\z_k)),\;\;i=1,\ldots, d; \;\; k=0,\ldots, M.
\end{equation}
To that end,  write the least squares residual
\begin{equation}\label{eq:LS:Upi}
\frac{1}{M+1}\sum_{k=0}^M |\rho_i(\z_k)|^2= 
\frac{1}{M+1}\sum_{k=0}^M | \sum_{j=1}^d \bfu_{ji} f_j(\z_k) - f_i(\DDS(\z_k))|^2  , 
\end{equation}
which is the $L^2$ residual with respect to the empirical measure defined as the sum of the Dirac measures concentrated at the $\z_k$'s, $\bfdelta_{M+1} = (1/(M+1))\sum_{k=0}^M\bfdelta_{\z_k}$. Hence, the columns of the matrix representation are defined as the solutions of the least squares problems 
\begin{equation}\label{eq:LS-Ui}
\int \left| \sum_{j=1}^d \bfu_{ji}f_j - f_i\circ\DDS\right|^2 \! d\bfdelta_{M+1} = \frac{1}{M+1}
\left\|\! \left[
\left(\begin{smallmatrix} f_1(\z_0) & f_2(\z_0) & \ldots & f_d(\z_0) \cr
\vdots & \vdots & \ldots & \vdots \cr
f_1(\z_M) & f_2(\z_M) & \ldots & f_d(\z_M)\end{smallmatrix}\right)\!\!\! 
\left(\begin{smallmatrix} \bfu_{1i}\cr \vdots \cr \bfu_{di}\end{smallmatrix}\right)\! - \! 
\left(\begin{smallmatrix} f_i(\DDS(\z_0)) \cr
\vdots  \cr
f_i(\DDS(\z_M)) \end{smallmatrix}\right)\!
\right]\!\!\right\|_2^2 \!\longrightarrow\!\! \min_{\bfu_{1i},\ldots,\bfu_{di}},
\end{equation}
for $i=1,\ldots, d$.
The solutions of the above algebraic least squares problems for all $i=1,\ldots, d$ are compactly written as the matrix $\KoopM\in\C^{d\times d}$ that minimizes $\| \X^T \KoopM - \Y^T\|_F$, i.e.  
\begin{equation}\label{eq:U=AT}
\KoopM = (\X^T)^\dagger \Y^T \equiv (\Y \X^\dagger)^T = \Aop^T,
\end{equation}
and the action of $\Koop$ can be represented, using (\ref{eq:Ufi(x)}), as
\begin{equation}\label{eq:Uf1_d(z)}
    \Koop \begin{pmatrix}
    f_1(z) & \ldots & f_d(z)
    \end{pmatrix} = \begin{pmatrix}
    f_1(z) & \ldots & f_d(z)
    \end{pmatrix} \KoopM + \begin{pmatrix}
    \rho_1(z) & \ldots & \rho_d(z)
    \end{pmatrix} .
\end{equation}
Similarly as with the computation  of $\Aop$ in \S \ref{SSS-Data}, $\KoopM$ is uniquely determined only if $\X^T$ is of full column rank. Otherwise,  we must proceed carefully when using the spectral data of $\KoopM$ to infer approximate eigenvalues of $\Koop$.  
In particular, if $d>M+1$, $\X^T$ has a nontrivial null-space, and 
 if $\widetilde{\KoopM}$ is another least squares solution, then $\X^T (\widetilde{\KoopM}-\KoopM)=\0$. On the other hand, along the linear manifold $\KoopM + \mathrm{Ker}(\X^T)=\{ \KoopM + {B}\; :\; \X^T{B}=\0\}$,  the Rayleigh quotient (matrix representation of the compression of $\KoopM$ onto the range of $\X$) $\X^\dagger (\KoopM+{B})\X=\X^\dagger \KoopM \X\in\C^{(M+1)\times (M+1)}$ remains uniquely determined. (Note that in the case of complex data we work with the adjoints $\X^*$ and $\Y^*$, instead of $\X^T$ (\ref{eq:XT}) and $\Y^T$ (\ref{eq:YT}), to obtain $\KoopM=\Aop^*=\overline{A}^T$, which is then the matrix representation in the basis of complex conjugate functions $\overline{f}_i$.)

The quality of this finite dimensional approximation of $\Koop$ depends on the selected dictionary of the observables (capturing a nearly invariant subspace that corresponds to the most relevant eigenvalues), as well as on the approximation level of the underlying measure by the empirical one, in particular on the distribution of the $\z_k$'s. For related numerical issues see \cite{SysId-KIG-2021} and for a theoretical study of convergence, see \cite{Korda_andMezic:2018}.
\subsubsection{Computation of eigenfunctions and the Koopman mode decomposition (\textsf{KMD})}\label{SSS=Eigfun+Eigmode}
Next, we describe the framework for practical computation of the modal decomposition from \S \ref{SS=Spectral-decomposition-observables}. It is the classical Rayleigh-Ritz extraction, based on (\ref{eq:Uf1_d(z)}) and the spectral decomposition of $\KoopM$.

Consider first the case $\mathrm{rank}(\X)=d$; then $d\leq M+1$, and $\KoopM$ is uniquely defined, column by column, from the solutions of the least squares problems (\ref{eq:LS-Ui}), for $i=1,\ldots, d$. In this case all $d$ eigenvalues (with the corresponding eigenvectors) are well determined by the data. 
For technical simplicity we assume that $\KoopM$ is diagonalizable, with the spectral decomposition 
$\KoopM = \bfSS \Lambda \bfSS^{-1}$, with $\Lambda=\mathrm{diag}(\lambda_i)_{i=1}^d$, $\bfSS=(\bfss_1,\ldots,\bfss_d)$, $\KoopM \bfss_i = \lambda_i \bfss_i$. We do not assume that the eigenvalues are simple, and in the case of multiple eigenvalues we list them as successive diagonal entries of $\Lambda$. Then, for $z\in\mathcal{X}$,  
\begin{equation}\label{eq:UfV}
\Koop \begin{pmatrix}
f_1(z)& \ldots & f_d(z)
\end{pmatrix} \bfSS = \begin{pmatrix}
f_1(z)& \ldots & f_d(z)
\end{pmatrix} \bfSS \Lambda + \begin{pmatrix}
\rho_1(z)& \ldots & \rho_d(z)
\end{pmatrix} \bfSS,
\end{equation}
and the approximate eigenfunctions of $\Koop$, extracted from the span of $f_1,\ldots, f_d$, are
$$
\begin{pmatrix}
\phi_1(z)& \ldots & \phi_d(z)
\end{pmatrix} \!=\! \begin{pmatrix}
f_1(z)& \ldots & f_d(z)
\end{pmatrix} \bfS \!=\! \begin{pmatrix}
\sum_{i=1}^d f_i(z) \bfSS_{i1} & \ldots & \sum_{i=1}^d f_i(z) \bfSS_{id}
\end{pmatrix},\;\;(\Koop \phi_i)(z) \!=\! \lambda_i \phi_i(z)\! +\! \sum_{j=1}^d \rho_j(z) \bfSS_{ji} .
$$
Following \S \ref{SS=Spectral-decomposition-observables}, we seek a decomposition of observables in terms of the $\phi_i$'s, similar to (\ref{eq:observable-spectral-decomp}).
In a numerical simulation, these eigenfunctions are accessible, as well as the observables, only as the tabulated values for $z\in\{\z_0,\ldots, \z_M\}$:
\begin{equation}\label{eq:eigs-tabulated}
\left(\begin{smallmatrix}
\phi_1(\z_0) & \phi_2(\z_0) & \phi_3(\z_0) & \ldots & \phi_d(\z_0) \cr
\phi_1(\z_1) & \phi_2(\z_1) & \phi_3(\z_1) & \ldots & \phi_d(\z_1) \cr
\vdots & \vdots & \vdots & \ldots & \vdots \cr
\phi_1(\z_{M+1}) & \phi_2(\z_{M+1}) &\phi _3(\z_{M+1}) & \ldots & \phi_d(\z_{M+1}) \cr
\end{smallmatrix}\right)
=
\left(\begin{smallmatrix}
f_1(\z_0) & f_2(\z_0) & f_3(\z_0) & \ldots & f_d(\z_0) \cr
f_1(\z_1) & f_2(\z_1) & f_3(\z_1) & \ldots & f_d(\z_1) \cr
\vdots & \vdots & \vdots & \ldots & \vdots \cr
f_1(\z_{M+1}) & f_2(\z_{M+1}) & f_3(\z_{M+1}) & \ldots & f_d(\z_{M+1}) \cr
\end{smallmatrix}\right) \bfSS =\snap^T \bfSS.
\end{equation}
Let now $\bfg(z)^T = (g_1(z),\ldots, g_d(z))$ be a vector valued observable and let $g_i(z) = \sum_{j=1}^d \gamma_{ji}f_j(z) + e_i(z)$, so that 
$\bfg(z)^T = (f_1(z),\ldots,f_d(z))\Gamma + E(z)$, $\Gamma = (\gamma_{ji})\in\C^{d\times d}$, $E(z)=(e_1(z),\ldots, e_d(z))$. (If $g_i=f_i$, then $\Gamma=\Id_d$ and $E=\0$. If $g_i\in \mathcal{F}_{\mathcal{D}}$, then $E=\0$.)
Hence
\begin{equation}\label{eq:gTVV-1}
\bfg(z)^T = \begin{pmatrix}
f_1(z) & \ldots & f_d(z)
\end{pmatrix} \bfSS \bfSS^{-1} \Gamma + E(z)=\begin{pmatrix}
\phi_1(z) & \ldots & \phi_d(z)
\end{pmatrix} \bfSS^{-1} \Gamma + E(z),\;\;z\in\mathcal{X} .
\end{equation}
Set $\bfV = \Gamma^T \bfSS^{-T} = \begin{pmatrix}\bfv_1 & \ldots & \bfv_d \end{pmatrix}$, where $\bfv_i$ is the $i$th column. Then  
$$
\begin{pmatrix}
g_1(z)\cr\vdots\cr g_d(z)
\end{pmatrix} = \underbrace{\Gamma^T \bfSS^{-T}}_{\bfV} \begin{pmatrix}
\phi_1(z)\cr\vdots\cr \phi_d(z)
\end{pmatrix} + E(z)^T = \sum_{i=1}^d \bfv_i \phi_i(z) + E(z)^T \approx \sum_{i=1}^d \bfv_i \phi_i(z) .
$$
Since $(\Koop \phi_i)(z)\approx \lambda_i \phi_i(z)$, we have 
\begin{equation}\label{eq:modal-decomp-d}
(\Koop_d^k \bfg)(z) =  \begin{pmatrix}
(\Koop^k g_1)(z)\cr\vdots\cr (\Koop^k g_d)(z)
\end{pmatrix} \approx \sum_{i=1}^d \bfv_i \phi_i(z)  \lambda_i^k .
\end{equation}
In the sequel, we use $\Gamma=\Id_d$; thus $\bfV=\bfSS^{-T}$.
We can assume that $\|\bfv_i\|_2=1$, since $\bfv_i\phi_i(z) = (\bfv_i/\|\bfv_i\|_2) (\|\bfv_i\|_2\phi_i(z))$, where $\|\bfv_i\|_2\phi_i$ is again an eigenfunction. To evaluate (\ref{eq:modal-decomp-d}) numerically at $\z_0$, use (\ref{eq:eigs-tabulated}).

If some eigenvalues are multiple, then with any block diagonal nonsingular matrix $D=\bigoplus_k D_k$, that commutes with $\Lambda$, we have $\KoopM = (\bfSS D) \Lambda (\bfSS D)^{-1}$. (The number of the blocks $D_k$ equals the number of different eigenvalues, and the dimensions corresponds to their multiplicities.) If we repeat the same construction with $\widetilde{\bfSS}=\bfSS D$, then the new approximate eigenfunction of $\Koop$ are $(\widetilde{\phi}_1,\ldots,\widetilde{\phi}_d) = (\phi_1,\ldots,\phi_d)D$, and the matrix of the modes is $\widetilde{\bfV}=\bfV D^{-1}$. At the end, we obtain another representation of the sum in (\ref{eq:modal-decomp-d}).

Using (\ref{eq:U=AT}), we conclude that $\Aop \bfSS^{-T} = \bfSS^{-T}\Lambda$, i.e. the columns of $\bfSS^{-T}$ are the  (right) eigenvectors of $\Aop$. Hence, for computing the Koopman modes, we can proceed with computing the eigenvectors of $\Aop$. The eigenvector matrix is necessarily of the form $\bfSS^{-T}D^{-1}$ with some $D=\bigoplus_k D_k$, as above.

Consider now the case $d >M+1=\mathrm{rank}(\X)$. We have $M+1 < d$ Ritz pairs of $\KoopM$, and in the decomposition (\ref{eq:UfV}) the matrix $\bfSS$ is tall rectangular, $d\times (M+1)$, so we cannot immediately insert $\bfSS \bfSS^{-1}$ as in (\ref{eq:gTVV-1}).
To replace the spanning set $f_1, \ldots, f_d$ with $(\phi_1, \ldots, \phi_{M+1})=(f_1, \ldots, f_d) \bfSS$, we must use $\bfSS \bfSS^\dagger \neq \Id_d$. If the full column rank $\bfSS$ is extracted from the range of $\X$, then $\bfSS\bfSS^\dagger\X=\X$. We proceed with the assumption that the data snapshots are real -- the additional goal is to point out that in that case all computation can be done (and in a software implementation it should) in real arithmetic, even if the eigenvalues and eigenvectors (the columns of $\bfSS$) are complex. Since the matrix $\KoopM$ is then real as well, the pair $\Lambda$, $\bfSS$ computed by the Rayleigh Ritz method will be closed under conjugation and can be indexed as follows: if $\lambda_i\in\R$, then $\bfss_i\in\R^d$, and if $\Im(\lambda_i) >0$, then $\lambda_{i+1}=\overline{\lambda}_i$, $\bfss_{i+1}=\overline{\bfss}_i$. Using the identity 
$$
\begin{pmatrix} \bfss_i & \overline{\bfss}_i\end{pmatrix} \begin{pmatrix} 1 & -\imunit \cr 1 & \imunit \end{pmatrix} = \begin{pmatrix} 2\Re(\bfss_i) & 2\Im(\bfss_i)\end{pmatrix}
$$
we immediately conclude that $\bfSS = \widetilde{\bfSS} J$, where $\widetilde{\bfSS}$ is real and $J$ nonsingular. 
(Here $\Re(\cdot)$ and $\Im(\cdot)$ denote the real and the imaginary parts of complex scalars or vectors.)
Hence, $\bfSS\bfSS^\dagger =\widetilde{\bfSS} \widetilde{\bfSS}^\dagger$ is real symmetric and $\X^T \bfSS\bfSS^\dagger = \X^T$. 
On the other hand, in a practical computation, we see the function values only at $z\in\{\z_0,\ldots, \z_{M+1}\}$, and for those values we can use $\bfSS\bfSS^\dagger$ instead of $\bfSS\bfSS^{-1}$ in relation (\ref{eq:gTVV-1}). The rest is straightforward, yielding the modal matrix $\bfV=\Gamma^T\bfSS^{\dagger T}$, and $\bfSS\bfSS^\dagger \Aop \bfSS^{\dagger T}= \bfSS^{\dagger T} \Lambda$. The latter reveals that $\bfSS^{\dagger}$, $\Lambda$ correspond to Ritz pairs of $\Aop$, extracted from the range of $\X$.

\noindent In the next section, we derive the \textsf{KMD} directly from an application of the Rayleigh Ritz procedure to the matrix $\Aop$.
%
\subsubsection{Krylov compression of \texorpdfstring{$\Koop_d$}{Ud} and the \textsf{KMD}}\label{SSS=Kryl-Compr}
Note that, for an $\bff\in\mathcal{F}$, (\ref{eq:seq-Uif}) naturally generates a Krylov sequence of functions $\bff, \Koop_d \bff, \ldots, \Koop_d^M \bff, \Koop_d^{M+1}\bff$,  and that 
\begin{eqnarray}
    \Koop_d \underbrace{\begin{pmatrix}
    \bff & \Koop_d \bff & \Koop_d^2 \bff & \ldots & \Koop_d^M \bff
    \end{pmatrix}}_{\mathcal{K}_{M+1}} &=& \begin{pmatrix}
    \bff & \Koop_d \bff & \Koop_d^2 \bff & \ldots & \Koop_d^M \bff
    \end{pmatrix} C_{M+1}
    + E_{M+1} , \\ 
\Koop_d \mathcal{K}_{M+1}  &=& \mathcal{K}_{M+1} C_{M+1}  + E_{M+1},\;\;C_{M+1} = \left(\begin{smallmatrix}
     0 & 0 & 0 & 0 & \alpha_0 \cr 
     1 & 0 & 0 & 0 & \alpha_1 \cr
     0 & 1 & 0 & 0 & \alpha_2 \cr
     0 & 0 & 1 & 0 & \alpha_3 \cr
     0 & 0 & 0 & 1 & \alpha_M 
    \end{smallmatrix}\right) , \label{eq:Krylov-dec-F}
\end{eqnarray}
where we have written $\mathcal{K}_{M+1}=\begin{pmatrix}
    \bff & \Koop_d \bff & \Koop_d^2 \bff & \ldots & \Koop_d^M \bff
    \end{pmatrix}$,
\begin{equation}\label{eq:UM+1f}
\Koop_d^{M+1}\bff = \sum_{i=0}^M \alpha_i \Koop_d^i \bff + \bfr_{M+1},
\end{equation}
and $E_{M+1} = \begin{pmatrix}
\0 & \bfr_{M+1}\end{pmatrix}$. In (\ref{eq:UM+1f}), $\bfr_{M+1}$ is the residual obtained after projecting $\Koop_d^{M+1}\bff$ onto the subspace spanned by $[\mathcal{K}_{M+1}]=\mathrm{span}\{\Koop_d^{i}\bff,\;\; i=0,\ldots, M\}$. Here we assume that $M+1 < d$ (possibly even $M+1\ll d$), so that we expect nonzero residual $\bfr_{M+1}$. Our earlier full rank assumption on $\X$ implies that its rank is $M+1$.

If $\Proj_{M+1}$ is the orthogonal projector onto $[\mathcal{K}_{M+1}]$, then the compression $\Proj_{M+1} \Koop_{d|[\mathcal{K}_{M+1}]}$ is represented by the matrix $C_{M+1}$. 
If $C_{M+1} \bfv = \lambda \bfv$, where $\bfv=(v_1,\ldots, v_{M+1})^T\neq\0$, then 
\begin{equation}\label{eq:Ritz-v}
\Koop_d (\sum_{i=0}^M v_{i+1} \Koop_d^i \bff ) = \lambda (\sum_{i=0}^M v_{i+1} \Koop_d^i \bff) + v_{M+1} \bfr_{M+1} .
\end{equation}
This means that $\lambda$ and the function $\bfh =\sum_{i=0}^M v_{i+1} \Koop_d^i \bff= \mathcal{K}_{M+1}\bfv$ satisfy $\Koop_d \bfh = \lambda \bfh + v_{M+1} \bfr_{M+1}$, i.e. 
$(\lambda, \bfh)$ is an approximate eigenpair with the residual 
\begin{equation}\label{eq:Ritz-v-r}
\| \Koop_d \bfh - \lambda \bfh\| = |v_{M+1}|  \|\bfr_{M+1}\|
\end{equation}
measured in the norm of the function space $\mathcal{F}$.

Given the data snapshots (\ref{eq:data-matrix}) as the only available numerical information, the coefficients $\bfalpha=(\alpha_1,\ldots, \alpha_M)$ in (\ref{eq:UM+1f}) can be determined using the discretized (algebraic) least squares projection and the notation from \S \ref{SSS-Data} as follows: The least squares error to be minimized is
\begin{equation}
\| (\Koop^{M+1}_d\bff)(\z_0) - \sum_{i=0}^M \alpha_i (\Koop_d^i \bff)(\z_0)\|_2^2 = \| \bff(\z_{M+1}) - \sum_{i=0}^M \alpha_i \bff(\z_i)\|_2^2 = \| \DOl_{M+1} - \X \bfalpha\|_2^2 . 
\end{equation}
If $\X$ is of full column rank, then $\bfalpha =\X^\dagger \DOl_{M+1}$ is the unique solution expressed using the Moore-Penrose pseudoinverse. Hence, for a particular initial $\z_0$, the relation (\ref{eq:Krylov-dec-F}) reads
\begin{equation}
(\Koop_d \mathcal{K}_{M+1})(\z_0) = \Y = \mathcal{K}_{M+1}(\z_0) C_{M+1} +  (\y_{M+1} - \X\X^\dagger \y_{M+1}) \bfe_{M+1}^T \;\; \bfe_{M+1}^T = \begin{pmatrix} 0 & \ldots& 0, 1 \end{pmatrix} .
\end{equation}
On the other hand, by (\ref{eq:Af=Uf}), $(\Koop_d \mathcal{K}_{M+1})(\z_0) = \Aop \mathcal{K}_{M+1}(\z_0)$, and, as a concrete numerical realization of (\ref{eq:Krylov-dec-F}) on the trajectory starting at $\z_0$, we obtain the Krylov decomposition
\begin{equation}
\Y = \Aop\X = \X C_{M+1} + (\y_{M+1} - \X\X^\dagger \y_{M+1}) \bfe_{M+1}^T .
\end{equation}
where $C_{M+1}=\X^\dagger \Aop \X=\X^\dagger\Y$ is the Rayleigh quotient. Note here that the full column rank assumption on $\X$ implies $\X^\dagger\X = \Id$. Also note that here we do not have $\Aop$ explicitly formed, nor we think of it as $\Aop=\Y\X^\dagger$.

Hence, since the residual $\widehat{\mathbf{r}}_{M+1}=\y_{M+1} - \X\X^\dagger \y_{M+1}$ is unlikely to be zero, we can extract from $\X$ only approximate (Ritz) eigenpairs of $\Aop$. To that end, we first compute the eigenvalues and eigenvectors of $C_{M+1}$.
Under the generic assumption that all eigenvalues of $C_{M+1}$ are algebraically simple,\footnote{Since $C_{M+1}$ is an unreduced Hessenberg matrix, its eigenvalues must be of geometric multiplicity one. If $C_{M+1}$ has multiple eigevaules, then its generalized eigenvector matrix is the inverse of the confluent Vandermonde matrix generated by the distinct eigenvalues. The Jordan structure of each multiple eigenvalue consists of a single Jordan block.} its spectral decomposition  is $C_{M+1} = \Vand_{M+1}^{-1}\bfLambda_{M+1} \Vand_{M+1}$, where
\begin{equation}\label{zd:eq:C-V-Lambda}
\bfLambda_{M+1} =\! \begin{pmatrix} \lambda_1 &  & \cr
& \ddots &   \cr 
&        & \lambda_{M+1}\end{pmatrix}\! ,\;\;
\Vand_{M+1} = \!
\begin{pmatrix} 
1 & \lambda_1 & \ldots & \lambda_1^{M} \cr
1 & \lambda_2 & \ldots & \lambda_2^{M} \cr
\vdots & \vdots & \ldots & \vdots \cr
1 & \lambda_{M+1} & \ldots & \lambda_{M+1}^{M} \cr
\end{pmatrix}\!,\;\; \mathrm{det}(\Vand_{M+1})\equiv \prod_{j>k}(\lambda_j-\lambda_k)\neq 0.
\end{equation}
In other words, the eigenvectors of $C_{M+1}$ are  the columns of the inverse of the Vandermonde matrix $\Vand_{M+1}$.
From $\Aop (\X\Vand_{M+1}^{-1}) = (\X\Vand_{M+1}^{-1})\Lambda_{M+1} + \widehat{\mathbf{r}}_{M+1}e^T\Vand_{M+1}^{-1}\approx (\X\Vand_{M+1}^{-1})\Lambda_{M+1}$, we see that the columns $\widehat{\bfv}_i$ of $\widehat{\bfV}=\X\Vand_{M+1}^{-1}=(\widehat{\bfv}_1,\ldots,\widehat{\bfv}_{M+1})$ are approximate eigenvectors of $\Aop$. With an eye towards (\ref{eq:modal-decomp-d}), we write  $\X = \widehat{\bfV} \Vand_{M+1}$, i.e., for $k=0,\ldots, M$,
\begin{equation}\label{eq:W-KMD}
(\Koop_d^k \bff)(\z_0)=\X e_{k+1} = \sum_{i=1}^{M+1} \frac{\widehat{\bfv}_i}{\|\widehat{\bfv}_i\|_2} \|\widehat{\bfv}_i\|_2 \lambda_i^k  = \sum_{i=1}^{M+1} \bfv_i \|\widehat{\bfv}_i\|_2 \lambda_i^k .
\end{equation}
It is precisely this structure that yields the spatio-temporal representation in \S \ref{SS=Spectral-decomposition-observables}. 
Indeed, if we set 
$$
\Phi = \begin{pmatrix}
    \bff & \Koop_d \bff & \Koop_d^2 \bff & \ldots & \Koop_d^M \bff
    \end{pmatrix}\Vand_{M+1}^{-1} =  \mathcal{K}_{M+1}\Vand_{M+1}^{-1} ,
$$
then $\lambda_j$ and the $j$th column $\Phi_{:j}$ of $\Phi$ are a Ritz pair of $\Koop_d$, $\Koop_d \Phi_{:j}\approx \lambda_j \Phi_{:j}$, see (\ref{eq:Ritz-v}), (\ref{eq:Ritz-v-r}). If $\Phi_{:j}=(\varphi_{1j},\ldots,\varphi_{dj})^T$, then $\Koop \varphi_{ij}\approx \lambda_j \varphi_{ij}$, $i=1,\ldots, d$. 
We have for $k=0,\ldots, M$
\begin{equation}\label{eq:Phi-KMD}
\Koop_d^k\bff(\z_0) =\begin{pmatrix}
\Koop^k f_1(\z_0)\cr\vdots\cr\Koop^k f_d(\z_0)
\end{pmatrix} = \sum_{i=1}^{M+1}\Phi_{:i}\lambda_i^k = 
\sum_{i=1}^{M+1} \begin{pmatrix}
\varphi_{1i}(\z_0)\cr\vdots\cr\varphi_{di}(\z_0)
\end{pmatrix}\lambda_i^k, 
\end{equation}
and we can extrapolate this to the future steps by increasing $k$ which amounts to rising the powers of $\lambda_i$. Also note that $\Phi$ evaluated at $\z_0$ equals precisely $\widehat{\bfV}$, so that (\ref{eq:W-KMD}) is a concrete numerical realization of (\ref{eq:Phi-KMD}).
In an ideal situation, $\lambda_j$ is geometrically simple eigenvalue and $\varphi_{ij}$ are nearly collinear for $i=1,\ldots, d$. However, this is not essential for the purposes of snapshots representation and prediction because the action of $\Koop_d$ is component-wise, and each $\varphi_{ij}$ is an approximate eigenfunction of $\Koop$.

This algebraically elegant process has a drawback that becomes apparent when we consider its numerical software implementation.
Vandermonde matrices are notoriously ill-condi\-ti\-oned. Moreover, in case of an off--attractor analysis the values $|\lambda_i|^j$ may vary in size over several orders of magnitude, which poses challenging problems for the finite precision computation. For that reason, an \textsf{SVD} based method of  Schmid \cite{schmid:2010}, designated as \textsf{DMD}, has become the main computational device for the \textsf{KMD}. However, we have recently shown in \cite{DMM-2019-DDKSVC-DFT} that this companion matrix based approach can be implemented more accurately using the \textsf{DFT} and specially tailored algorithms for the Vandermonde and the related Cauchy matrices.

\subsubsection{Schmid's dynamic mode decomposition (\textsf{DMD})}\label{SS=Schmid-DMD}
The Rayleigh-Ritz procedure outlined in \S \ref{SSS=Kryl-Compr} is based on a Krylov sequence, which naturally fits the dynamics of a discrete dynamical system driven by $\Koop$ (in the space of observables). However, it yields a numerically ill-conditioned problem, as a consequence of that very representation.
From a numerical point of view, 
the Rayleigh-Ritz procedure is best executed in unitary/orthonormal bases, so that $\X$ should be replaced with an orthonormal matrix spanning the same subspace. Since $\X$ can be nearly numerically rank deficient (its columns are actually generated by the power method), Schmid \cite{schmid:2010} used the \textsf{PCA} \cite{PCA} with prescribed cutoff threshold to construct the best lower dimensional subspace (i.e. a \textsf{POD} basis) that captures the data, and then used the Rayleigh-Ritz extraction from that subspace.  For the readers' convenience, we briefly review the \textsf{DMD} algorithm; we assume the more general setting where the snapshots are generated with several initial conditions, so that the input data are not necessarily of the form (\ref{eq:seq-Uif}). That is, the matrices $\X$ and $\Y$ are such that, column-wise,  $\x_k=\bff(\z_k)$, $\y_k=\bff(\DDS(\z_k))=\Aop \x_k$; see \cite{Tu-On-DMD-TheoryAppl-2014}. The total number of snapshots (column dimension) is in this general case denoted by $m$; in the case of a single trajectory (\ref{eq:seq-Uif}), $m=M+1$.

The theoretical underpinning is the classical matrix theorem on best low rank approximations.
\begin{theorem}\label{zd:TM:SVD-EYM}(Eckart-Young \cite{Eckart-Young-1936}, Mirsky \cite{Mirsky-1960})
Let the \textsf{SVD} of $\X\in\mathbb{C}^{d\times m}$ be 
$$
\X = U\Sigma \mcV^*,\;\;\Sigma = \mathrm{diag}(\sigma_i)_{i=1}^{\min(d,m)},\;\;\sigma_1\geq\cdots\geq\sigma_{\min(d,m)}\geq 0 .
$$	
For $r\in\{ 1, \ldots, \mathrm{rank}(\X)\}$, define $U_r = U(:,1:r)$, $\Sigma_r=\Sigma(1:r,1:r)$, $\mcV_r=\mcV(:,1:r)$, and $\X_r=U_r\Sigma_r \mcV_r^*$. Then, $\X_r$ is closest matrix of rank at most $r$ to $\X$, in $\|\cdot\|_2$ and the Frobenius norm $\|\cdot\|_F$, i.e.
	\begin{equation}
			\min_{\mathrm{rank}(\Xi)\leq r} \| \X - \Xi \|_2 = \| \X - \X_r\|_2 = \sigma_{r+1} ; \;\;\;\;
\min_{\mathrm{rank}(\Xi)\leq r} \| \X - \Xi \|_F = \| \X - \X_r\|_F = \sqrt{\sum_{i=r+1}^{\min(d,m)} \sigma_i^2} .		
			\label{zd:eq:EYM:2}\\
	\end{equation}
\end{theorem}
Hence, we can replace $\X$ with its best low rank approximation by truncating its \textsf{SVD} $\X = U\Sigma \mcV^* \approx U_r \Sigma_r \mcV_r^*$, where $U_r=U(:,1:r)$ is $d\times r$ orthonormal ($U_r^* U_r=\Id_r$), $\mcV_r=\mcV(:,1:r)$ is $m\times r$, also orthonormal ($\mcV_r^* \mcV_r=\Id_r$), and $\Sigma_r=\mathrm{diag}(\sigma_i)_{i=1}^r$ contains the largest $r$ singular values of $\X$. 
In brief, $U_r$ is the \textsf{POD} basis for the snapshots $\x_1,\ldots, \x_{m}$,
$$
\sum_{i=1}^m \| \x_i - U_r U_r^*  \x_i \|_2^2 = \min_{\Theta^* \Theta=\Id_r}\sum_{i=1}^m \| \x_i - \Theta\Theta^* \x_i \|_2^2.
$$
The index $r$ is selected so that the approximation error (\ref{zd:eq:EYM:2}) is below a user prescribed threshold value, and it is a numerical rank \cite{Golub:1976:RDL:892104} of $\X$.
Now, DMD uses the range of $U_r$ for the Rayleigh-Ritz extraction. The Rayleigh quotient ${A}_r = U_r^* \Aop U_r$ is computed, using 
\begin{equation}\label{zd:eq:A*SVDX}
\Y = \Aop \X \approx \Aop U_r\Sigma_r \mcV_r^* ,\;\;\mbox{and}\;\; \Aop U_r = \Y \mcV_r \Sigma_r^{-1}, 
\end{equation}
as
\begin{equation}\label{zd:eq:Schmid-S}
A_r = U_r^* \Y \mcV_r\Sigma_r^{-1} ,
\end{equation}
which is suitable for data driven setting because it does not use $\Aop$ explicitly. Clearly, (\ref{zd:eq:A*SVDX}, \ref{zd:eq:Schmid-S}) only require that 
$\Y = \Aop \X$; it is not necessary that $\Y$ is shifted $\X$ as in \S \ref{SSS=Kryl-Compr}. 
Each eigenpair $(\lambda, w)$ of $A_r$ generates the corresponding Ritz pair $(\lambda, U_r w)$ for $\Aop$. 
This is the essence of the Schmid's method \cite{schmid:2010}, summarized in  Algorithm \ref{zd:ALG:DMD} below.

\begin{algorithm}[hbt]
	\caption{{$[\bfV_r, \bfLambda_r]=\textsf{DMD}(\X,\Y)$}}
	\label{zd:ALG:DMD}
	\begin{algorithmic}[1]
		\Require \  		
		$\bullet$ $\X=(\x_1,\ldots,\x_m), \Y=(\y_1,\ldots,\y_m)\in \mathbb{C}^{d\times m}$ that define a sequence of snapshots  pairs $(\x_k,\y_k\equiv \Aop \x_k)$. (Tacit assumption is that $d$ is large and that $m \ll d$.)
		%
		\State $[U,\bfSigma, \mcV]=svd(\X)$ ; \Comment{\emph{The thin \textsf{SVD}: $\X = U \bfSigma \mcV^*$, $U\in\mathbb{C}^{d\times m}$, $\bfSigma=\mathrm{diag}(\sigma_i)_{i=1}^m$, $\mcV\in\mathbb{C}^{m\times m}$.}}
		\State Determine numerical rank $r$ ;
		\State Set $U_r=U(:,1:r)$ ;  $\mcV_r=\mcV(:,1:r)$ ;  $\bfSigma_r=\bfSigma(1:r,1:r)$ ; 	
		\State ${A}_r = (({U}_r^* \Y) \mcV_r)\bfSigma_r^{-1}$ ; \Comment{\emph{Schmid's formula for the Rayleigh quotient $U_r^* \Aop U_r$}.}
		\State $[\bfW_r, \bfLambda_r] = \mathrm{eig}(A_r)$ ; \Comment{$\bfLambda_r=\mathrm{diag}(\lambda_i)_{i=1}^r$; $A_r \bfW_r(:,i)=\lambda_i \bfW_r(:,i)$; $\|\bfW_r(:,i)\|_2=1$}
		\State $\bfV_r = U_r \bfW_r$ . \Comment{\emph{Ritz vectors.}}
		\Ensure $\bfV_r=(\bfv_1 \; \ldots \; \bfv_r)$, $\bfLambda_r$ .
	\end{algorithmic}
\end{algorithm}
\noindent Schmid's DMD algorithm has been  notably successful in \textsf{CFD} applications. For more on interesting applications and modifications of the \textsf{DMD}, see e.g. \cite{williams2015-EDMD}, 
\cite{Modal-analysis-fluid-flow-overview-2017}, \cite{Chen:2012jh}, 
\cite{2015arXiv150203854H}, \cite{Dawson2016}, 
\cite{Hemati:2014jm}, \cite{Takeishi-PhysRevE.96.033310}, \cite{Takeishi-ijcai2017-392}, \cite{Takeishi-8296769}, \cite{Takeishi-NIPS2017_6713}, \cite{dmd-control-doi:10.1137/15M1013857}, \cite{varpro-opt-dmd-doi:10.1137/M1124176}.

\subsubsection{Refined Rayleigh-Ritz Data Driven Modal Decomposition (\textsf{RRRDDMD})}\label{S=DDMD-RRR-2018}
Recently, in \cite{DDMD-RRR-2018}, we revisited DMD and introduced several modifications.
First, we show that the residuals $\|\Aop \bfv_i - \lambda_i \bfv_i\|_2$ can be computed in a data driven scenario as well. This allows for selecting good Ritz pairs, with small residuals, which proved to be the key for selecting good modes for the prediction algorithm; see \S \ref{SS=Predict-Idea-Limits}.
Further, we show that the Ritz vectors can be improved by using the well known refinement technique, which we have adapted to the data driven setting of the \textsf{DMD}. 

\begin{algorithm}[ht]
	\caption{$[\bfV_r, \bfLambda_r, \mathrm{rez}_r, \rho_r]\!=\!\mathrm{DDMD\_RRR}(\X,\Y; \epsilon)$ \{\emph{Refined Rayleigh-Ritz Data Driven Modal Decomposition} \cite{DDMD-RRR-2018}\}}
	\label{zd:ALG:DMD:RRR}
	\begin{algorithmic}[1]
		\Require \noindent	
		\begin{itemize} 
			\item $\X=(\x_1,\ldots,\x_m), \Y=(\y_1,\ldots,\y_m)\in \mathbb{C}^{d\times m}$ that define a sequence of snapshots pairs $(\x_k,\y_k\equiv \Aop \x_k)$. (Tacit assumption is that $d$ is large and that $m \ll d$.)
			\item Tolerance level $\epsilon$ for numerical rank determination.		
		\end{itemize}
		\State $\DD_x=\mathrm{diag}(\|\X(:,i)\|_2)_{i=1}^m$; $\X^{(1)}= \X \DD_x^{\dagger}$; $\Y^{(1)} =\Y \DD_x^{\dagger}$;
		\State $[U,\bfSigma, \mcV]=svd(\X^{(1)})$ ; \Comment{\emph{The thin \textsf{SVD}: $\X^{(1)} = U \bfSigma \mcV^*$, $U\in\mathbb{C}^{d\times m}$, $\bfSigma=\mathrm{diag}(\sigma_i)_{i=1}^m$}.}
		\State Determine numerical rank $r$, with the threshold $\epsilon$. See \cite[\S 3.1.1]{DDMD-RRR-2018} .
		\State Set $U_r=U(:,1:r)$; $\mcV_r=\mcV(:,1:r)$; $\bfSigma_r=\bfSigma(1:r,1:r)$; 	
		\State ${B}_r = \Y^{(1)} (\mcV_r\bfSigma_r^{-1})$; \Comment{\emph{Schmid's data driven formula for $\Aop U_r$}.}
		\State $[Q_r, R]=qr(\begin{pmatrix} U_r, & B_r\end{pmatrix})$; \Comment{\emph{The thin \textsf{QR} factorization: $\begin{pmatrix} U_r, & B_r\end{pmatrix} = Q_rR$; $Q_r$ not computed}.}
		\State $A_r = \mathrm{diag}(\overline{R_{ii}})_{i=1}^r R(1:r,r+1:2r)$; \Comment{\emph{$A_r = U_r^* \Aop U_r$ is the Rayleigh quotient}.}
		\State $\bfLambda_r = \mathrm{eig}(A_r)$ \Comment{$\bfLambda_r=\mathrm{diag}(\lambda_i)_{i=1}^r$; \emph{Ritz values, i.e. eigenvalues of $A_r$}.}
		\For{$i=1,\ldots, r$}
		\State $[\sigma_{\lambda_i},w_{\lambda_i} ]=svd_{\min}(\left( \begin{smallmatrix} R(1:r,r+1:2r)-\lambda_i R(1:r,1:r)\cr R(r+1:2r,r+1:2r) \end{smallmatrix}\right))$;
				 \Comment{\emph{Min. singular value and the corr. right sing. vector, see \cite[\S 3.3]{DDMD-RRR-2018}.}}
		\State $\bfW_r(:,i)=w_{\lambda_i}$ ; 
		$\mathrm{rez}_r(i) = \sigma_{\lambda_i}$; \Comment{\emph{Optimal residual.}}
		\State $\rho_r(i)=w_{\lambda_i}^* A_r w_{\lambda_i}$; \Comment{\emph{Rayleigh quotient, $\rho_r(i)= (U_r w_{\lambda_i})^* \Aop (U_r w_{\lambda_i})$}.}
		\EndFor
		\State $\bfV_r = U_r \bfW_r$; \Comment{\emph{Refined Ritz vectors}.}
		\Ensure $\bfV_r=(\bfv_1 \; \ldots \; \bfv_r)$, $\bfLambda_r$, $\mathrm{rez}_r$, $\rho_r$.
	\end{algorithmic}
\end{algorithm}


\subsubsection{Spatio-temporal  representation of the snapshots}\label{SSS=Space-time-represent-LS}


In general, a \textsf{DMD} algorithm will compute $r\leq m$ Ritz vectors (modes) with the corresponding eigenvalues. In particular, in the Schmid's \textsf{DMD}, $r$ may be considerably smaller than $m$, as e.g. in the case of an off-attractor analysis of a dynamical systems, after removing peripheral eigenvalues, see \cite[\S 4.1]{DDMD-RRR-2018}. In any case, the most important coherent structures of the process are determined by a subset of the modes; so we may want to express the available data snapshots by $r<m$ modes. It is desirable that such modes can represent the snapshots reasonably well, and that they have small residuals which, as we shall see below, is essential for the prediction of the evolution of the sequence (\ref{eq:seq-Uif}), with $m=M+1$.
Assume that we have such a selection of $r$ numerically linearly independent modes and, to ease the notation, assume that we have indexed the Ritz pairs so that the selected ones are indexed with $j=1,\ldots, r$. With this setup,  a modal decomposition of $\bff_k=\bff(\z_k)$ can be written as 
\begin{equation}\label{eq:hi-reconstruct-appr}
\bff_k \approx \sum_{j=1}^{r} \lambda_j^{k}\alpha_j \rv_j,\;\;k=0,\ldots, M+1.
\end{equation}
If $r=M+1=m$, then the coefficients $\alpha_1, \ldots, \alpha_{m}$ can be computed as 
\begin{equation}\label{eq:alpha-pinv-V}
(\alpha_j)_{j=1}^m=\RV_m^\dagger \bff_0, 
\end{equation}
so that this reconstruction is exact for $k=0,\ldots, M$.
In matrix notation, if we define $\RV_r =\begin{pmatrix} \rv_1 & \ldots & \rv_r\end{pmatrix}$ then we have 
\begin{equation}\label{eq:reconstruction-formula-ell}
\begin{pmatrix} \bff_0 & \bff_1 & \ldots & \bff_{M+1}\end{pmatrix}\approx \begin{pmatrix} \rv_{1} & \rv_{2} & \ldots & \rv_{r} \end{pmatrix} \begin{pmatrix} 
{\alpha}_{1} &  &  & \cr
& {\alpha}_{2} &  &   \cr 
&  & \ddots &        \cr 
&        &  & {\alpha}_{r}\end{pmatrix}
\begin{pmatrix} 
1 & \lambda_{1} & \ldots & \lambda_{1}^{M+1} \cr
1 & \lambda_{2} & \ldots & \lambda_{2}^{M+1} \cr
\vdots & \vdots & \ldots & \vdots \cr
1 & \lambda_{r} & \ldots & \lambda_{r}^{M+1} \cr
\end{pmatrix} \equiv \RV_{r} D_{\alpha} \Vanderm_{r, M+2} .
\end{equation}
To compensate for the truncation error, the coefficients $\alpha_j$ can be recomputed by solving the weighted least squares problem
\begin{equation}\label{eq:SP-DMD}
(\alpha_1,\ldots,\alpha_r) = \mathrm{arg min}_{\alpha_j}\sum_{k=0}^{M+1} w_k^2 \| \sqrt{\Omega}(\bff_k - \sum_{j=1}^r \rv_j \alpha_j \lambda_j^{k})\|_2^2 , 
\end{equation}  
where $w_k \geq 0$ are the weights that can be used to emphasize importance of some time indicies or to introduce forgetting factors, $\Omega$ is positive definite matrix,\footnote{In fact, we allow also a diagonal semidefinite matrix $\Omega$ as a mean to exclude selected components of the $\bff_k$'s from the minimization (\ref{eq:SP-DMD}).} and $\sqrt{\Omega}$ stands for the positive definite square root or the Cholesky factor of $\Omega$. 
For numerical methods for this optimization problem we refer to \cite{Jovanovic:2014ft}, \cite{LS-Vandermonde-Khatri-Rao-2018}. Here, for the reader's convenience, we provide an explicit formula for $\Omega=\Id_d$:
\begin{equation}\label{eq:alpha:normal:eq:orig-data}
(\alpha_1,\ldots,\alpha_r)^T = [ (\RV_{r}^*\RV_{r}) \odot (\overline{\Vanderm_{r,M+2}\WW^2 \Vanderm_{r,M+2}^*}) ]^{-1} [ (\overline{\Vanderm_{r,M+2}\WW}\odot (\RV_{r}^* \X\WW))\eb ] ,
\end{equation}
where $\WW=\mathrm{diag}(w_k)$, $\eb = ( 1, 1, \ldots, 1)^T$, and $\odot$ denotes the Hadamard matrix product; see \cite[\S 3.2]{LS-Vandermonde-Khatri-Rao-2018}.

\subsection{Global Koopman prediction algorithm}\label{S=Global_Prediction}
Here we give the details of the new proposed algorithm, designated as \emph{Global Koopman Prediction} (\textsf{GKP}) algorithm.  The basic idea is to extract the intrinsic eigenvalues and modes of the dynamical system under consideration, and then to predict the evolution of the system by using the principles outlined in \S \ref{SS=Spectral-decomposition-observables}, \S \ref{SSS=Eigfun+Eigmode}.  In order to reveal the relevant eigenvalues and the corresponding modes that capture the dynamics of the system on a larger time interval and not only locally, one has to use large training sets. This strategy is at risk if the algorithm is oblivious to unusual and unexpected changes (perturbations) that can be classified as Black Swan events. If such data are used in a learning window, the long term prediction is doomed to fail. We use the numerically computed spectral information on $\Koop$ to develop an additional device to equip the algorithm with a \emph{litmus test} for detecting Black Swan events {(a posteriori, of course)}, and, moreover, with a retouching scheme to restore the global prediction capability (see \S 2.2 in the main text).  Furthermore, an important feature is that the data snapshots are lifted in a Hankel matrix structure, as described in \S \ref{SS=Lift+Hankel}.


\subsubsection{Setting the scene - the prediction task}\label{SS=Prediction-task}

 Consider a discrete dynamical system $\z_{k+1}=\DDS(\z_k)$
that is accessible through a sequence of snapshots
\begin{equation}\label{eq:observables-ui}
\bff_0, \bff_1, \bff_2, \ldots  \;\;\; ( \bff_k\in\R^{\dimn},\;\; k = 0, 1,  2, \ldots, M, M+1, \ldots )
\end{equation}
where $d\geq 1$ is the dimension of the scalar or vector-valued \emph{system observable} $\bff : \mathcal{X}\rightarrow \R^{\dimn}$, and $\bff_k=\bff(\z_k)$ is its value for the (possibly unknown) state $\z_k$, with the time stamp $t_k$, $k=0, 1, 2, \ldots$. The goal is to learn the dynamics from the available data and then to predict the future values.

More precisely, suppose that the \emph{present time moment} is $t_{p-1}$, and that,  up to that moment, the data is readily available; we seek a prediction of the data value at the next time moments $t_p, t_{p+1}, \ldots, t_{p+\tau(p)}$. We call that future time moments $t_p$ the \emph{prediction moments}.
The prediction  will be based on a \emph{sliding window} of size $\win$ in the sequence (\ref{eq:observables-ui}), i.e.  we will use the $\bff_k$'s starting from the index $b=p-\win$ that defines 
 the \emph{active window} $\mathcal{W}(p,\win)$ of consecutive data (\emph{training set}) with indices $b=p-\win, b+1, \ldots, b+\win-1=p-1$. 
In terms of the system mapping $\DDS$, these values can be represented as
 \begin{equation}\label{eq:training-set}
 \bff_b= (\bff\circ\DDS^b)(\z_0),\; \bff_{b+1} = (\bff\circ\DDS^{b+1})(\z_0), \ldots,\; \bff_{b+\win-1}=(\bff\circ\DDS^{b+\win-1})(\z_0).
 \end{equation}

\subsubsection{Lifting the data into a Hankel matrix structure and the \textsf{H-DMD}}\label{SS=Lift+Hankel}
The key for a successful application of the prediction framework from \S \ref{SS=Spectral-decomposition-observables} is that the finite dimensional numerical approximation from  \S \ref{SSS=Eigfun+Eigmode} captures the spectral information accurately enough. 
To that end, we adopt the Hankel-\textsf{DMD} approach of \cite{Tu2014jcd}, \cite{ArbabiMezic2016}.
For an active window $\mathcal{W}(p,\win)$, first conveniently split $\win = m_H + n_H$, and then
 lift the observables  into the higher dimensional space $\R^{\ell}$, $\ell=\dimn\cdot n_H$, and arrange them as columns of a $\ell \times (m_H+1)$ (block) Hankel matrix as follows:

\begin{equation}\label{eq:Hankel}
\Han = 	\left(
	\begin{array}{ccccc}
		\bff_{b} & \bff_{b+1} & \cdots & \bff_{b+m_H-1} & \bff_{b+m_H} \\
		\bff_{b+1} & \bff_{b+2} & \cdots & \bff_{b+m_H} & \bff_{b+m_H+1} \\
		\vdots & \vdots & \ddots &  \vdots & \vdots \\ 
		\bff_{b+n_H-1} & \bff_{b+n_H} & \cdots & \bff_{b+n_H+m_H-2} & \bff_{b+n_H+m_H-1} \\
	\end{array}
	\right) = \begin{pmatrix} \hb_1 & \hb_2 & \ldots & \hb_{m_H+1} \end{pmatrix}.
\end{equation}
We can think of the $\hb_i$'s as the values of the vector-valued observable $\hb : \mathcal{X}\longrightarrow \R^{\ell}$ composed with the powers of $\DDS$ analogously to (\ref{eq:training-set}), i.e. $\hb = \begin{pmatrix} \bff\circ \DDS^b & \bff\circ \DDS^{b+1} & \ldots &  \bff\circ\DDS^{b+n_H-1}\end{pmatrix}^T$ and 
\begin{eqnarray}
  \begin{pmatrix} \hb_1 & \hb_2 & \hb_3 & \ldots & \hb_{m_H+1} \end{pmatrix} &=&   
\begin{pmatrix} \hb(\z_b) & \hb\circ \DDS(\z_b) & \hb\circ \DDS^2(\z_b) & \ldots &  \hb\circ \DDS^{m_H}(\z_b)\end{pmatrix} \\
&=&   
\begin{pmatrix} \hb(\z_b) & (\Koop_{\ell}\hb)(\z_b) & (\Koop_{\ell}^2\hb)(\z_b) & \ldots & (\Koop_{\ell}^{m_H}\hb) (\z_b)\end{pmatrix} ,
\end{eqnarray}
which can be interpreted as a Krylov sequence for the Koopman operator $\Koop_{\ell}=\otimes_{1}^{\ell}\Koop$, $\Koop_{\ell}\hb = \hb\circ\DDS$.
The techniques from \S \ref{S=Numerical-Computation} now apply in this new setting simply by setting $\hb$ instead of $\bff$, and $\ell$ instead of $d$. The matrix $\Han$ plays the role of the snapshot matrix\footnote{Note that this is different from a system identification technique based on the SVD decomposition of $\Han$.} $\snap$, and we have $\X=\Han(:,1:m_H)$, $\Y=\Han(:,2:m_H+1)$.  We will attempt predicting the $\hb_k$'s, and from the obtained results extract the predictions of the original observable $\bff$. The starting points are the \textsf{DMD} of $\Han$ (\textsf{H-DMD}), and the corresponding \textsf{KMD}. Since the \textsf{KMD} changes with the sliding active data window $\mathcal{W}(p,\win)$, we use the term active \textsf{KMD} (\textsf{AKMD}) when we refer to the computation used for prediction.

\subsubsection{Prediction - the basic idea and its limitations}\label{SS=Predict-Idea-Limits}
%

Suppose that the \textsf{DMD} algorithm, applied to (\ref{eq:Hankel}),  has extracted  $r=m_H$ Ritz pairs, and that the Ritz vector span the range of $\Han$.  Then we can determine the coefficients $\alpha_j$ such that  
\begin{equation}\label{eq:hi-reconstruct}
\hb_k = \Aop^{k-1}\hb_1 = \sum_{j=1}^{m_H} \lambda_j^{k-1}\alpha_j \rv_j + \bfdelta_{k,m_H+1} \widehat{\bfr}_{m_H+1},\;\;k=1,\ldots, m_H+1,
\end{equation}
where $\bfdelta_{k,m_H+1}$ is the Kronecker delta symbol, and $\widehat{\bfr}_{m_H+1}$ is the residual of the orthogonal projection of $\hb_{m_H+1}$ onto the range of $\X$. This means that the decomposition of the snapshots in terms of the modes is exact, except for the last one, which may not belong to the range of $\X$, and the residual $\widehat{\bfr}_{m_H+1}$ represents its decomposition error.
For details we refer to \cite[\S 2.4, \S 3.2]{DMM-2019-DDKSVC-DFT}.
If we want to extend the above relation beyond the index $k=m_H+1$ (i.e. to extrapolate into the future the evolution of the sequence $\hb_k$), we can apply the appropriate power of $\Aop$ and use the approximation $\Aop \rv_j = \lambda_j \rv_j + \rez_j \approx \lambda_j \rv_j$. This is a fairly simple operation - it amounts to increasing the power of the $\lambda_j$'s.  
Of course, the residuals will accumulate with each such iteration, e.g.
\begin{eqnarray}
\!\!\!\!\!\!\!\!\!\Aop \hb_{m_H+1} &\!\!\!\!=\!\!\!\!& \sum_{j=1}^{m_H} \lambda_j^{m_H+1}\alpha_j \rv_j + \sum_{j=1}^{m_H} \lambda_j^{m_H}\alpha_j \rez_j + \Aop \widehat{\bfr}_{m_H+1}, \label{eq:pred+rezd-1}\\
\!\!\!\!\!\!\!\!\!\Aop^2 \hb_{m_H+1} &\!\!\!\!=\!\!\!\!& \sum_{j=1}^{m_H} \lambda_j^{m_H+2}\alpha_j \rv_j + 
\sum_{j=1}^{m_H} \lambda_j^{m_H+1}\alpha_j \rez_j + \sum_{j=1}^{m_H} \lambda_j^{m_H}\alpha_j \Aop \rez_j + \Aop^2 \widehat{\bfr}_{m_H+1},\;\;\Aop^3 \hb_{m_H+1} = \ldots \label{eq:pred+rezd-2}
\end{eqnarray}
So, using the first sums above (and ignoring the residual terms) to predict future of the $\hb_k$'s has limited range, except in the case of small $\widehat{\bfr}_{m_H+1}$ and small residuals $\rez_j$, which are not too much amplified under the action of the powers of $\Aop$. 
Hence, it is desirable to have a \textsf{KMD} that uses only the selected modes corresponding to the Ritz pairs with small residuals, and that we can have an accurate decomposition of the type (\ref{eq:reconstruction-formula-ell}), using the selected modes. This selection is possible in data driven scenarios using the methods from \cite{Jovanovic:2014ft} and \cite{DDMD-RRR-2018}, outlined in \S \ref{S=DDMD-RRR-2018}.
The desire to have a high fidelity representation 
of the data snapshots
with as few as possible  modes $\rv_j$ is motivated by revealing latent coherent structures of the e.g. flow field; small residuals allow for extrapolation of the dynamics forward in time.

If our goal is solely the prediction, the weight factors $w_i$ in (\ref{eq:SP-DMD}) can be tuned to favor most recent snapshots, and the weighting matrix $\Omega$ can emphasize particular block rows in the $\hb_k$'s;  see \S \ref{SSS=Space-time-represent-LS} and \cite[\S 3]{LS-Vandermonde-Khatri-Rao-2018}. In particular, with a suitable choice of $\WW$ and $\Omega$, we can focus the reconstruction of the $\hb_k$'s to the present snapshot $\bff_{n+n_H+m_H-1}=\bff_{t_{p-1}}$.
%
%

Since the $\bff_k$'s, starting from the past time stamp index $b$ and ending at the present index $p-1=b+n_H+m_H-1$, are in the last block row of  $\Han$ (see (\ref{eq:Hankel})), the corresponding formulas are obtained by taking the last $\dimn$ components of the Ritz vectors $\rv_j$. To that end, define $\widehat{\rv}_j = \rv_j((n_H-1)\dimn +1: n_H \dimn)$ as  the trailing $\dimn$ components of $\rv_j$.
Hence, from the \textsf{AKMD} of the lifted observables, we read off approximate decomposition of the snapshots $\bff_k$ as 
\begin{equation}\label{eq:prediction}
	\widetilde{\bff}_k = \sum_{j=1}^r \widehat{\rv}_j \alpha_j \lambda_j^{k-1}. 
\end{equation} 
For $k=b+n_H,...,b+n_H+m_H-1$, (\ref{eq:prediction}) is a reconstruction of the acquired data, while for $k=p,p+1,\ldots, p+\tau(p)$, (\ref{eq:prediction}) is an extrapolation of the \textsf{AKMD}, and it gives us  predictions for future data snapshots. We say that $\widetilde{\bff}_{p+\tau}$, $\tau > 0$ is the prediction of the observable at the lead time $\tau$.  

If the number of rows $n_H$ of Hankel matrix is smaller than number of columns $m_H$ $(n_H < m_H)$ then the \textsf{KMD} gives some sort of regression function for the data in the reconstruction window
$b+n_H, \ldots, b+n_H+m_H-1$.
The reason why the first part of active window is not declared as reconstruction window is that when the \textsf{KMD} of the form (\ref{eq:prediction}) is used, the {Koopman eigenfunction values}
 are determined such that for $t=0$ the sum of the right hand side is equal to the first snapshot. When applied to Hankel matrix this means that the data in the beginning of the active window $b,\ldots, b+n_H-1$, which form the first column of Hankel matrix are reconstructed with high accuracy. 
On the other hand, if $n_H \ge m_H$ and if the Hankel matrix has full column rank, the data in the whole active window are reconstructed with high accuracy.
{
\subsection{A worked example}\label{SS=worked-example}
We now illustrate the key elements of the procedure outlined in \S \ref{SS=Lift+Hankel}, \S \ref{SS=Predict-Idea-Limits} using a worked example. The problem under study is the spread of the  coronavirus disease (COVID-19). The data consists of reported cummulative daily cases in some European countries. It should be stressed that the algorithm uses only the raw data -- no other information on the nature of the data or on modelling parameters is assumed. Further, the data itself is clearly not ideal, as it depends on the reliability of the tests, testing policies in different countries (triage, number of tests, reporting intervals,  
reduced testing during the weekends), contact tracing strategies, the number of asymptomatic transmissions  etc.  Moreover, using the data from different countries in the same vector observable poses an additional difficulty for a data driven revealing of the dynamics, because the countries independently and in an uncoordinated manner impose different restrictions, thus changing the dynamics. 

For an analysis of a particular country, it is better to define the observables as the reported cases on local level, e.g. provinces, counties, cities with similar conditions. Clearly, the dynamics of the spread of the infection depends on the population density as well. This is best seen e.g. by comparing the heat map of the reported cases in the USA with the image of the USA from space at night. In the numerical examples in this section, we purposely use data from different countries to make the prediction task more challenging,    which makes it an excellent stress test example.

Our goal with this example is twofold. First, 
we show the potentials and the limits of the proposed prediction algorithm. Secondly, we discuss technical details of the computational scheme. 

We use the following  datasets:
\begin{itemize}
\item[\textbf{DS1}] The numbers of reported COVID-19 cases in Germany, France and the United Kingdom in the period February 29 to November 18. The ordered triple of reported cases is an observable.
\item[\textbf{DS2}] The dataset \textbf{DS1} augmented by the numbers of reported cases in Denmark, Czechia, Slovenia, Austria and Slovakia.
\item[\textbf{DS3}] The numbers of reported COVID-19 cases in a selected European country in the period February 29 to November 18, augmented with two sequences of filtered data.
\end{itemize} 

The test of the prediction algorithm runs on the lifted data ($265$ observables from $\R^d$: $d=3$ for \textbf{DS1} and \textbf{DS3}; $d=8$ for \textbf{DS2}) i.e. on the columns of the $94 d \times 172$ Hankel matrix $\Han=\begin{pmatrix} \hb_1 & \hb_2 & \ldots & \hb_{172} \end{pmatrix}$ (see (\ref{eq:Hankel})) with the block partition $94\times 172$, each block being $d\times 1$. 
The matrix $\Han$ is used as a historical record, encoding the period February 29 -- November 18, and we run the prediction algorithm starting from some past index and test its accuracy by comparison with the historical data. 
We use simple increasing window starting at the index $b=1$ and ending at $p-1$, where we choose different values of $p$. Then we predict the next $\tau+1$ values  from the moment $p$ on. 

In the first experiment, we use \textbf{DS1} and attempt prediction for $35$ days ahead. We take the first $40$ columns of $\Han$ as available data and set  $\X=\Han(1:282,1:39)$, $\Y=\Han(1:282,2:40)$. (This corresponds to the period February 29 -- July 10, and the prediction for $35$ days ahead starts July 11.) The prediction relative error is shown on the left panel in Figure \ref{fig:covid_DS1_40_35}. The right panel shows the Koopman Ritz values computed in the algorithm; note that the algorithm has revealed the eigenvalue $1$, and that all other Ritz values are inside the unit circle. The  quite satisfactory prediction skill (recall, no information whatsoever on the nature of the data is used) and well behaved Ritz values 
are related to the nature of the dynamics of the infection during the summer. 

In this example, it is instructive to check the Government Response Stringency Index
(GRSI) \cite{GRSI}\footnote{For an interactive exploration of the GSRI see \url{https://ourworldindata.org}.} for the entire time interval involved in the computation. The three indexes behave differently: it can be  noticed that France had sharper changes than Germany and the United Kingdom (e.g. around June 20),  and sometimes similar increase of stringency but a week earlier than the other two countries (e.g. beginning to mid March). On the other hand, the GRSI for Germany and the United Kingdom were not that much different throughout the observed period; see the left panel in Figure \ref {fig:covid_GRSI}. 
It should be noted, however, that the GRSI does not measure the quality of the implementation of the imposed restriction and that for a particular country  it does not necessarily indicate the trends in the dynamics of the disease spreading. 
\begin{figure}[ht]
	\centering
	\includegraphics[width=0.48\linewidth]{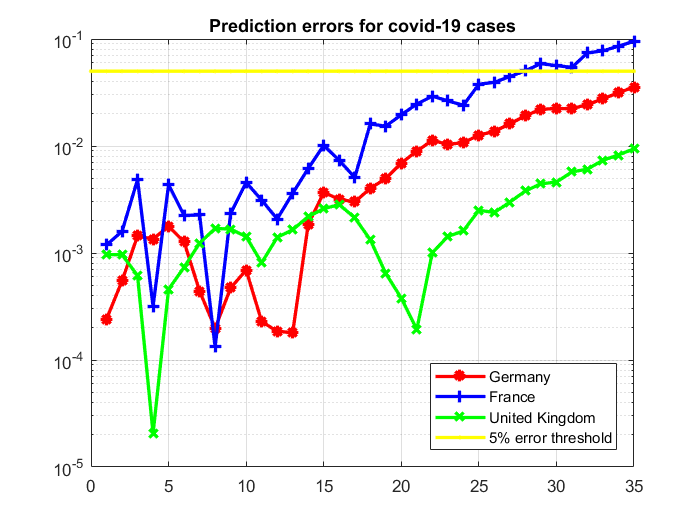}
	\includegraphics[width=0.48\linewidth]{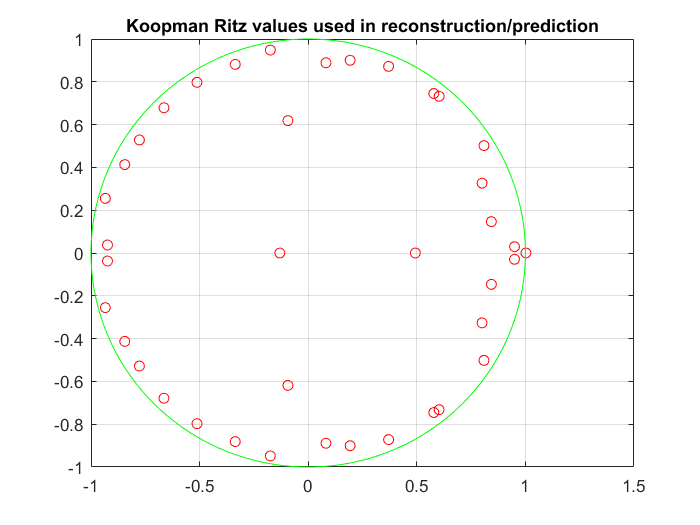}
  \caption{\emph{Left panel}: Relative prediction error for Germany, France and the United Kingdom for a $35$ days prediction starting after the data window $\hb_{1:40}$. (In terms of the original data, prediction starts on July 11.) \emph{Right panel}: The Koopman Ritz values used in the \textsf{KMD}.}
	\label{fig:covid_DS1_40_35}
\end{figure}


\begin{figure}[ht]
	\includegraphics[width=0.48\linewidth, height=2.3in]{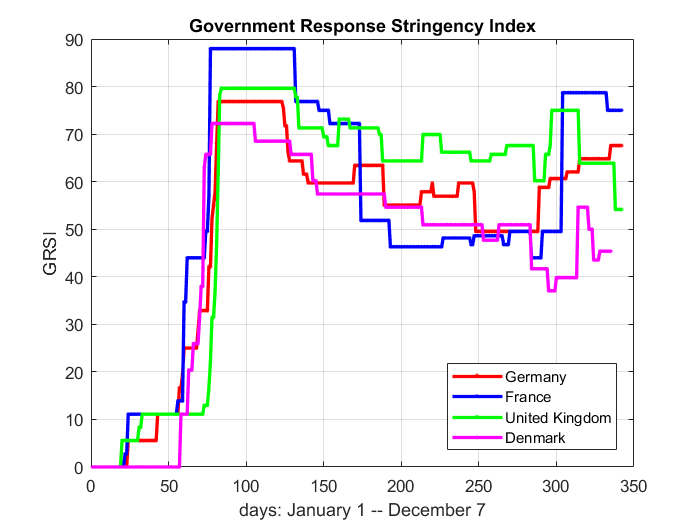}\hspace{-2mm}
	\includegraphics[width=0.48\linewidth, height=2.3in]{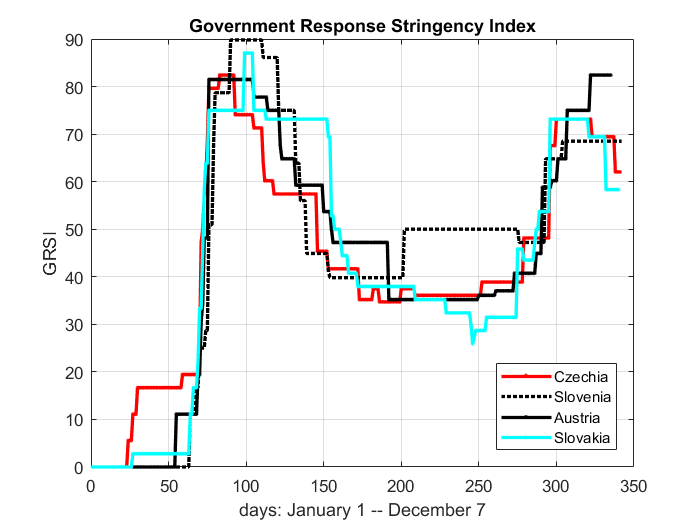}
	\vspace{-2mm}
  \caption{Government Response Stringency Index measures response indicators (OxCGRT indicators) such as school closing, workplace closings, cancelling public events, restrictions on gathering size, closing public transport, stay at home requirement, restriction on internal movement and international travel. The GRSI data for Germany, France, United Kingdom, Denmark, Czechia, Slovenia, Austria and Slovakia are taken from \cite{GRSI}. For more details see \cite{GRSI-2}.}
	\label{fig:covid_GRSI}
\end{figure}

Now, in the same interval, we add new observables by including the data from five more countries: Danemark, Czechia, Slovenia, Austria and Slovakia.
Hence, the matrix $\Han$ is $752\times 172$.
The prediction errors for a $28$ days prediction are given in the left panel in Figure \ref{fig:covid_DS2_40_28}.
Remarkably, the computed Koopman Ritz values nearly match the one computed in the first test with only three countries, see the right panel in 
Figure \ref{fig:covid_DS2_40_28}. Note that even with the differences shown in Figure \ref{fig:covid_GRSI}, the  main trend of the implementation of the measures is similar. This might help explain the robustness of the spectrum indicated in Figure \ref{fig:covid_DS2_40_28}, where such differences do not seem to lead to drastic change in the spectral behavior. We believe this indicates the robustness of our methodology.

\begin{figure}[ht]
	\centering
	\includegraphics[width=0.48\linewidth]{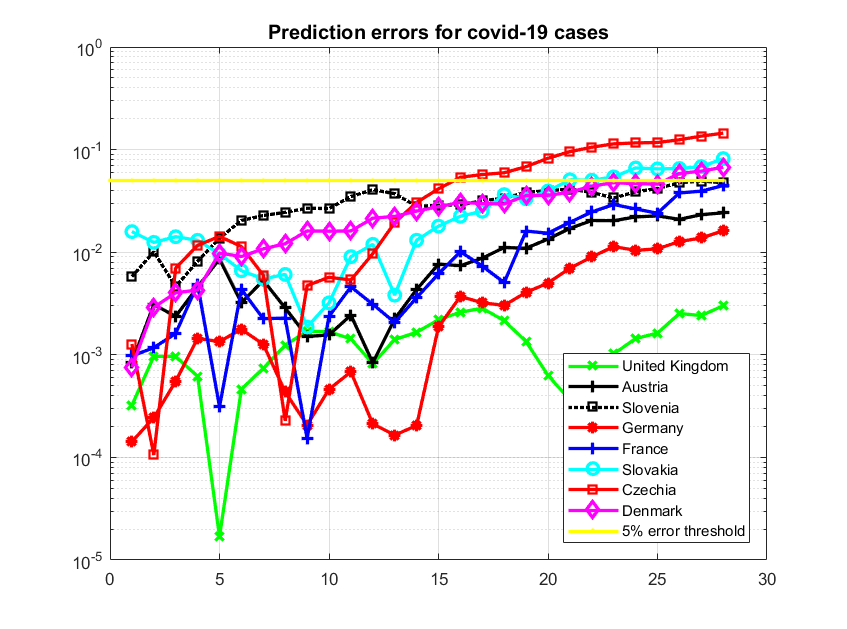}
	\includegraphics[width=0.48\linewidth]{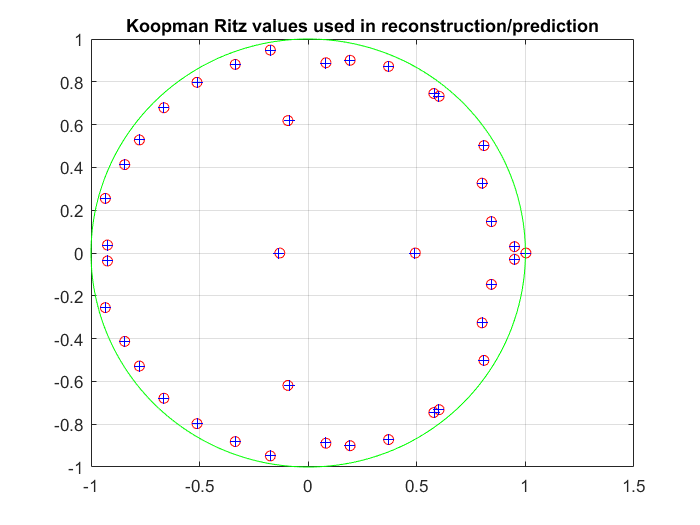}
	\vspace{-4mm}
  \caption{\emph{Left panel}: Relative prediction error for eight European countries for a $28$ days prediction starting after the data window $\hb_{1:40}$. (In terms of the original data, prediction starts on July 11.) \emph{Right panel}: The Koopman Ritz values used in the \textsf{KMD}, denoted as blue pluses (\textcolor{blue}{+}). The red circles (\textcolor{red}{$\circ$}) denote the Ritz values computed using only three countries as shown in Figure \ref{fig:covid_DS1_40_35}.}
	\label{fig:covid_DS2_40_28}
\end{figure}

We proceed with the numerical experiment using the dataset \textbf{DS1}. We further expand the learning window and then  consider three consecutive steps with $\hb_{1:105}$, $\hb_{1:106}$, $\hb_{1:107}$. The relative errors for $35$ days prediction are shown in the first row of Figure \ref{fig:covid_DS1_105_35}.
In the context of policy changes that affected the dynamics of the infection spreading, and the fact that the algorithm is purely data driven, the results can be considered satisfactory: in the first graph, the error is below five percent for $16$ days and below ten percent for three weeks for all three countries (first graph), below six percent for almost entire $35$ days period (second graph), below five percent for more than three weeks and below ten percent for $30$ days (third graph). 

\begin{figure}[ht]
	\centering
	\includegraphics[width=0.32\linewidth]{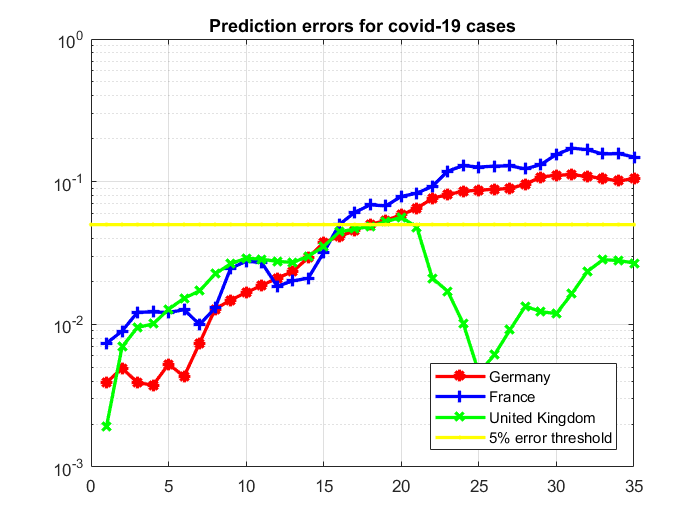}
	\includegraphics[width=0.32\linewidth]{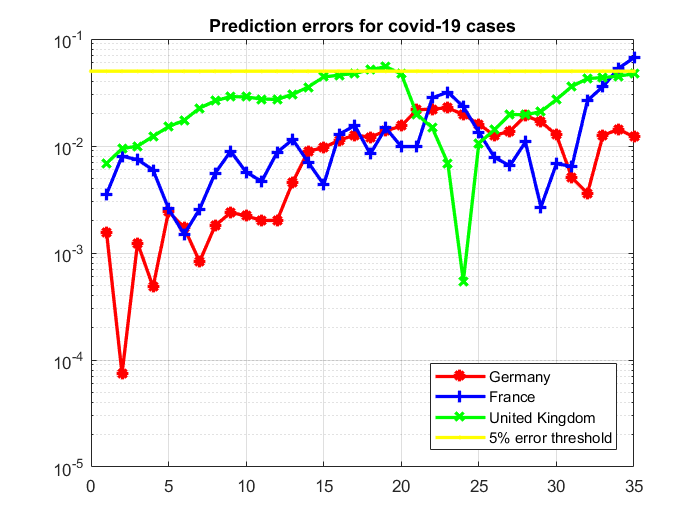}
	\includegraphics[width=0.32\linewidth]{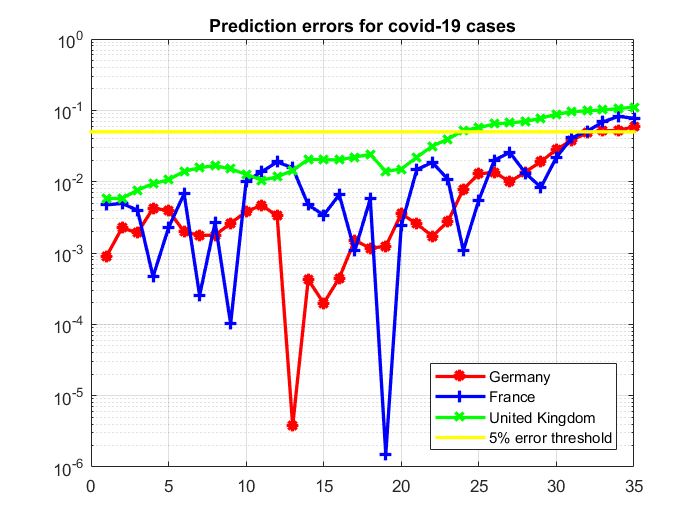}	
\includegraphics[width=0.32\linewidth]{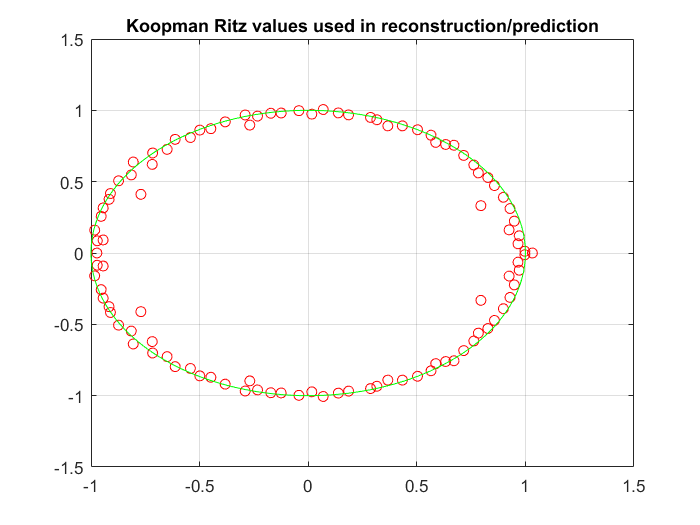}
	\includegraphics[width=0.32\linewidth]{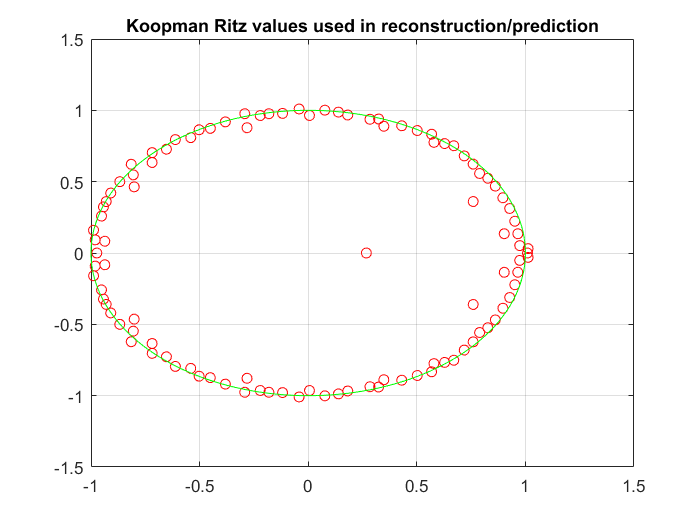}
	\includegraphics[width=0.32\linewidth]{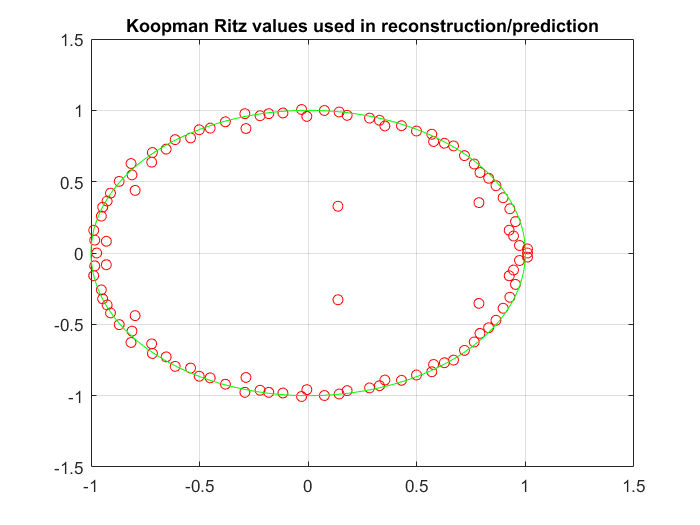}	
	\vspace{-3mm}
  \caption{\emph{First row}: Prediction error for Germany, France and the United Kingdom for a $35$ days prediction starting from the data windows $\hb_{1:105}$ (prediction for September 14 -- October 18), $\hb_{1:106}$ (prediction for September 15 -- October 19), $\hb_{1:107}$ (prediction for September 16 -- October 20), respectively. \emph{Second row}: The corresponding Koopman Ritz values used in the \textsf{KMD}.}
	\label{fig:covid_DS1_105_35}
\end{figure}

In the next test, we use the data windows $\hb_{1:132},\ldots, \hb_{1:140}$ for $28$ days predictions for the time intervals October 11 - November 7, October 12 - November, ..., October 19 - November 15. The results are shown in Figure \ref{fig:covid_DS1_132_28}.

\begin{figure}[ht]
	\centering
	\includegraphics[width=0.32\linewidth]{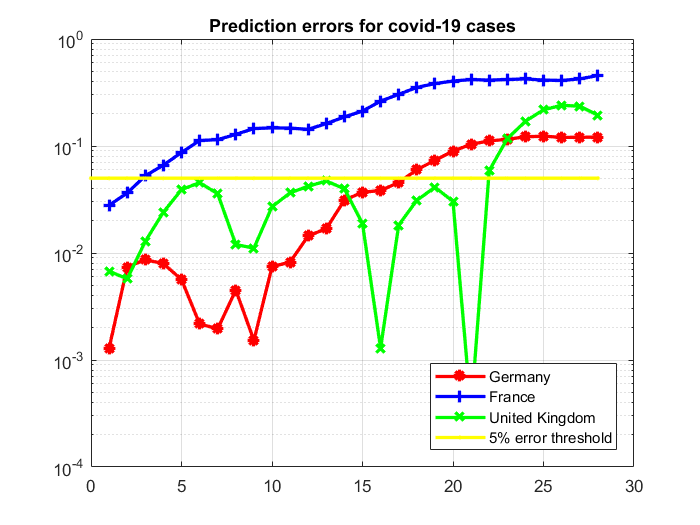}
	\includegraphics[width=0.32\linewidth]{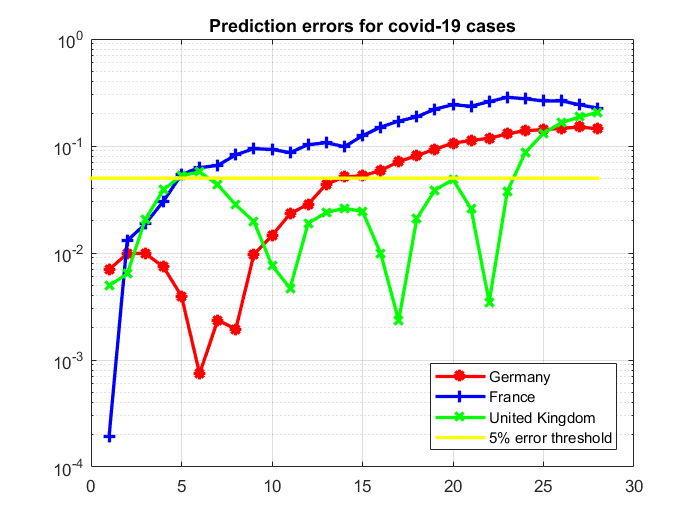}
\includegraphics[width=0.32\linewidth]{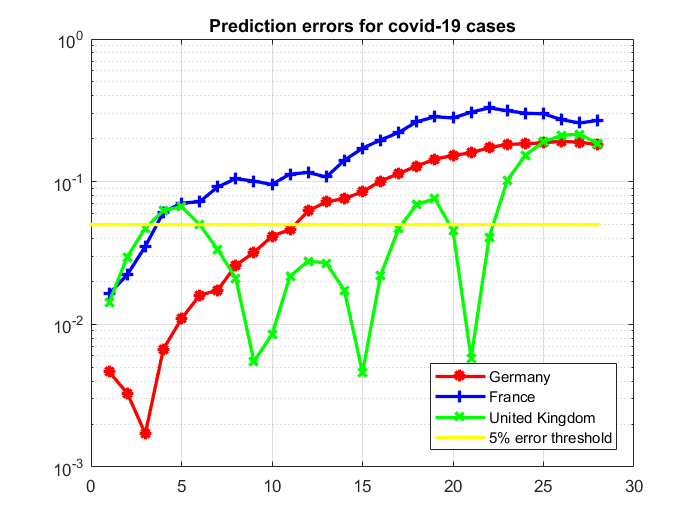}
\includegraphics[width=0.32\linewidth]{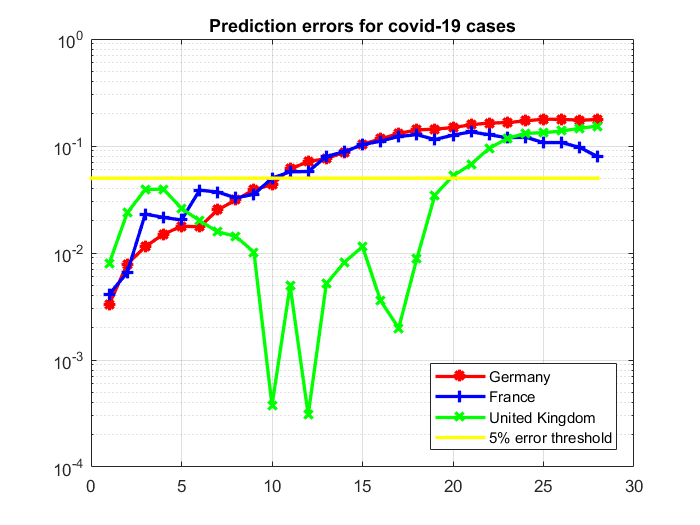}
	\includegraphics[width=0.32\linewidth]{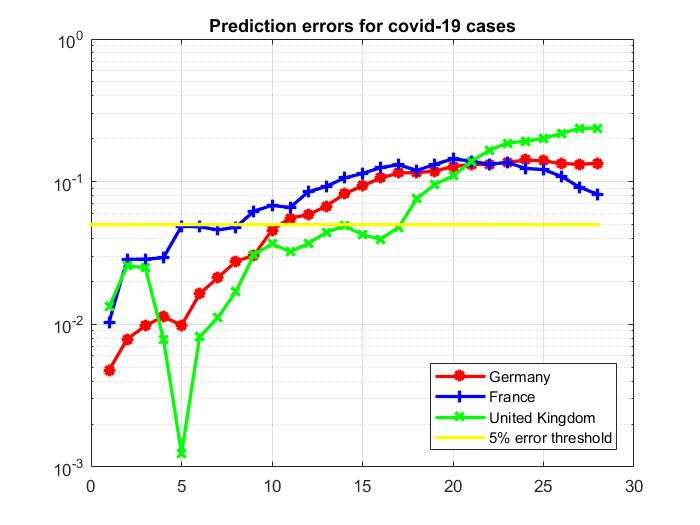}
\includegraphics[width=0.32\linewidth]{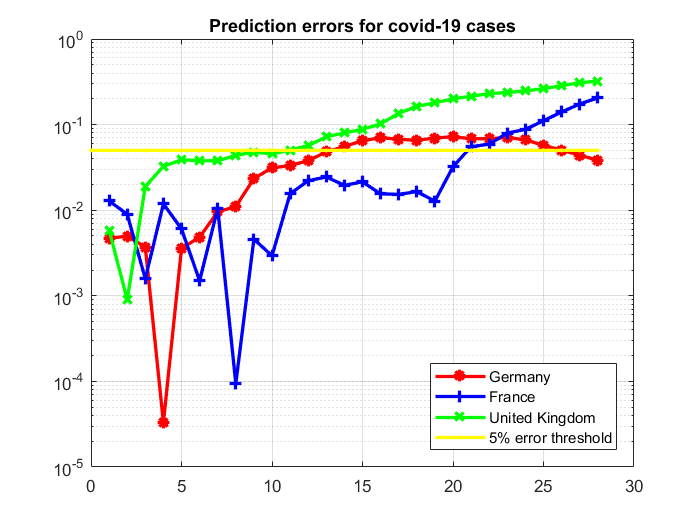}
\includegraphics[width=0.32\linewidth]{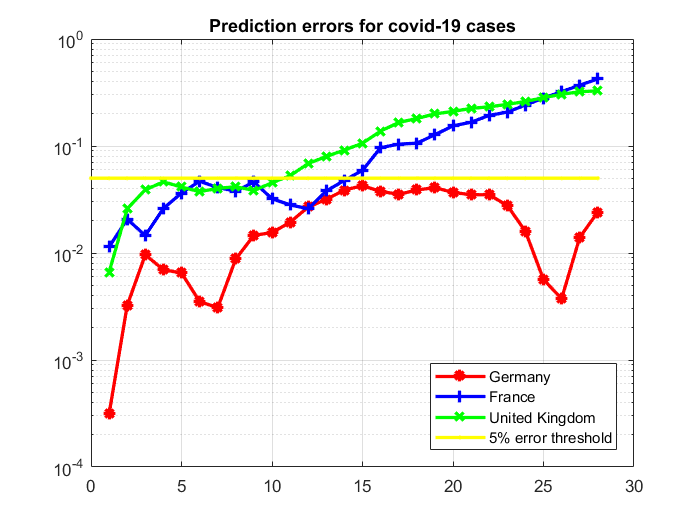}
	\includegraphics[width=0.32\linewidth]{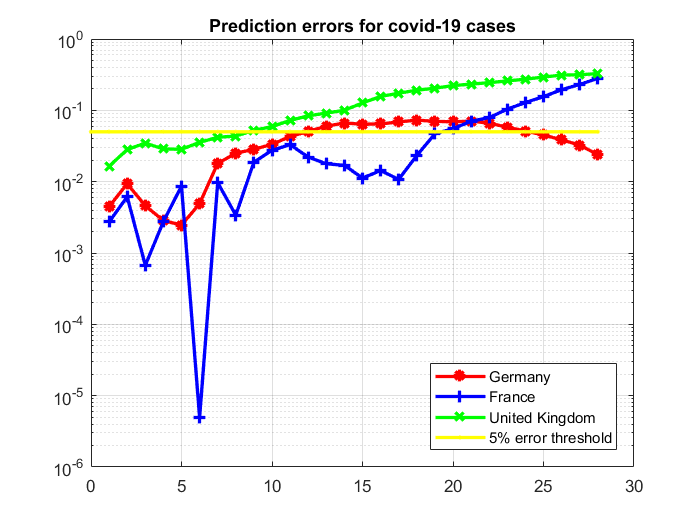}
\includegraphics[width=0.32\linewidth]{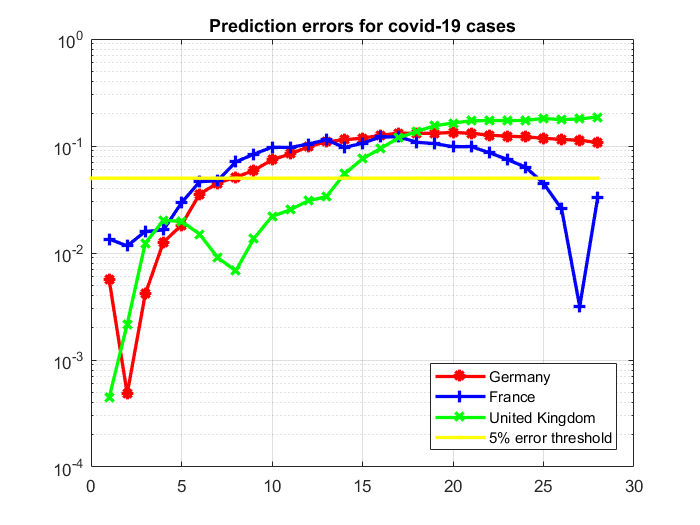}
  \caption{Prediction error for Germany, France and the United Kingdom for a $28$ days prediction,  based on the windows $\hb_{1:132},\ldots, \hb_{1:140}$, respectively. The prediction intervals are, respectively, October 11 - November 7, October 12 - November, ..., October 19 - November 15.  
  }
	\label{fig:covid_DS1_132_28}
\end{figure}

Now, we go to the datased \textbf{DS3}. The focus is on some computational details related to the two main ingredients -- the \textsf{DMD} and the \textsf{KMD}. The datased \textbf{DS3} is constructed by a single and a double application of the Savitzky-Golay filter to the Germany data, so that $d=3$. (The filter uses cubic polynomial and data window of width 5. On the left boundary, we add zero values, and on the right boundary we leave the original data. The filtered data differ from the original  at most five to ten percent relative error in the first $30$ days and at most $O(10^{-3})$ afterwards.) The purpose of the test is to create a situation that one could encounter when deploying the Koopman/\textsf{DMD} framework for data driven prediction or for a discovery and analysis of latent coherent structures. 
\begin{figure}[H]
	\centering
\includegraphics[width=0.32\linewidth]{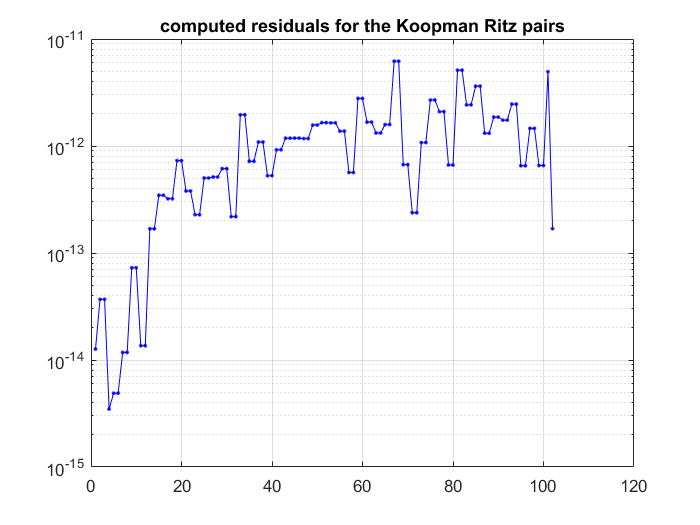}	
\includegraphics[width=0.32\linewidth]{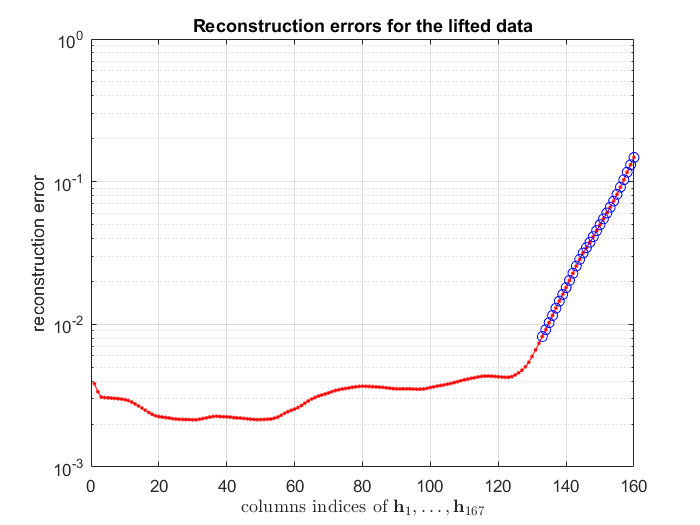}
\includegraphics[width=0.32\linewidth]{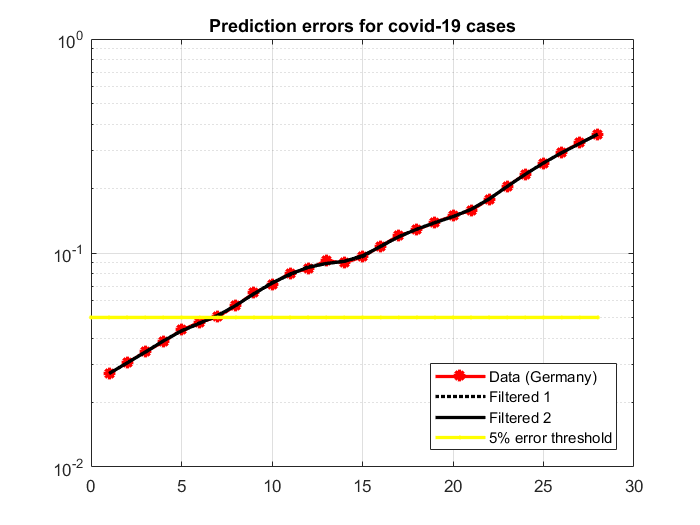}	
\vspace{-2mm}
\caption{Prediction experiment with \textbf{DS3} with data from Germany. \emph{Left panel}: the computed residuals for the computed $102$ Koopman Ritz pairs (extracted from a subspace spanned by $132$ snapshots $\hb_{1:132}$). Note that all residuals are small. The corresponding Ritz values are shown in the first panel in Figure \ref {fig:DE_132_28_W}. \emph{Middle panel}: \textsf{KMD} reconstruction error for $\hb_{1:132}$ and the error in the predicted values $\hb_{133:160}$ (encircled with $\textcolor{blue}{\circ}$). The reconstruction is based on the coefficients $(\alpha_j)_{j=1}^r=\mathrm{arg\min}_{\alpha_j}\sum_{k} \| \hb_k - \sum_{j=1}^{r} \lambda_j^{k}\alpha_j \rv_j\|_2^2$.
\emph{Right panel}: Prediction errors for the period October 11 -- November 7. Compare with the first graph in Figure \ref{fig:covid_DS1_132_28}.
}
\label{fig:DE_132_28_UW}
\end{figure}
\noindent We recall that a \textsf{DMD} algorithm uses a rank revealing decomposition with some threshold value and that the number of the computed Ritz pairs may be smaller than the column dimension of the matrix $\X$; see \S \ref{SS=Schmid-DMD}. Then the reconstruction formula (\ref{eq:alpha-pinv-V}) for the coefficients is not valid, and one has to satisfy  (\ref{eq:hi-reconstruct-appr}) by solving the least squares problem
$\sum_{k} \| \hb_k - \sum_{j=1}^{r} \lambda_j^{k}\alpha_j \rv_j,\|_2^2 \longrightarrow\min_{\alpha_j}$,
where the reconstruction error is not necessarily small, and it introduces noise into the extrapolation process outlined in \S \ref{SS=Predict-Idea-Limits}. (Recall, if we have full set of modes, then the reconstruction is perfect and the only error is from the finite precision arithmetic.) The prediction skill based on this \textsf{KMD} is shown in Figure \ref{fig:DE_132_28_UW}. 
Now we change the reconstruction strategy and state the problem as the 
weighted least squares problem (\ref{eq:SP-DMD}) with the weights that favour the four most recent snapshots with weights set to one, and with the weights of all other snapshots set to the machine round-off unit $\bfvareps\approx 2.2\cdot 10^{-16}$.  
\begin{figure}[ht]
	\centering
\includegraphics[width=0.32\linewidth]{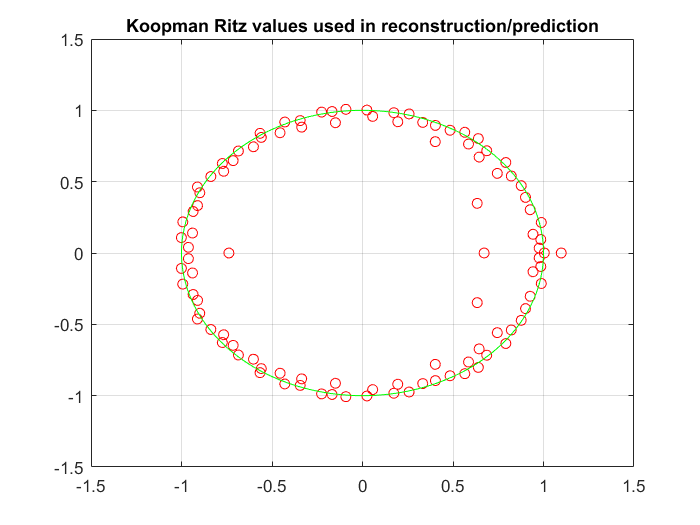}	
\includegraphics[width=0.32\linewidth]{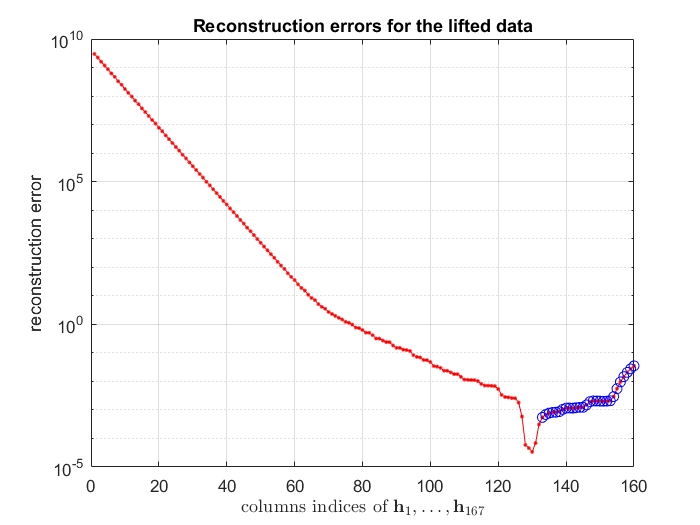}
\includegraphics[width=0.32\linewidth]{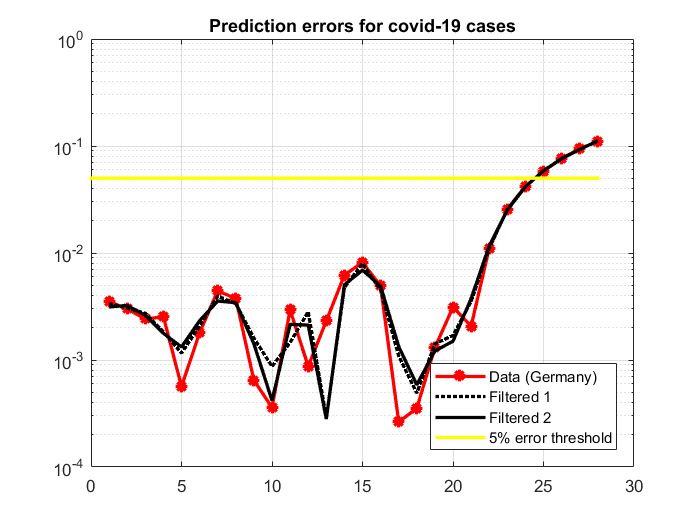}	
\caption{Prediction experiment with \textbf{DS3} with data from Germany. \emph{Left panel}: the computed $102$ Koopman Ritz values (extracted from a subspace spanned by $132$ snapshots $\hb_{1:132}$). The corresponding residuals are shown in the first panel in Figure \ref{fig:DE_132_28_UW}. \emph{Middle panel}: \textsf{KMD} reconstruction error for $\hb_{1:132}$ and the error in the predicted values $\hb_{133:160}$ (encircled with $\textcolor{blue}{\circ}$). The reconstruction is based on the coefficients $(\alpha_j)_{j=1}^r=\mathrm{arg\min}_{\alpha_j}\sum_{k} w_k^2 \| \hb_k - \sum_{j=1}^{r} \lambda_j^{k}\alpha_j \rv_j\|_2^2$.
\emph{Right panel}: Prediction errors for the period October 11 -- November 7. Compare with the first graph in Figure \ref{fig:covid_DS1_132_28}, and with the third graph in Figure \ref{fig:DE_132_28_UW}.
}
\label{fig:DE_132_28_W}
\end{figure}
\noindent The effect of weighting is best seen by comparing the middle graphs in Figures \ref{fig:DE_132_28_UW} and \ref{fig:DE_132_28_W}. In the case of weighting, the reconstruction of almost all leading snapshots is bad, but the ones more important for the prediction task have much smaller error. 

}
\subsection{Comments on SIR type models}

The key coefficient in SIR-type models, the so-called reproduction number $R_0$ can  be estimated using Koopman operator techniques. Namely, the classic SIR model reads
\begin{eqnarray}
\dot s&=&-\beta si \nonumber \\
\dot i&=&\beta si -\nu i\nonumber \\
\dot r&=&\nu i \nonumber \\
\end{eqnarray}
Under condition $s=1$ (infinite reservoar of susceptibles), the exponential growth happens when $\beta-\nu>0$, i.e.
$\beta/\nu>1$. The reproduction number is defined by $R_0=\beta/\nu$. Thus, $R_0$ is related to the coefficient of exponential growth. Since we know that eigenvalues of the linearized system are eigenvalues of the Koopman operator, the largest real Koopman eigenvalue is related to $R_0$.
Another number commonly estimated for use in tracking of epidemics is the instantaneous reproduction number defined by \cite{gostic2020practical}
\[
R_t^{inst}=\beta s D,
\]
where $D$ is the duration of the infectiousness.
For small $\nu$, and constant $s=1$, we have
\begin{equation}
   R_t^{inst} \approx \beta D,
\end{equation}
providing another connection between the SIR models, Koopman operator spectrum and Reproduction numbers.
\subsection{The framework for prediction and Black Swan event detection}\label{SS-BS-detect}
In \textsf{GKP} algorithm, the spectral information is extracted from a sequence of active windows -- for each window, the snapshots are  arranged in a  Hankel matrix whose columns define a new set of snapshots and  approximate eigenvalues and eigenvectors are computed using Algorithm \ref{zd:ALG:DMD:RRR}.
%
In the case when the dynamics of the system is not coupled with some other dynamical system, we expect that, in the absence of unexpected disturbances, the \textsf{AKMD} will capture at least the basic trends of the dynamics. In particular, the spectral radius of the active window (the maximal absolute value of the selected Ritz values) should not change too much. Further, the \textsf{DMD} algorithm should compute Ritz pairs with reasonably small residuals. This is plausible, because the sequence (\ref{eq:seq-Uif}) can be, at any moment, interpreted as an excerpt from a power method generated sequence, and the power method in the limit reveals the absolutely dominant eigenvalues.

However, if the dynamical system data are hit by disturbance, this could be recognized e.g. by detecting the active windows whose spectral radii  change significantly, or by the absence of Ritz pairs with small residuals (see Figures \ref{Fig:L:1b}, \ref{Fig:L:1c}, \ref{fig:Lorenz-global-predict} ). This enables us to pinpoint the discrete time moments/subintervals at which  disturbances interfere with the original dynamics.  
For the chosen reference interval $\mathcal{I}$, if $\max_j |\lambda_j| \not\in\mathcal{I}$ or if there are no Ritz pairs with reasonably small residuals, we flag the observed active window as the window  which possibly contains a Black Swan event. By sliding the active windows along the computational domain, using the flagged windows, we determine the time sub-intervals containing  disturbances whose dynamics is not well captured by the corresponding \textsf{AKMD} models. 

{The reference interval $\mathcal{I}$ can be determined (and dynamically adjusted) e.g. by first computing the Ritz values for many active windows, and then by trial and error, including a statistical reasoning and information theoretic  techniques (see  e.g. \cite{metzner2012}) learn to differentiate between the acceptable interval $\mathcal{I}$ for spectral radii and the values that are considered outliers. This is best done on a case-by-case basis.} 

%
This scheme can be implemented with different sizes of the Hankel matrices (see \S \ref{SS=Lift+Hankel}) and with different sizes of active windows and then the Black Swan event intervals can be determined by taking into account all determined intervals.


\subsubsection{The retouching trick to process Black Swan events}
If the Black Swan event data are included in the training set, the dynamics of the original system (decoupled from this disturbance) cannot be revealed, which means that the prediction of the dynamics after the Black Swan event will be damaged, if not impossible. 
{However, instead of using the original data we can replace them with the data obtained by the prediction based on the information from the previous active windows, preceding flagged intervals.}
This is illustrated in (\ref{eq:bs-bypass}): the perturbed value $(\Koop^{k}\bff)(\z_0) + \bfvarepsilon_k$ is replaced with $\widehat{(\Koop^{k}\bff)(\z_0)}$, which is a predicted value based on the previous undisturbed data. The same can be done for the remaining data in the flagged window.
\begin{equation}\label{eq:bs-bypass}
    \begin{array}{l|c|c|c|c|c}
\mbox{undisturbed} &     \ldots   &  (\Koop^{k-1}\bff)(\z_0) & (\Koop^{k}\bff)(\z_0) & (\Koop^{k+1}\bff)(\z_0) & \ldots\\\hline
{\mbox{disturbance at $k$, $k+1$}} & \ldots   &  (\Koop^{k-1}\bff)(\z_0) & {(\Koop^{k}\bff)(\z_0) + \bfvarepsilon_k}& {(\Koop^{k+1}\bff)(\z_0) + \Koop\bfvarepsilon_k + \bfvarepsilon_{k+1}} & \ldots\\\hline
{\mbox{use prediction at $k$, $k+1$}} & \ldots   &  (\Koop^{k-1}\bff)(\z_0) & {\widehat{(\Koop^{k}\bff)(\z_0)}} & {\widehat{(\Koop^{k+1}\bff)(\z_0)}} & \ldots
    \end{array}
\end{equation}
The prediction after the Black Swan event then becomes more stable and in most cases quite successfully predicts data after the Black Swan event; see \S 2.2 and Fig. 2 in the main text. 

There are many variations of this scheme. For instance, it could happen that the Black Swan interval detected in the algorithm is too long and possibly unrealistic. Therefore we limit the length of the interval on which the data are replaced in order to prevent the algorithm from changing the dynamics too much. Then we apply  the algorithm again and detect if the replacements result with decreasing of maximum of the absolute value of the eigenvalues over the active windows. The whole process can be repeated more times to remove eventual  Black Swan events that are not taken into account in the previous steps. Finally, the retouched data, cleaned from the Black Swan event disturbances, are used for the prediction.  

\begin{algorithm}[ht]\label{alg:global}
	\caption{\bf (Global Koopman Prediction (\textsf{GKP}) with Black Swan event detection and switching to local prediction)}	
	
	\begin{algorithmic}[1]
		\Require{$\bullet$ Data snapshots $\bff_0,\bff_1,...,\bff_{end} $; $\bullet$  the size of the training data window $\win=n_H+m_H$; $\bullet$ the dimensions $n_{H}$, $m_{H}$ ($n_{H}> m_{H}$) of the Hankel matrices; $\bullet$ the threshold $\eta > 0$ for the maximal acceptable value of the residual of Ritz vectors; $\bullet$ the sliding step $\Delta p$; $\bullet$ the maximal number $N_{rep}$ of iterative retouching of the perturbed data; $\bullet$ the maximal time length $L_{BS}$ for one step replacement of the Black Swan event data with the predicted values; $\bullet$ the reference interval $\mathcal{I}$ for spectral radius used for the detection of Black Swan event moments; $\bullet$ lead times $\tau_g$ and $\tau_l$ for the global and local prediction.}
		\Ensure{Predicted system observables $\widetilde{\bff}_{n_{H}+m_{H}},\widetilde{\bff}_{n_{H}+m_{H}+1}, \ldots, \widetilde{\bff}_{\win+n_{dmd}\cdot\Delta p+\tau_g}$ (or $\widetilde{\bff}_{\win+n_{dmd}\cdot\Delta p+\tau_l})$}
		\State{$j_{rep}=0$; $n_{dmd} = \lfloor\frac{end-\win}{\Delta p}\rfloor$; $BS_{event} = False$, $n_{BS}=0$}	
		\While{($j_{rep} < N_{rep}$) and ($j_{rep}=0$ or $n_{BS} \neq 0$)}
        \State{$n_{BS}=0$}
		\For {$p=\win,\win+\Delta p,...,\win+n_{dmd}\cdot\Delta p$} 
		\State{For the active window $\mathcal{W}(p,\win)$ apply KMD algorithm to obtain AKMD   
	  	using $n_{H} \times m_{H}$ Hankel matrices.}
		\State{If there are no Ritz values for which the associated residual is smaller than $\eta$, set $\max_j |\lambda_j| = \infty.$}
		\If{$\max_j |\lambda_j| \notin \mathcal{I}$ }
		    \If{$BS_{event}=False$}
		    \State{Set $t_{BSbegin} = t_{\max(0,p-\Delta p)}$; $BS_{event}=True$; $n_{BS} = n_{BS} + 1 $ \Comment {\emph{New disturbance appears.}}}
		    \EndIf
		    \State{Flag the time interval $[t_{p-\Delta p},t_p]$ as a Black Swan event interval; }
			\State{In the interval $[t_{p-\Delta p},t_p]$ use Local Koopman Prediction with lead time $\tau_l$ (Algorithm \ref{alg:local}) } 
			\State\Comment{\emph{Remark: Local Koopman Prediction algorithm (Algorithm \ref{alg:local}) can be applied on the whole domain and then associated local prediction is used on the detected critical intervals.}}
		    \State{Store the data from last active window not including Black Swan event for retouching the data in $[t_{p-\Delta p},t_p]$ using the prediction obtained by global AKMD .}
		\Else    
		    \If{$BS_{event} = True$} 
		        \State{Set $t_{BSend} = t_{p-\Delta p}$; $BS_{event} = False$  \Comment {\emph{End of the Black Swan event.}}}
		        \State{Replace the original data in the Black Swan event interval $[t_{BSbegin}, \min(t_{BSbegin}+L_{BS}, t_{BSend}]$ with the stored retouched data.}
		     \Else
		     \State{Using the \textsf{AKMD} associated with $\mathcal{W}(p,\win)$ and (\ref{eq:prediction}), extrapolate to obtain the predictions  $\widetilde{\bff}_{p+\tau}$, $\tau \le \tau_g$.}
		     \EndIf
		\EndIf
		\EndFor
		\State{$j_{rep}= j_{rep}+1$}
		\EndWhile
	\end{algorithmic} 
\end{algorithm}
\subsubsection{Local Koopman prediction}\label{SS=Local-prediction}
{In some cases, the global prediction algorithm is not feasible. For instance, when we  just start collecting the data, we have not enough information for a \textsf{GKP} analysis. Or, in the situation when \textsf{GKP} recognizes the beginning of a Black Swan event, as discussed at the beginning of \S \ref{SS-BS-detect}. Then, the available data cannot be used for prediction, because the dynamical system has changed. The new model must be built from scratch, as if we just started getting new data. The best we can do is to create a new local algorithm that needs less data, but also with a much shorter reach into the future. } 

In the Local Koopman Prediction (\textsf{LKP}) algorithm we change the size of the active window depending on the success of the previous prediction. 
The idea is to assimilate as much acquired data as possible, so we set Hankel matrix dimension variable with prediction moment, i.e. $n_H=n_H(p)$ and $m_H=m_H(p)$. 
We also choose the minimal Hankel matrix dimension
\begin{equation}\label{eq:Hankel-mindim}
	(n_{H,min}\cdot d) \times(m_{H,min}+1), 
\end{equation}
and we start predictions with such minimal Hankel matrix i.e. for first prediction $p=p_0$ we set
\begin{equation}\label{eq:Hankel-restart}
	n_H(p_0)=n_{H,min},\;\; m_H(p_0)=m_{H,min}.  
\end{equation}
When data $\bff_p$ at prediction time $t_p$ becomes available, we can compute the error of the prediction $\widetilde{\bff}_p$, using suitable norm, as 
\begin{equation}\label{eq:relerr}
	\epsilon_p=\|\widetilde{\bff}_p - \bff_p\| / \|\bff_p\| .
\end{equation}
At other prediction moments, if the prediction error (\ref{eq:relerr}) is smaller than the referent error $\epsilon_{ref}$ we assimilate the newly acquired data into the active window by increasing the Hankel matrix size
\begin{equation}\label{eq:Hankel-resize}
	n_H(p)=n_H(p-1)+1 \mbox{, or } m_H(p)=m_H(p-1)+1.
\end{equation} 
Otherwise, i.e. if the prediction error (\ref{eq:relerr}) is larger then the referent one, we reset the Hankel matrix dimension to the minimal one: 
$n_H(p)=n_{H,min}$, $m_H(p)=m_{H,min}$.  

In both cases we recompute the Hankel matrix and the \textsf{AKMD}
for each new prediction.  

\begin{algorithm}
	\caption{\bf (Local Koopman Prediction (\textsf{LKP}) with resizing Hankel matrix) \label{alg:local}}	
	
	\begin{algorithmic}[1]
		\Require{$\bullet$ Data snapshots $\bff_0,\bff_1, \ldots, \bff_{end}$; $\bullet$ indices of time moments for the begin and the end of the local prediction $k_0, k_f$  (optionally) $\bullet$ minimal Hankel matrix dimension $n_{H,min}$, $m_{H,min}$;  $\bullet$ error threshold $\epsilon_{ref}$; $\bullet$ lead time $\tau_l$. }
		\Ensure{Predicted system observables $\widetilde{\bff}_{n_{H,min}+m_{H,min}},\widetilde{\bff}_{n_{H,min}+m_{H,min}+1}, \ldots, \widetilde{\bff}_{end+\tau_l}$ (or $\widetilde{\bff}_{k_0}$, \ldots, $\widetilde{\bff}_{k_f+\tau_l}$)}
		\If{$k_0$ and $k_f$ not defined}
		\State{$k_0=n_{H,min}+m_{H,min}, k_f=end$}
		\EndIf
		\For {$p=k_0,k_0+1,...,k_f$ let $t_{p-1}$ be the time of the last known data} 
		\If {the error $\epsilon_p$ (\ref{eq:relerr}) is larger than referent error $\epsilon_{ref}$}
		\State {Resize the Hankel matrix to the minimal size (\ref{eq:Hankel-restart}).}
		\Else
		\State {Increase the size of the Hankel matrix using (\ref{eq:Hankel-resize}).} 
		\EndIf
		\State {Form the Hankel matrix (\ref{eq:Hankel}) and compute the \textsf{AKMD}}.  
		\State {Using the \textsf{AKMD}, extrapolate to obtain the prediction (\ref{eq:prediction}).} 
		\EndFor
	\end{algorithmic} 
\end{algorithm}

\section{Supplementary material -- Discussion.}
In this supplementary material we validate our approach by three case studies.
In \S \ref{S=Lorenz}, we use the Lorenz system to illustrate the main idea of monitoring the Koopman Ritz values and the prediction skills of the proposed method. The model free aspect of the proposed method is  further illustrated in applications in two entirely different settings: physiological processes (resonant breathing) in Supplementary section \ref{SS=Res-Breath} and geomagnetic substorms (prediction of the \textsf{AL} index) in Supplementary section \ref{S=GeoMag}. 
Finally, in Supplementary section \ref{SS=Influenza_more} we provide additional numerical results related to the influenza epidemics studied in section 2.2 of the paper.

\subsection{Case study: Lorenz system}\label{S=Lorenz}
\begin{figure}[ht]
\begin{subfigure}{.35\textwidth}
	\centering
	\includegraphics[width=\linewidth,height=1.8in]{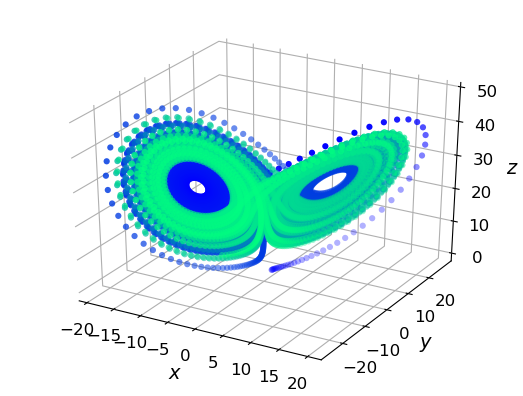}
	\vspace{-5mm}
	\caption{\emph{Lorenz system (\ref{eq:LorenzSystem}).}}
	\label{Fig:L:1a}
\end{subfigure}
\begin{subfigure}{.32\textwidth}
	\includegraphics[width=\linewidth,height=1.8in]{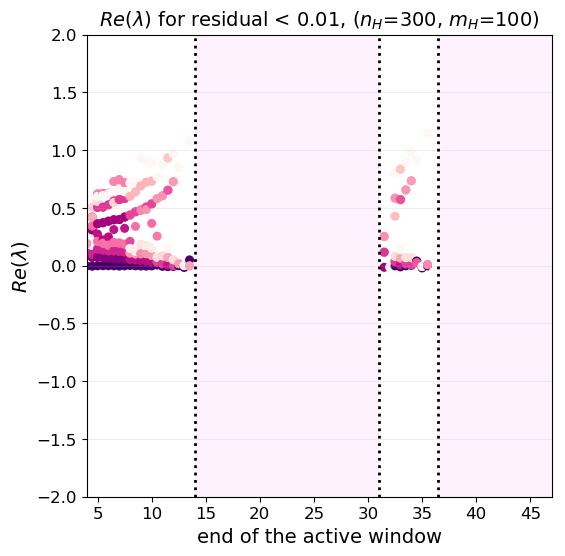}
	\vspace{-5mm}
	\caption{ $\Re(\lambda)$}
	\label{Fig:L:1b}
\end{subfigure}
\begin{subfigure}{.32\textwidth}
\includegraphics[width=\linewidth,height=1.8in]{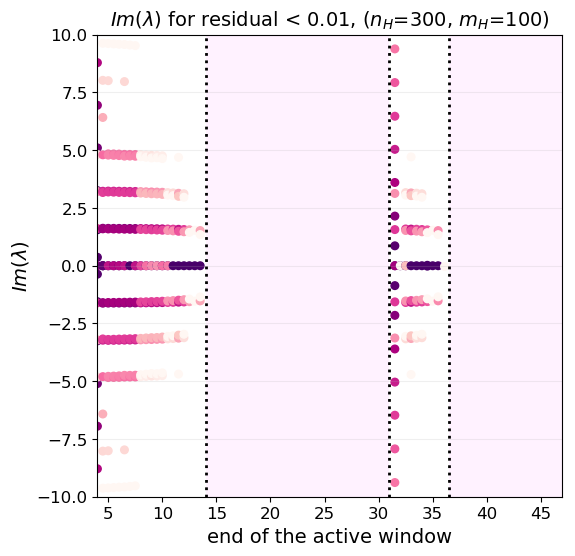}
\vspace{-5mm}
\caption{$\Im(\lambda)$}
\label{Fig:L:1c}
\end{subfigure}
\begin{subfigure}{\textwidth}
\includegraphics[width=\linewidth, height=2.0in]{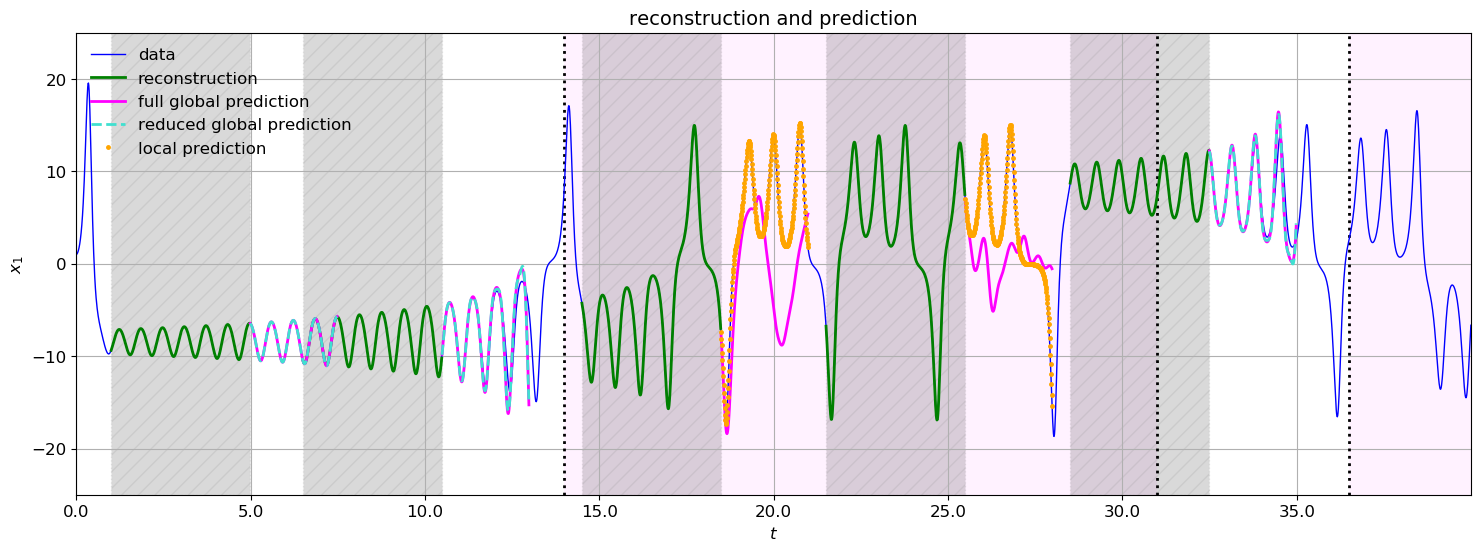}
\vspace{-7mm}
\caption{Global and local prediction for $2.5 s$ ahead.\label{fig:Lorenz-global-predict}}
\end{subfigure}
\vspace{-3mm}
	\caption{
(\ref{Fig:L:1a}): the Lorenz system (\ref{eq:LorenzSystem}). For each active window, the Ritz pairs with the residual below $\eta_r=0.01$ are selected (see \cite[\S 3.2]{DDMD-RRR-2018}); the real and the imaginary parts of the corresponding Ritz values are shown in (\ref{Fig:L:1b}, \ref{Fig:L:1c}).
The color intensity of the eigenvalues indicates  the amplitudes of the corresponding modes.	
	(\ref{fig:Lorenz-global-predict}): \textsf{KMD} reconstruction and prediction of the observable {$x_1 \equiv x$} for the Lorenz system (\ref{eq:LorenzSystem}). The data are collected in five active windows (time intervals $[1,5]$, $[6.5,10.5]$, $[14.5,18.5]$, $[21.5,25.5]$, $[28.5,32.5]$, marked by shadowed rectangles) and then the dynamics is predicted for the time moments ahead of the active window.  Note - by comparing the positions of the intervals with poor prediction with the eigenvalue-free pink rectangles in Figures \ref{Fig:L:1b}, \ref{Fig:L:1c} - that the failure of the global prediction occurs after the active windows which do not contain Ritz pairs with sufficiently small residuals, as indicated by magenta curves.    The local algorithm recovers the prediction capability, {using a sequence of shorter moving \textsf{KMD}'s } {and prediction for $10$ time steps ahead, as indicated by orange curves}.
}	
\label{Fig:L}
\end{figure}

The critical underlying concept in chaotic dynamics is that of sensitivity to initial conditions and the associated positivity of Lyapunov exponents that measure long term exponential deviation of nearby trajectories \cite{lorenz1963deterministic}. Namely, the long term exponential divergence of nearby trajectories leads to unpredictability due to the finite nature of (any) prediction algorithm. Even the implementation of exact equations of a dynamical system on any computing device leads to finite precision calculations and therefore ultimate exponential divergence of prediction from true trajectory. However, this neglects the finer aspects of chaotic dynamics that are exhibited in the 
most paradigmatic of chaotic systems -- the Lorenz dynamical systems, 
 modeled by Lorenz equations
\begin{equation}\label{eq:LorenzSystem} 
\dot\x = \begin{pmatrix}
\dot{x} \cr
\dot{y} \cr
\dot{z}
\end{pmatrix}  = \begin{pmatrix}
\sigma (y-x) \cr
x (\rho-z) -y \cr 
xy-\beta z 
\end{pmatrix}
\end{equation} 
\noindent with $\sigma=10$, $\beta = 8/3$, and $\rho=28$ for which the system exhibits chaotic behavior.
For understanding of the prediction capability for the Lorenz system, more important than the long term exponential divergence of trajectories is the short term divergence typically induced by switching between the two \emph{wings} of the butterfly (see Figure \ref{Fig:L:1a}).

At the core of our approach is the observation that, while inside one of the butterfly wings, the system behaves in a predictable manner. The exponential divergence is ultimately due to switching between the two butterfly wings.
The first time such a switch happens, the situation resembles a {Black Swan} event \cite{Taleb-2007} (although there is an ontological difference highlighted in the main text): the trajectory suddenly wonders off to a different part of the state space and starts exploring there, until another switch happens taking it back to the known part of the state space. This fits our paradigm of splitting the state space into domains over which prediction is possible and monitoring for the switch between such domains.   \\
\indent The current theory is thus an extension of the ideas in \cite{mezic:2005}, where deterministic components of stochastic dynamical systems were extracted using Koopman operator methods,  and \cite{brunton:2017} where it was shown that Lorenz system can be described well by a set of linear evolution equations driven by stochastic term that induces switching. In both of these, the detection of the switching moment and the precise interaction of local and global behavior on subdomains of state space was not accounted for; {we address it here}.\\
\indent We use (\ref{eq:LorenzSystem}) to  test the prediction potential of the \textsf{KMD}. 
%
We have generated data using numerical simulation (\textsf{ODE} solver) of (\ref{eq:LorenzSystem}) with  the time resolution $\delta t=0.01 s$, thus obtaining a discrete dynamical system. For a present moment (index) $\now$, an \emph{active window} of length $\win$ is selected as in \S \ref{SS=Lift+Hankel}, with $b=\now-\win$, and
the selected data are lifted in the Hankel structure.   The \textsf{KMD} is computed for the corresponding vector valued observables $\bfh_i$, and used for their prediction  as explained in \S 1 of the paper.
%
For computing the \textsf{KMD} {for the global prediction algorithm}, the active windows of size $\win=400$ and $300 \times 100$ Hankel matrices are used. By sliding the active windows along the computational domain we get prediction at different times. 
\captionsetup[subfigure]{position=b}
{When the actual data and the prediction errors become available, we either continue forecasting with the same \textsf{KMD}, or a switching device invokes the local prediction scheme with $21\times 11$ Hankel matrices if the error is above a preset threshold.} {The prediction is then with shorter forecast lead time, and predicted data are based on a sequence of local \textsf{KMD}'s.} The local algorithm keeps increasing the active windows {and the lead time}, whilst monitoring the error; see the \textbf{Methods} section.\\
\indent In Figure \ref{Fig:L},  we show the \textsf{KMD} reconstruction  and prediction results of the observable $x_1=x$ of the system (\ref{eq:LorenzSystem}) for a selection of five active windows. While the reconstruction -indicated by green traces -  works well (as expected, see e.g. \cite{Jovanovic:2014ft}, \cite{DMM-2019-DDKSVC-DFT}, \cite{LS-Vandermonde-Khatri-Rao-2018}), the prediction capability -indicated by magenta traces -  is lost for the third and the fourth active windows. An inspection of the quality of the approximate Ritz pairs computed by the \textsf{DMD} \cite{DDMD-RRR-2018} shows that for the time interval containing those two windows none of the computed pairs has the residual below $0.01$, i.e. no useful spectral information, which is essential for the \textsf{KMD},  could be extracted from the available data snapshots, see Figures \ref{Fig:L:1b}, \ref{Fig:L:1c}. As a consequence, the prediction using  numerical realization of \textsf{KMD}  cannot give satisfactory results. \\
\indent On the other hand, in that part of the domain where the trajectory behaves chaotically, and the intensive change of the nature of the eigenpairs precludes accurate numerical approximations, local prediction scheme quickly adapts to the new data, forgets the previously acquired information, and delivers better results. See Figure \ref{fig:Lorenz-global-predict}. \\
\indent The reconstruction with a reduced number of modes (see Figure \ref{fig:Lorenz-global-predict}) uses only the Ritz pairs with small residuals (see \cite[\S 3.2]{DDMD-RRR-2018}). {The number of modes used for prediction after the  first, the second and the fifth active window (gray rectangles) were 10, 18 and 10, respectively.}  
One can observe that the reconstruction and prediction capabilities - shown in dashed green - are comparable with using the full \textsf{KMD}.  
%
\subsubsection{Remark}
Regarding the question of detecting the switching moments, Figure \ref{fig:Lorenz-global-predict}  provides  an insight. If we look at the first pink zone with no ``good" eigenvalues (Figures \ref{Fig:L:1b}, \ref{Fig:L:1c}), we see that it starts close to the switching moment. Also, this zone is quite long because the switching moments in that zone are too close to each other and no learning data window can fit in-between. Only when two switching moments are distant enough, we can find ``good" eigenvalues, the learning data window exits the pink zone, and the global prediction recovers. {{It is remarkable that recent works \cite{kordaetal:2020,giannakisetal:2018} have found spectral objects - pseudoeigenfunctions - that govern quite regular short term dynamics inside  the wing of the Lorenz attractor. This dynamical feature - discovered by careful analyses of the continuous Koopman operator spectrum for the Lorenz system - seems to enable the prediction algorithm performance.}}
\subsubsection{Discussion}
{Historically, the most discussed way in which a substantial change in dynamics can occur in dynamical systems is due to a change of a value of a bifurcation parameter \cite{iooss2012elementary}.
The prediction method that we propose is not necessarily related to a change of parameter in the system. Namely, the original description of the black swan event does not relate to a change in parameter, just to travel to another part of the space (here considering the dynamical system to be the ecological system). White swans were known to exist in Europe, but explorers found black swans in Australia. The prediction that an explorer would make when traveling to Australia might have been existence of a white swan. Upon observation, they concluded that black swans exist. The bird had all other  properties of the white swan, except for the color. There were no parameter changes, no bifurcation that occurred. Similarly, the prediction of the dynamics while on one wing of the Lorenz butterfly attractor is based on the eigenvalues of the Koopman operator detected   while sampling that wing. Once the dynamics "travels" to the other wing, the change in dynamics is recognized (although there are no parameter changes), but as the dynamics continues on the other wing, the same eigenvalues are obtained. The difference is in the resulting local \cite{Mezic:2019} eigenfunctions, (or pseudoeigenfunctions, as in \cite{kordaetal:2020,giannakisetal:2018}), that are related by the symmetry $(x,y)\rightarrow (-x,-y)$. 

As is well known, chaotic dynamics is an asymptotic property of a dynamical system, and the associated unpredictability is not due to local passage near saddles, but to long term repeat of such events, that ultimately leads to mixing dynamics \cite{luzzatto2005lorenz}. Switch in dynamics is here due to internal effects, and thus ontologically different from the Black Swan situation. The switch is due - in the Lorenz case -  precisely to the local saddle event, that transitions the dynamics from one wing of the butterfly to the other. We presented a method by which such passage can actually be detected, and accounted for, inside a prediction algorithm.

We note there are methods of prediction of chaotic dynamical systems that can predict the evolution over several Lyapunov times of chaotic systems \cite{pathak2017using}. 
 
Note that our purpose is somewhat different than in \cite{pathak2017using}, We are more interested in detection of failure to predict accurately, then establishing a method for long-term (climate) prediction. In separate work \cite{SysId-KIG-2021} we pursue the question of long-term prediction of the Lorenz model and provide evidence of ability of Koopman based methods to predict over many Lyapunov time-scales.
}    
\subsection{Case study: Resonant breathing}\label{SS=Res-Breath}
The mathematical model of the human cardiovascular system was developed by Ursino and Magosso in (\cite{magosso2002cardiovascular}, \cite{ursino1998interaction}, \cite{ursino2000acute}, \cite{ursino2003role}). This model includes mathematical descriptions of a pulsating heart, as well as the mechanics of blood flow (\cite{stefanovska1999physics}) and baroreflex activation (\cite{magosso2002cardiovascular}, \cite{ursino1998interaction}, \cite{ursino2000acute}, \cite{ursino2003role}). 
It includes more than $90$ parameters and $21$ states (pressures, flows, volumes, resistances, and elastances). Twenty-one delay differential equations reflect conservation of mass and balance of forces at arteries and veins, as well as delayed physiological responses to vagal and sympathetic neural activity. This allows for simulation of high-resolution blood pressure and heart period as a function of time. In  \cite{fonoberova2014computational} we modified the Ursino and Magosso model to use experimentally derived respiration period as a model input. In addition, we set external noise from the Ursino and Magosso model to zero, because of the noise in the respiratory input used in our model.

Data for model validation were provided by $12$ men and $12$ women who were healthy college students between $21$ and $24$ year of age. They were participants in an experiment one of the aims of which was to develop a computational physiology approach to model how cardiovascular processes change when the baroreflex mechanism is challenged. This study was approved by the Rutgers University Institutional Review Board for the protection of human subjects involved in research. One of the tasks that the participants completed, was a $5$-minutes resonance breathing task (6P) (\cite{lehrer2003heart}, \cite{lin2012heart}, \cite{vaschillo2002heart}, \cite{vaschillo2012measurement}), during which they breathed at a rate of approximately $6$ breaths/min following a visual pacer (Easy Air, Biofeedback Foundation of Europe, Montreal, Canada). The specific details on the participants’ selection/exclusion process and experimental procedure can be found in \cite{fonoberova2014computational}.

In \cite{fonoberova2014computational}, to find the optimal set of parameters for each subject, we selected as model output the cost function that takes into account power spectral densities and time averages of several observables, such as heart period. Instead of doing brute-force optimization on the cost function with over $90$ parameters in the model, we used the following procedure: 1) an initial sensitivity analysis was performed to select the most important parameters to tune, and 2) optimization of these most important parameters was performed to minimize the cost function. The details of the calibration procedure can be found in \cite{fonoberova2014computational}.

We use the results of the chosen simulation to analyze the global prediction algorithm on it. The numerical solutions were obtained using AIMdyn's \textsf{GOSUMD} software. The used time step for numerical simulations was $\Delta t_0 = 0.003$. The parameters in the simulations, with the exception of the function modeling breathing, were chosen as obtained in \cite{fonoberova2014computational}. The input breathing function in this simulation was chosen so that in the first part of simulation, the period of input breathing function was constant and equal to $10$ seconds. The period of $10$ second simulates the rhythm of resonant breathing. In the second half of the simulation, experimentally determined normal breathing function was used as an input. 

As already mentioned, in the global prediction algorithm one should provide long enough set of data in order to extract the dynamical system parameters related to the phenomena we want to capture with the algorithm. On the other hand, the time step between the neighboring snapshots should be chosen so that the balance between the numerical complexity and of the length of dynamical phenomena one want to reveal with used \textsf{KMD} algorithm is achieved.   

Since in the system there are no frequencies larger than $10$, it is enough to take $\Delta t =0.03$ between neighboring snapshots. By using time-lagged embedding for each variable separately, we form the Hankel matrices and apply the \textsf{GKP} for the reconstruction and prediction. In the computations we present here, we use the training sets that consists of $900$ snapshots and Hankel matrices of dimension $600 \times 300$ sliding along the computational domain with the chosen step. The length of the training sets was chosen so that at least two time periods of the global disturbance that we want to reveal are included in them.  

\begin{figure}[ht]
\centering
\begin{subfigure}{.35\textwidth}
	\includegraphics[width=\linewidth,height=1.5in]{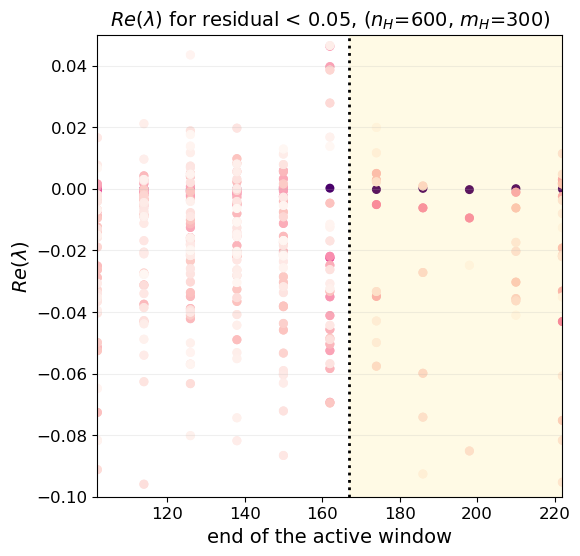}
	\end{subfigure}
\begin{subfigure}{.35\textwidth}	
	\includegraphics[width=\linewidth,height=1.5in]{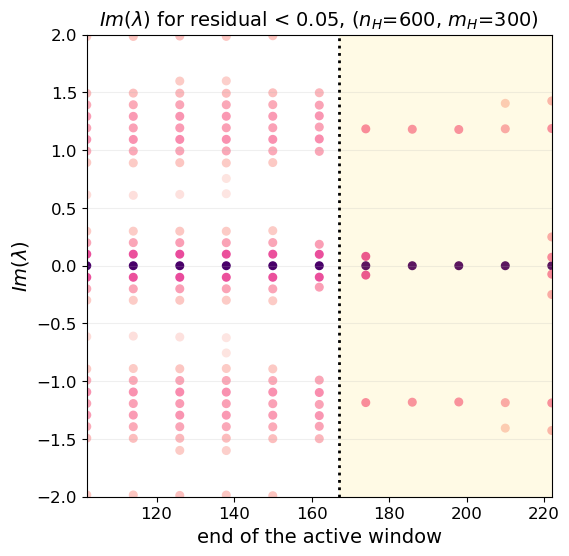}
	\end{subfigure}
\vspace{-1mm}
		\caption{Physiology model. The real and the imaginary parts of eigenvalues for sliding active windows for which the residuals are smaller than the threshold $\eta_r = 0.025$.
			The intensity of color of eigenvalues is associated with the amplitude of modes.}
		\label{fig:Physio1}
\end{figure}

\begin{figure}[ht]
	\centering
	\includegraphics[width=0.72\linewidth,height=1.5in]{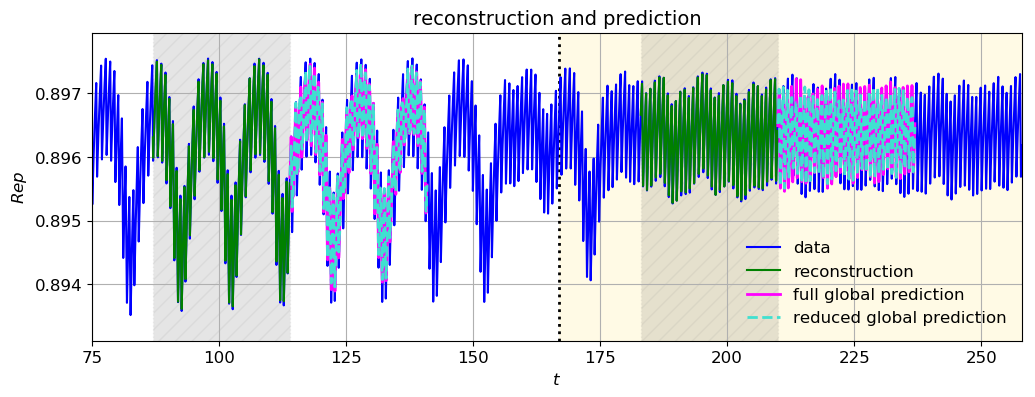}
	\includegraphics[width=0.25\linewidth,height=1.5in]{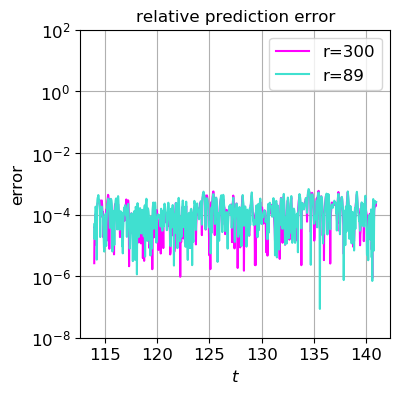}
	\hfill
	\includegraphics[width=0.72\linewidth,height=1.5in]{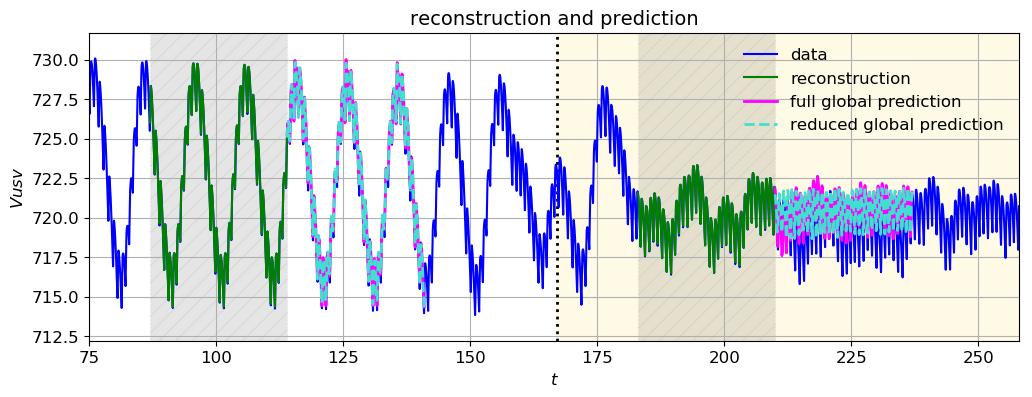}
	\includegraphics[width=0.25\linewidth,height=1.5in]{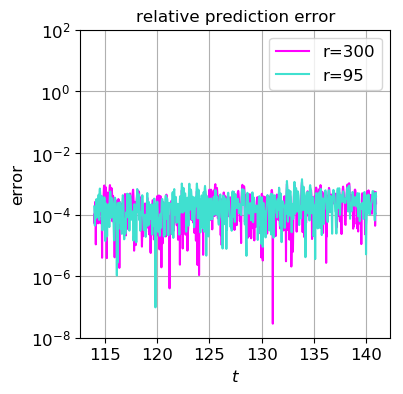}
	\caption{Physiology model. Extrasplanchnic peripheral resistance and splanchnic venous unstressed
		volume.  First column: Reconstruction and prediction obtained using \textsf{GKP} for the chosen active windows.  The full ($r=300$) and reduced ($r<300$) prediction obtained with \textsf{GKP} on active windows in the first part of simulation in the breathing zones captures accurately the dynamics, while even the full prediction obtained with the training set in the resonant breathing zone does not capture well the dynamics. Second column: The prediction errors for the full and reduced prediction (with $r$ modes) for the chosen active window in the first part of simulation.}
	\label{fig:Physio2}
\end{figure}

The switching moment from resonant to normal breathing is nicely detected by eigenvalues provided by \textsf{DMD} algorithm in Figure \ref{fig:Physio1}. The nature of eigenvalues changes significantly in the active window beginning at $t=150$ when switching of the dynamics occurs. It is nicely visible from Figures \ref{fig:Physio1} -- \ref{fig:Physio2} that when the training set is in the zone of the resonant breathing, the \textsf{GKP} results are perfect in that same zone, and then deteriorate as we move into the normal breathing zone. This is as expected since after the moment of transition from the resonant to normal breathing the dynamical system is not governed by the same set of parameters. When the training set is in the zone of the normal breathing, the \textsf{GKP} in that zone is much less accurate then in the resonant breathing zone. It catches well the global behavior but it is poor in the details.

\noindent What we can conclude from the presented results is the following. When the training set is in the zone of the resonant breathing, the \textsf{GKP} algorithm results are perfect in that same zone, and then deteriorate as we move into the normal breathing zone. This is as expected since after the moment of transition form the resonant to normal breathing the dynamical system is not governed by the same set of parameters. When the training set is in the zone of the normal breathing, the \textsf{GKP} in that same normal breathing zone is much less accurate then in the resonant breathing zone. It catches well the global behavior but it is poor in the detail, most of the error value comes from the difference in the phase. This is also logical since normal breathing is much more irregular than the resonant breathing and it turns out that it can not be learned with high accuracy.

\subsection{Case study: Geomagnetic substorms}\label{S=GeoMag}

\begin{figure}[ht]
	\centering
	\includegraphics[width=\linewidth,height=2in]{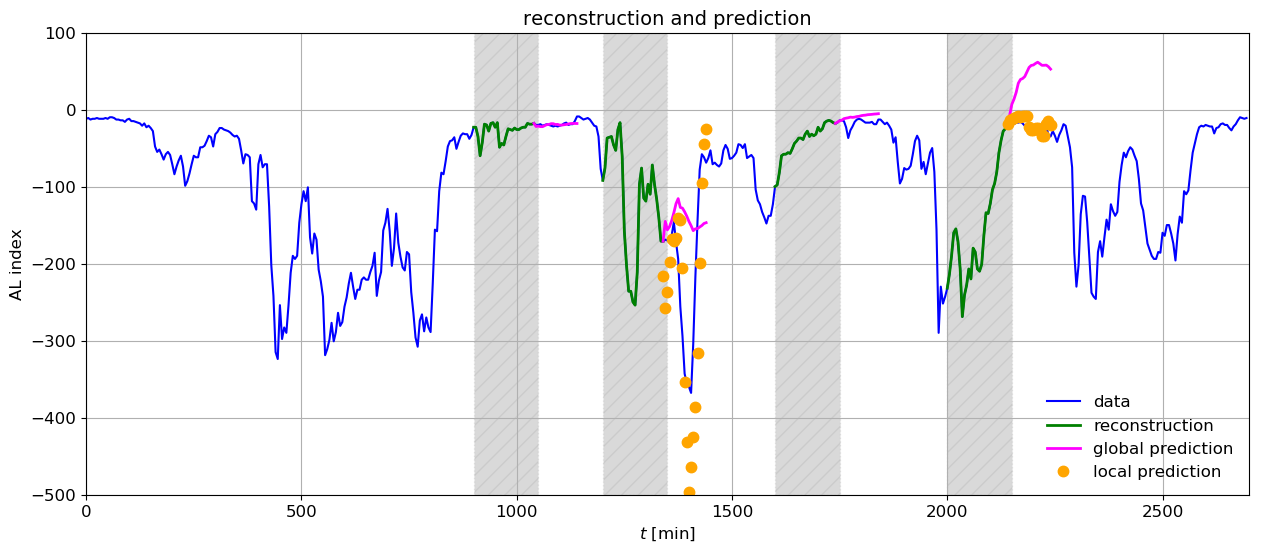}
	\vspace{-4mm}
	\caption{Geomagnetic substorms data:  \textsf{KMD} reconstruction and prediction of the \textsf{AL} index. {For reconstruction and global prediction, the data are collected in four active windows (time intervals $[900,1050]$, $[1200,1350]$, $[1600,1750]$, $[2000,2150]$, indicated by shadowed rectangles), the algorithm  uses $20 \times 10$ Hankel matrices and then the dynamics is predicted for $20$ time steps ahead.} {The time resolution of the collected data is $\delta t = 5$ min.} The local prediction algorithm uses  $3 \times 2$ matrices  and the error threshold for {resizing the Hankel matrix to minimal size (switching to  the local algorithm)} is set to $10$. {The dynamics with the local prediction algorithm is predicted two time steps ahead.}}
	\label{fig:geomag1}
\end{figure}

Geomagnetic storms and substorms are violent disturbances of the Earth's magnetosphere, caused by energy transfer of the solar wind into the planets magnetosphere, with potentially severe impact on the human civilization 
\cite{NOAA} \cite{Hapgood-2019}. \\
\indent Physics-based modeling (see e.g. \cite{spencer-geomag-substorm-2018}, \cite{Storm-substorm-book}, \cite{kamide1998current}) of geomagnetic substorms and storm/substorm interaction is a challenging task and the subject of intensive study. It  must cover multiscale, nonlinear interactions of plasmas that are not in equilibrium, or are in an unstable equilibrium, which makes such modeling difficult to apply when prediction is needed, see e.g. \cite{Morley-2007}. \\
\indent On the other hand, given an abundance of observation data, a data-driven approach is an attractive alternative; see e.g. \cite{giannakis-2016arXiv161207272G}, \cite{maimaiti-geomag-substorm-deep-learn-2019}. 
The intensity of a substorm is quantified by the Auroral Electrojet (\textsf{AE})  index, the \textsf{AL}, which is a measure of the magnitude of the geomagnetic field disturbances on the ground induced by ionospheric currents developed during substorm. Other information such as e.g. solar wind data \cite{NEWELL201628}, the \emph{Dst} index,  and other substorm signature indices may be available \cite{Nose-substorm-Wp-index-2012} and used as observables.\\
\indent For the purpose of this case study of the proposed approach as a purely data driven black-box methodology, we choose to use the \textsf{AL} index as the only observable; the data are downloaded  from the Kyoto Geomagnetism Data Service (\texttt{http://wdc.kugi.kyoto-u.ac.jp/}). \\
%
\indent The presented results are obtained by using global prediction algorithm with the active windows of size $30$ and the Hankel matrices of dimension $20 \times 10$. By sliding the active windows along the computational domain we get prediction at different time moments. In Figure \ref{fig:geomag1} we present the obtained reconstruction and prediction results for four active windows. Note that the modal representation of the signal is good, and it could provide a valuable insights to the experts in magnetic storm physics. \\
\indent In the framework of our theory, the poor global prediction results are to be expected -- almost all eigenvectors used in the \textsf{KMD} have large residuals. 
This once more justifies our approach, based on using the residuals of the Ritz pairs, computable even in the data driven setting, using the method from \cite{DDMD-RRR-2018}.
However, large prediction errors trigger the switch to the local prediction algorithm, which delivers more accurate predictions, at least for shorter lead time,  as shown in Figure \ref{fig:geomag1}.

\subsection{Additional numerical results for prediction of influenza cases}\label{SS=Influenza_more}
\begin{figure}[ht]
	\centering
	\includegraphics[width=.99\linewidth,height=1.0in]{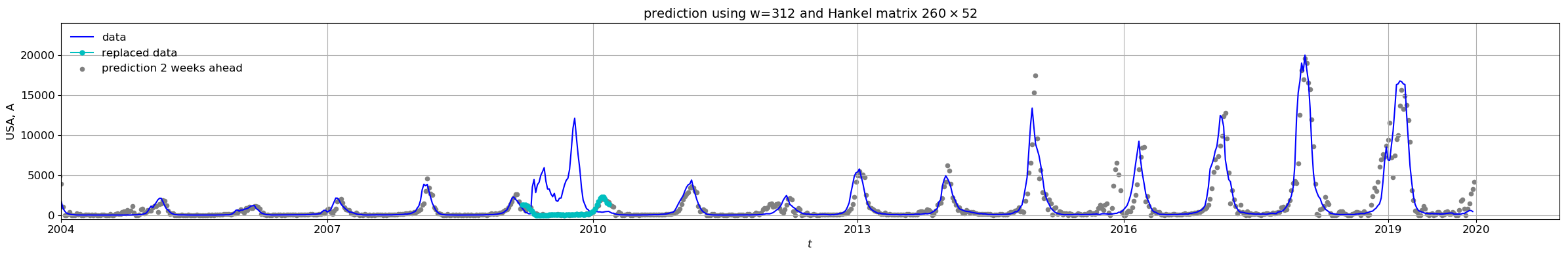}
	\includegraphics[width=.99\linewidth,height=1.0in]{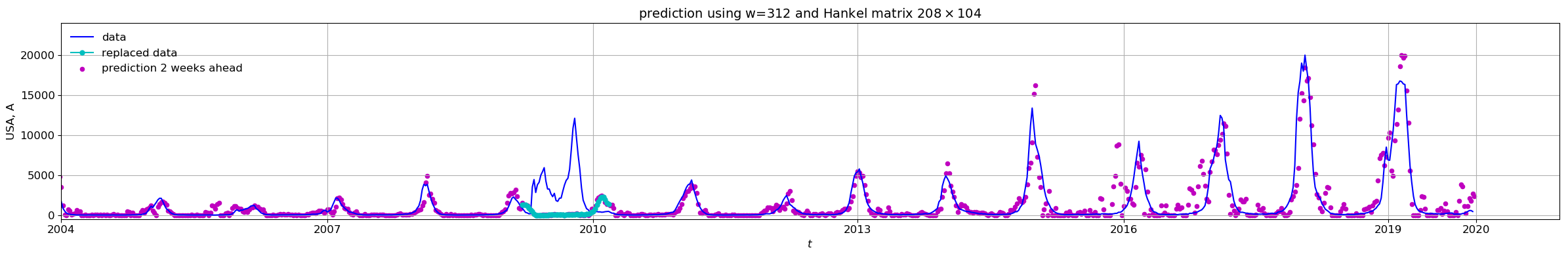}
	\includegraphics[width=.99\linewidth,height=1.0in]{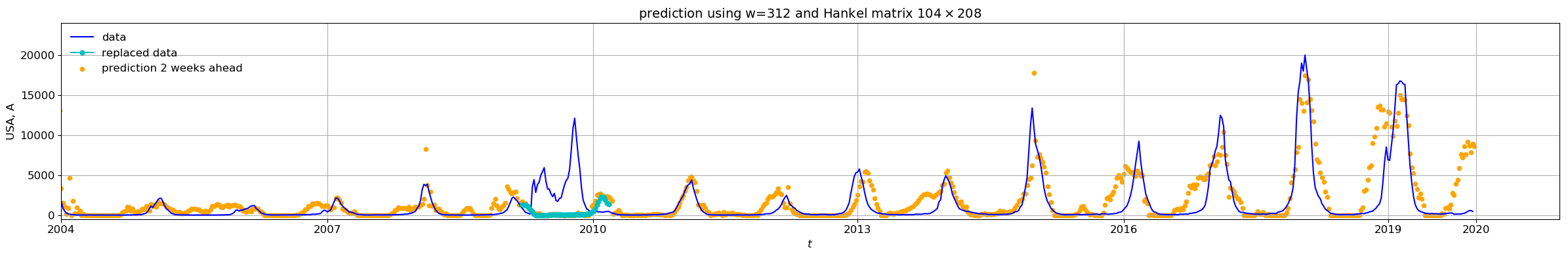}
	\vspace{-3mm}
	\caption{Influenza data (USA). Global Koopman prediction on influenza with the size of active windows 312 and different sizes of Hankel matrices. The prediction 2 weeks ahead by using KMD's from the active windows sliding along computational domain with sliding step $\Delta p=1$.}
	\label{fig:influenzaUSA}
\end{figure}

\begin{figure}[ht]
	\centering
	\includegraphics[width=.99\linewidth,height=1.0in]{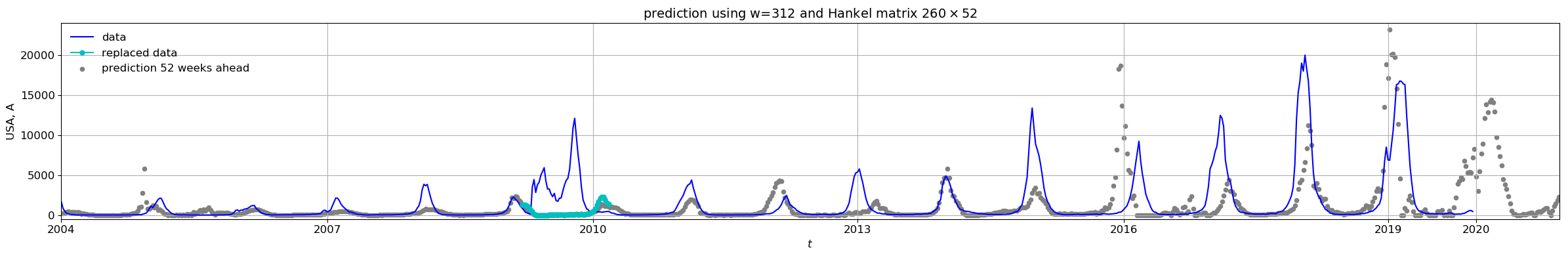}
	\includegraphics[width=.99\linewidth,height=1.0in]{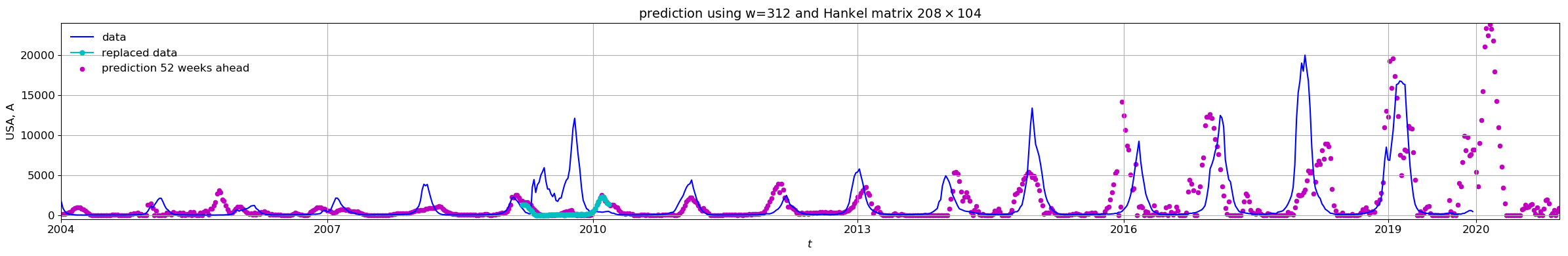}
	\includegraphics[width=.99\linewidth,height=1.0in]{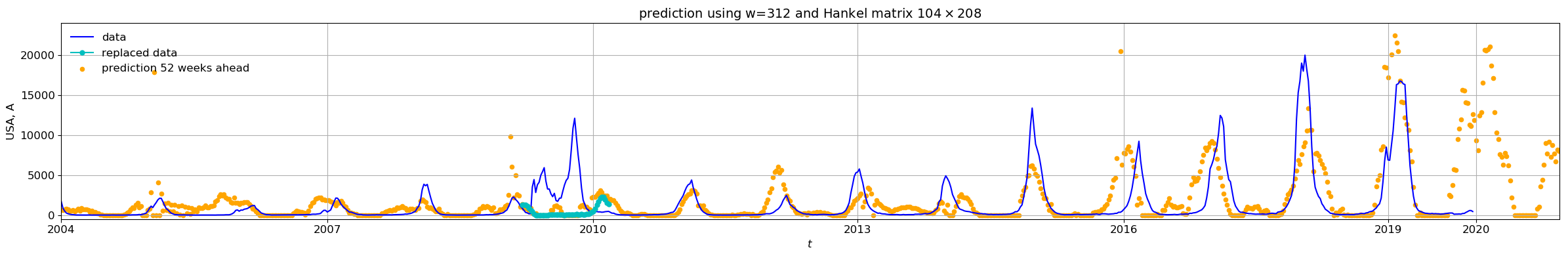}
	\vspace{-3mm}
	\caption{Influenza data (USA). Global Koopman prediction on influenza with the size of active windows 312 and different sizes of Hankel matrices. The prediction 52 weeks ahead by using KMD's from the active windows sliding along computational domain  with sliding step $\Delta p=1$.}
	\label{fig:influenzaUSA-2}
\end{figure}

\begin{figure}[ht]
	\centering
	\includegraphics[width=.99\linewidth,height=1.in]{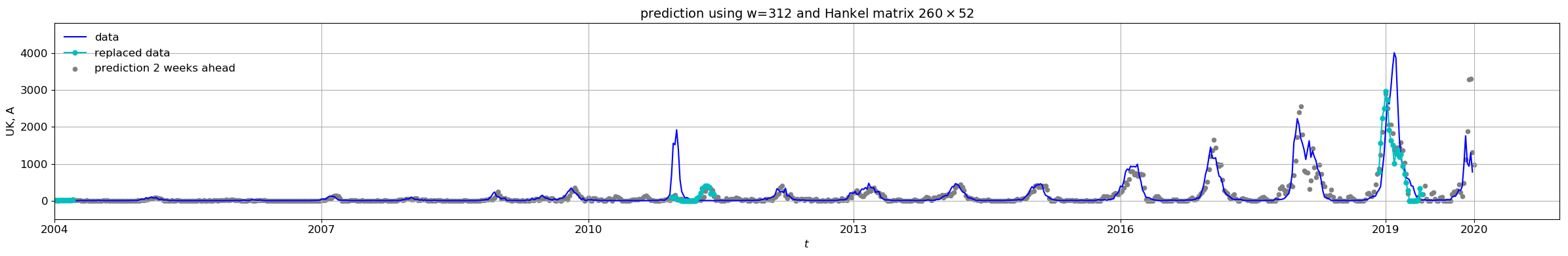}
	\includegraphics[width=.99\linewidth,height=1.in]{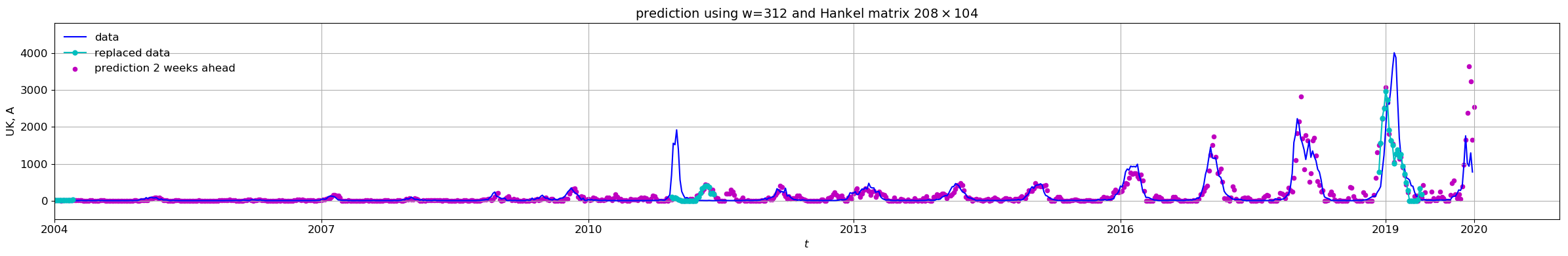}
	\includegraphics[width=.99\linewidth,height=1.in]{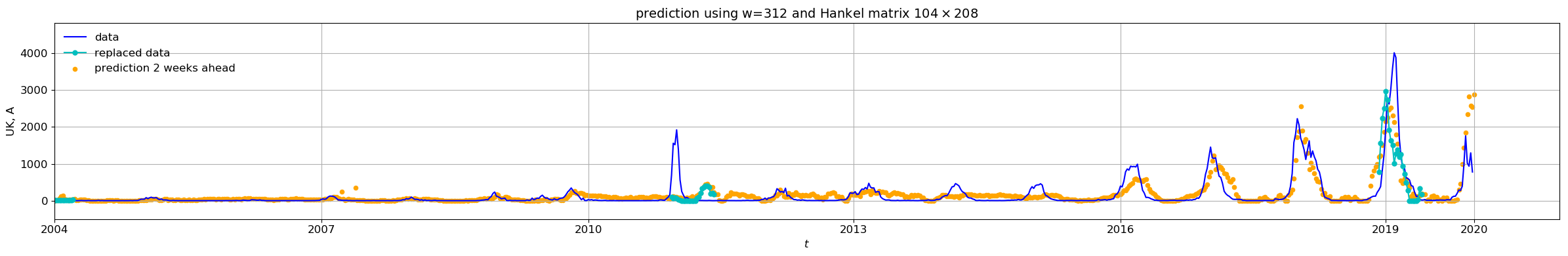}
	\vspace{-3mm}
	\caption{Influenza data (UK). Global Koopman prediction on influenza with the size of active windows 312 and different sizes of Hankel matrices. The prediction 2 weeks ahead by using KMD's from the active windows sliding along computational domain  with sliding step $\Delta p=1$.}
	\label{fig:influenzaUK}
\end{figure}

\begin{figure}[ht]
	\includegraphics[width=.99\linewidth,height=1.in]{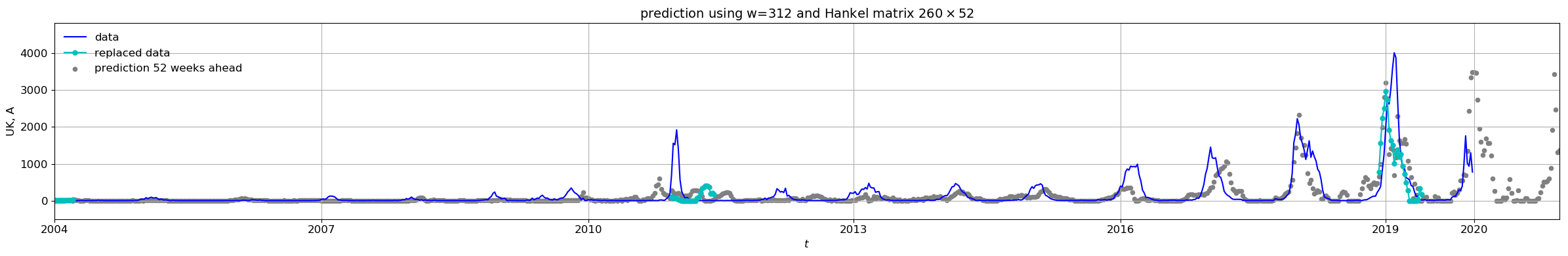}
	\includegraphics[width=.99\linewidth,height=1.in]{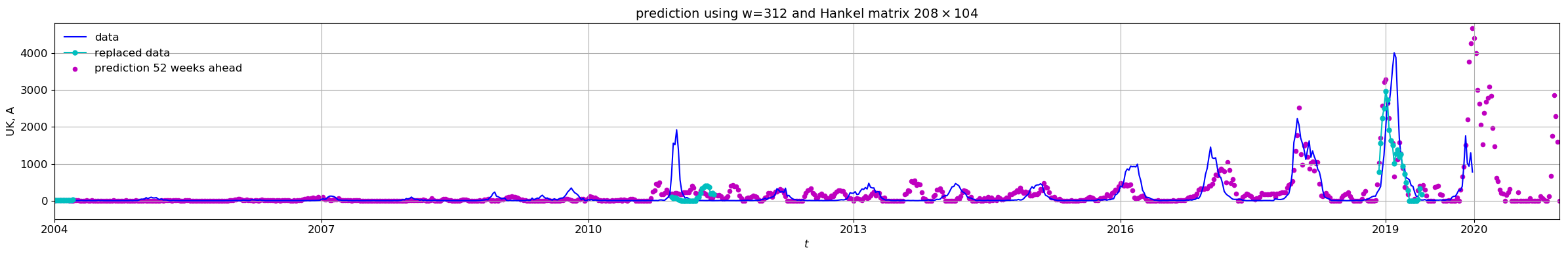}
	\includegraphics[width=.99\linewidth,height=1.in]{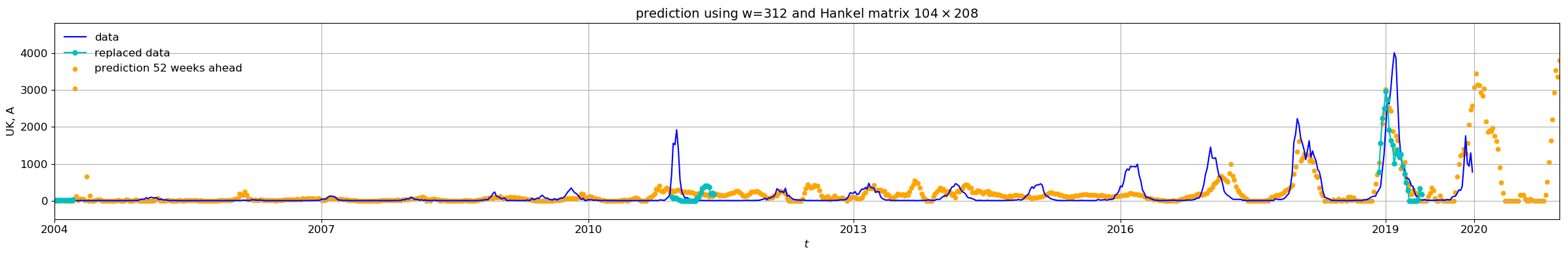}
	\vspace{-3mm}
	\caption{Influenza data (UK). Global Koopman prediction on influenza with the size of active windows 312 and different sizes of Hankel matrices. The prediction 52 weeks ahead by using KMD's from the active windows sliding along computational domain with sliding step $\Delta p=1$.}
	\label{fig:influenzaUK-2}
\end{figure}

Here we provide some additional results for the material in \S 2 of the paper.
In Figures \ref{fig:influenzaUSA}--\ref{fig:influenzaUK-2} we show prediction of the dynamics of influenza for USA and UK for $2$ and $52$ weeks ahead, obtained with the \textsf{KMD} decompositions in the global prediction algorithm, using sliding active windows of size 312.

\begin{figure}[ht]
\begin{subfigure}{0.38\textwidth}
\includegraphics[width=0.49\linewidth,height=1.4in]{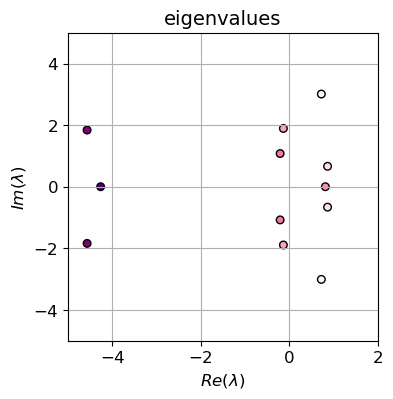}
\includegraphics[width=0.49\linewidth,height=1.4in]{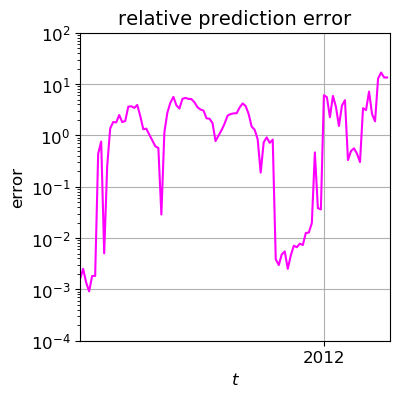}
\vspace{-6mm}
\caption{Before retouching}
\end{subfigure}
\begin{subfigure}{0.38\textwidth}
\includegraphics[width=0.49\linewidth,height=1.4in]{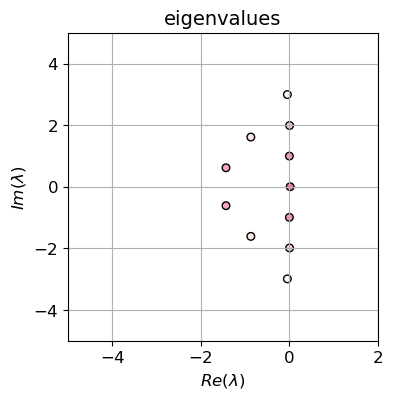}	
\includegraphics[width=0.49\linewidth,height=1.4in]{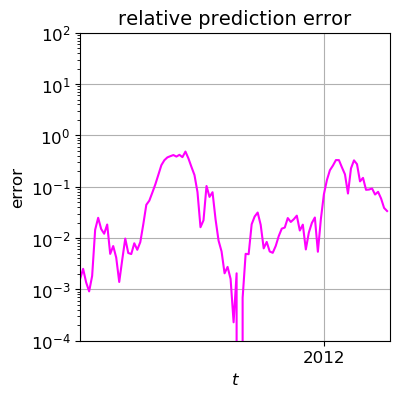}
\vspace{-6mm}
\caption{After retouching}
	\end{subfigure}
\begin{subfigure}{0.22\textwidth}
\includegraphics[width=\linewidth,height=1.4in]{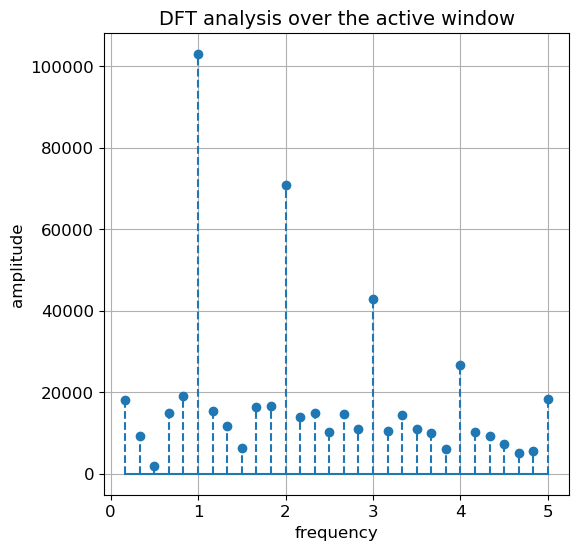}
\vspace{-6mm}
\caption{DFT analysis.}
	\end{subfigure}
\vspace{-2mm}
	\caption{Influenza data (USA). The most relevant eigenvalues and the prediction errors in the global algorithm (using $208 \times 104$ Hankel matrices) for the active window as in Figure 1b in the main paper. Note how the unstable eigenvalues ($\Re(\lambda)>0$) impact the  prediction performance, and how the retouching moves them to the left. Compare with Figures 1b and 2 in the main paper. Note how the dominant frequencies from the \textsf{DFT} analysis correspond to the imaginary parts of the eigenvalues computed after the retouching and selected by the residual criterion. }
	\label{fig:influenza-eigs-error-replace}
\end{figure}

In Figure \ref{fig:influenza-eigs-error-replace},  we provide additional numerical illustration (most relevant eigenvalues before and after retouching, with the corresponding prediction errors, and relation with the dominant frequencies from the \textsf{DFT} analysis) related to Figure 1b in the main paper.

\end{document}